\documentclass[a4paper]{article}%
\usepackage{amsmath}
\usepackage{amsfonts}
\usepackage{amssymb}
\usepackage{graphicx}%
\providecommand{\U}[1]{\protect\rule{.1in}{.1in}}
\newtheorem{theorem}{Theorem}

\newtheorem{corollary}[theorem]{Corollary}

\newtheorem{lemma}[theorem]{Lemma}

\newtheorem{proposition}[theorem]{Proposition}
\newtheorem{remark}[theorem]{Remark}

\newenvironment{proof}[1][Proof]{\noindent\textbf{#1.} }{\ \rule{0.5em}{0.5em}}
\begin{document}

\title{Construction of KdV flow\\-a unified approach-}
\author{Shinichi KOTANI\\Nanjing University, Osaka University }
\date{}
\maketitle

\begin{abstract}
A KdV flow is constructed on a space whose structure is described in terms of
the spectrum of the underlying Schr\"{o}dinger operators. The space includes
the conventional decaying functions and ergodic ones. Especially any smooth
almost periodic function can be initial data for the KdV equation.

\end{abstract}

\section{Introduction}

This article is a continuation of \cite{k2}, where a KdV flow was constructed
on a space of potentials with reflectionless property on an energy interval
$[\lambda_{1},\infty)$. Since the KdV equation is closely related with 1D
Schr\"{o}dinger operators, we use the terminology potentials to describe
initial data for the KdV equation. When the previous paper was written, the
author intended to remove this reflectionless property by approximating
general potentials by reflectionless potentials, which made the procedure
rather involved. However he has recognized that a direct extension is possible
independently of the last paper. Therefore the present paper is readable
without \cite{k2}, although its knowledge would be very helpful for prompt
understanding of the whole context.

Our approach to this problem is essentially based on Sato's philosophy
\cite{sa}, whose analytical version was given by Segal-Wilson \cite{s-w}. From
our point of view that is an analysis on eigen-spaces of underlying
Schr\"{o}dinger operators which seems quite natural due to GGKM and Lax.

To give perspective and state the main results several terminologies and
notations have to be prepared. For positive odd integer $n$ let $\Gamma_{n}$
be%
\[
\Gamma_{n}=\left\{  g=e^{h}\text{; \ }h\text{ is a real odd polynomial of
degree}\leq n\right\}
\]
and $C$ be a simple smooth closed curve in $\mathbb{C}\cup\left\{
\infty\right\}  $ defined by%
\[
C=\left\{  \pm\omega\left(  y\right)  +iy\text{; \ }y\in\mathbb{R}\right\}
\]
with a smooth positive function $\omega$ on $\mathbb{R}$ satisfying
$\omega\left(  y\right)  =\omega\left(  -y\right)  $, hence $C$ satisfies%
\[
C=-C\text{, \ \ \ }C=\overline{C}\text{.}%
\]
$D_{\pm}$ are the interior and exterior domains separated by the curve $C$
defined by%
\[
D_{+}=\left\{  z\in\mathbb{C}\text{; }\left\vert \operatorname{Re}z\right\vert
<\omega\left(  \operatorname{Im}z\right)  \right\}  \text{,\ }D_{-}=\left\{
z\in\mathbb{C}\text{; \ }\left\vert \operatorname{Re}z\right\vert
>\omega\left(  \operatorname{Im}z\right)  \right\}  \text{.}%
\]
The curve $C$ is chosen so that $g\in\Gamma_{n}$ remains bounded on $D_{+}$,
or more concretely%
\[
\omega\left(  y\right)  =O\left(  y^{-\left(  n-1\right)  }\right)  \text{
\ as }\left\vert y\right\vert \rightarrow\infty\text{.}%
\]
The \textbf{Hardy spaces} associated with curve $C$ is defined by%
\[
\left\{
\begin{array}
[c]{l}%
H\left(  D_{+}\right)  =\text{the closure in }L^{2}\left(  C\right)  \text{ of
rational functions with no poles in }D_{+}\\
H\left(  D_{-}\right)  =\text{the closure in }L^{2}\left(  C\right)  \text{ of
rational functions with no poles in }D_{-}%
\end{array}
\right.  \text{.}%
\]
It is known that%
\[
L^{2}\left(  C\right)  =H\left(  D_{+}\right)  \oplus H\left(  D_{-}\right)
\text{ \ (not necessarily orthogonal),}%
\]
and elements of $H\left(  D_{\pm}\right)  $ can be extended as analytic
functions on $D_{\pm}$ respectively. The projections to $H\left(  D_{\pm
}\right)  $ are given by%
\[
\left\{
\begin{array}
[c]{l}%
\mathfrak{p}_{+}u\left(  z\right)  =\dfrac{1}{2\pi i}%
{\displaystyle\int_{C}}
\dfrac{u\left(  \lambda\right)  }{\lambda-z}d\lambda\text{ \ for }z\in
D_{+}\smallskip\\
\mathfrak{p}_{-}u\left(  z\right)  =\dfrac{1}{2\pi i}%
{\displaystyle\int_{C}}
\dfrac{u\left(  \lambda\right)  }{z-\lambda}d\lambda\text{ \ for }z\in D_{-}%
\end{array}
\right.  \text{ if }u\in L^{2}\left(  C\right)  \text{.}%
\]
We enlarge the space $H\left(  D_{+}\right)  $ to admit polynomials. Namely
for $N\in\mathbb{Z}_{+}$ set%
\[
H_{N}\left(  D_{+}\right)  =\left(  z-b\right)  ^{N}H\left(  D_{+}\right)
\]
with $b\in D_{-}$, and define a norm in $H_{N}\left(  D_{+}\right)  $ by%
\[
\left\Vert u\right\Vert _{N}=\sqrt{\int_{C}\left\vert u\left(  \lambda\right)
\right\vert ^{2}\left\vert \lambda\right\vert ^{-2N}\left\vert d\lambda
\right\vert }\text{.}%
\]
Clearly $H_{N}\left(  D_{+}\right)  $ does not depend on the choice of $b$,
and $z^{m}\in H_{N}\left(  D_{+}\right)  $ if $m\leq N-1$.

In the previous paper we constructed the KdV flow as an action of $\Gamma_{n}$
on a Grassmann manifold consisting of $z^{2}$-invariant subspaces of
$L^{2}\left(  \left\vert z\right\vert =r\right)  $. In the present case we
construct an extension of the flow not on a Grassmann manifold of subspaces of
$z^{N}L^{2}\left(  C\right)  $ but on a space of vector functions
$\boldsymbol{a}\left(  \lambda\right)  =\left(  a_{1}\left(  \lambda\right)
,a_{2}\left(  \lambda\right)  \right)  $ on $C$. An analogue of a $z^{2}%
$-invariant subspace is%
\[
W_{\boldsymbol{a}}=\left\{  \boldsymbol{a}\left(  \lambda\right)  u\left(
\lambda\right)  \text{; \ }u\in H_{N}\left(  D_{+}\right)  \right\}  \text{,}%
\]
where%
\[
\left\{
\begin{array}
[c]{l}%
u_{e}\left(  \lambda\right)  =\dfrac{1}{2}\left(  u\left(  \lambda\right)
+u\left(  -\lambda\right)  \right)  \text{, \ }u_{o}\left(  \lambda\right)
=\dfrac{1}{2}\left(  u\left(  \lambda\right)  -u\left(  -\lambda\right)
\right)  \smallskip\\
\boldsymbol{a}\left(  \lambda\right)  u\left(  \lambda\right)  =a_{1}\left(
\lambda\right)  u_{e}\left(  \lambda\right)  +a_{2}\left(  \lambda\right)
u_{o}\left(  \lambda\right)
\end{array}
\right.  \text{.}%
\]
In the present paper, however, spaces $W_{\boldsymbol{a}}$ will not appear explicitly.

For $L\in\mathbb{Z}_{+}$ a space of symbols of Toeplitz operators is
introduced:%
\[
A_{L}\left(  C\right)  =\left\{
\begin{array}
[c]{c}%
a\text{; }a\left(  \lambda\right)  \text{ is bounded on }C\text{ and there
exists a bounded analytic}\\
\text{function }f\text{ on }D_{+}\text{ such that }\lambda^{L}\left(  a\left(
\lambda\right)  -f\left(  \lambda\right)  \right)  \text{ is bounded on }C
\end{array}
\right\}  \text{.}%
\]
The number $L$ is related to the degree of differentiability of the flow. The
\textbf{Toeplitz operator} with symbol $a$ is defined by%
\[
\left(  T\left(  a\right)  u\right)  \left(  z\right)  =f\left(  z\right)
u\left(  z\right)  +\left(  \mathfrak{p}_{+}\left(  a-f\right)  u\right)
\left(  z\right)
\]
for $u\in H_{N}\left(  D_{+}\right)  $, which is possible if $L\geq N$. This
$T(a)$ does not depend on the choice of $f$ and defines a bounded operator on
$H_{N}\left(  D_{+}\right)  $. We have to treat vector symbols $\boldsymbol{a}%
\left(  \lambda\right)  $ and the vector version $\boldsymbol{A}_{L}\left(
C\right)  $ of $A_{L}\left(  C\right)  $ essentially due to the fact that the
underlying Schr\"{o}dinger operators are second order. The associated Toeplitz
operator is defined by%
\[
\left(  T\left(  \boldsymbol{a}\right)  u\right)  \left(  z\right)  =\left(
T\left(  a_{1}\right)  u_{e}\right)  \left(  z\right)  +\left(  T\left(
a_{2}\right)  u_{o}\right)  \left(  z\right)  \text{.}%
\]
Let%
\[
\boldsymbol{A}_{L}^{inv}\left(  C\right)  =\left\{  \boldsymbol{a}%
\in\boldsymbol{A}_{L}\left(  C\right)  \text{; \ }T\left(  \boldsymbol{a}%
\right)  \text{ is invertible on }H_{L}\left(  D_{+}\right)  \right\}
\text{.}%
\]
Since $1\in H_{1}\left(  D_{+}\right)  $, $z\in H_{2}\left(  D_{+}\right)  $,
one can define%
\[
u=T\left(  \boldsymbol{a}\right)  ^{-1}1\in H_{1}\left(  D_{+}\right)  \text{,
\ }v=T\left(  \boldsymbol{a}\right)  ^{-1}z\in H_{2}\left(  D_{+}\right)
\text{,}%
\]
if $\boldsymbol{a}\in\boldsymbol{A}_{L}^{inv}\left(  C\right)  $ and $L\geq
2$.$\boldsymbol{\ }$Set%
\[
\varphi_{\boldsymbol{a}}\left(  z\right)  =\boldsymbol{a}\left(  z\right)
u\left(  z\right)  -1\text{, \ }\psi_{\boldsymbol{a}}\left(  z\right)
=\boldsymbol{a}\left(  z\right)  v\left(  z\right)  -z\in H\left(
D_{-}\right)  \text{.}%
\]
Then there exist a constant $\kappa_{1}\left(  \boldsymbol{a}\right)  $ and
$\phi_{\boldsymbol{a}}\in H\left(  D_{-}\right)  $ such that%
\[
\varphi_{\boldsymbol{a}}\left(  z\right)  =\kappa_{1}\left(  \boldsymbol{a}%
\right)  z^{-1}+z^{-1}\phi_{\boldsymbol{a}}\left(  z\right)  \text{.}%
\]
We call the functions $\left\{  \varphi_{\boldsymbol{a}},\psi_{\boldsymbol{a}%
}\right\}  $ as \textbf{characteristic functions} for $\boldsymbol{a}%
\in\boldsymbol{A}_{L}^{inv}\left(  C\right)  $, since $\boldsymbol{a}$ is
uniquely determined by them. Define%
\[
m_{\boldsymbol{a}}\left(  z\right)  =\dfrac{z+\psi_{\boldsymbol{a}}\left(
z\right)  }{1+\varphi_{\boldsymbol{a}}\left(  z\right)  }+\kappa_{1}\left(
\boldsymbol{a}\right)  \text{ }\left(  =z+o(1)\right)  \text{.}%
\]

$\Gamma_{n}$ naturally acts on $\boldsymbol{A}_{L}\left(  C\right)  $, but not
always on $\boldsymbol{A}_{L}^{inv}\left(  C\right)  $. Schr\"{o}dinger
operators and the KdV equation are obtained by applying the group $\Gamma_{n}$
to $\boldsymbol{A}_{L}^{inv}\left(  C\right)  $. Let $e_{x}(z)=e^{xz}\in
\Gamma_{1}$ and suppose $e_{x}\boldsymbol{a}\in\boldsymbol{A}_{L}^{inv}\left(
C\right)  $ for any $x\in\mathbb{R}$. Then%
\[
f_{\boldsymbol{a}}\left(  x,z\right)  =e^{-xz}\left(  1+\varphi_{e_{x}%
\boldsymbol{a}}\left(  z\right)  \right)
\]
satisfies a Schr\"{o}dinger equation%
\[
-\partial_{x}^{2}f_{\boldsymbol{a}}\left(  x,z\right)  +q\left(  x\right)
f_{\boldsymbol{a}}\left(  x,z\right)  =-z^{2}f_{\boldsymbol{a}}\left(
x,z\right)
\]
with $q(x)=-2\partial_{x}\kappa_{1}\left(  e_{x}\boldsymbol{a}\right)  $. One
can recover $m_{\boldsymbol{a}}\left(  z\right)  $ by%
\[
m_{\boldsymbol{a}}\left(  z\right)  =-\dfrac{\left.  \partial_{x}%
f_{\boldsymbol{a}}\left(  x,z\right)  \right\vert _{x=0}}{f_{\boldsymbol{a}%
}\left(  0,z\right)  }\text{.}%
\]
A solution to the KdV equation is obtained by another family of functions
$e_{t,x}(z)=e^{xz+tz^{3}}$ of $\Gamma_{3}$, namely%
\[
q\left(  t,x\right)  =-2\partial_{x}\kappa_{1}\left(  e_{t,x}\boldsymbol{a}%
\right)
\]
satisfies%
\begin{equation}
\partial_{t}q\left(  t,x\right)  =\dfrac{1}{4}\partial_{x}^{3}q\left(
t,x\right)  -\dfrac{3}{2}q\left(  t,x\right)  \partial_{x}q\left(  t,x\right)
\text{ \ (KdV equation).} \label{e21}%
\end{equation}
Solutions to the higher order KdV equations can be obtained similarly. This is
the core of Sato's theory.

The basic quantity $m_{\boldsymbol{a}}$ is closely related to the Weyl
functions of Schr\"{o}dinger operators. If $q$ takes real values, one can
associate a Schr\"{o}dinger operator%
\[
L_{q}=-\partial_{x}^{2}+q
\]
with potential $q$. Throughout the paper we assume%
\begin{equation}
\left(  L_{q}u,u\right)  _{L^{2}\left(  \mathbb{R}\right)  }\geq\lambda
_{0}\left(  u,u\right)  _{L^{2}\left(  \mathbb{R}\right)  }\text{ \ for any
}u\in C_{0}^{\infty}\left(  \mathbb{R}\right)  \label{e18}%
\end{equation}
with some $\lambda_{0}<0$. Under this condition it is known that $L_{q}$ has a
unique self-adjoint extension, and there exist non-trivial functions $f_{\pm
}\left(  x,z\right)  \in L^{2}\left(  \mathbb{R}_{\pm}\right)  $ satisfying%
\[
-\partial_{x}^{2}f_{\pm}+qf_{\pm}=-z^{2}f_{\pm}\text{.}%
\]
These functions are unique up to constant multiple. The Weyl functions
$m_{\pm}$ are defined by%
\[
m_{+}(z)=\dfrac{\left.  \partial_{x}f_{+}\left(  x,z\right)  \right\vert
_{x=0}}{f_{+}\left(  0,z\right)  }\text{, \ \ }m_{-}(z)=-\dfrac{\left.
\partial_{x}f_{-}\left(  x,z\right)  \right\vert _{x=0}}{f_{-}\left(
0,z\right)  }\text{.}%
\]
$m_{\pm}$ are analytic functions on $\mathbb{C}\backslash\lbrack\lambda
_{0},\infty)$ and satisfy%
\[
\dfrac{\operatorname{Im}m_{\pm}(z)}{\operatorname{Im}z}>0\text{.}%
\]
Such an analytic function on $\mathbb{C}_{+}$ is called a \textbf{Herglotz
function}. The functions $m_{\pm}$ contain every information of the spectral
properties of $L_{q}$. The simplest one is the coincidence of the
\textrm{sp}$L_{q}$ with the domain of analyticity of $m_{\pm}$, hence $m_{\pm
}$ are analytic on $\mathbb{C}\backslash\lbrack\lambda_{0},\infty)$ and the
interior domain $D_{+}$ for the curve $C$ is supposed to contain the interval
$\left[  -\mu_{0},\mu_{0}\right]  $ with $\mu_{0}=\sqrt{-\lambda_{0}}$. One
thing which should be stressed here is that $m_{\pm}$ can be defined for any
potential $q$ regardless of decaying or oscillating. Moreover, since
$f_{\boldsymbol{a}}\left(  x,z\right)  \in L^{2}\left(  \mathbb{R}_{\pm
}\right)  $ holds depending on $\operatorname{Re}z\gtrless0$ under a certain
condition on $\boldsymbol{a}$, one see that $m_{\boldsymbol{a}}$ coincides
with the Weyl functions $m_{\pm}$, that is,%
\begin{equation}
m_{\boldsymbol{a}}\left(  z\right)  =\left\{
\begin{array}
[c]{cc}%
-m_{+}\left(  -z^{2}\right)  & \text{if \ }\operatorname{Re}z>0\\
m_{-}\left(  -z^{2}\right)  & \text{if \ }\operatorname{Re}z<0
\end{array}
\right.  \text{.} \label{e20}%
\end{equation}
Hence in this case $q$ is determined by $m_{\boldsymbol{a}}$ owing to the
inverse spectral theory. We call $m_{\boldsymbol{a}}$ as $m$\textbf{-function}
of $\boldsymbol{a}$, which will be the fundamental object in this paper, and
call $f_{\boldsymbol{a}}\left(  x,z\right)  $ as \textbf{Baker-Akhiezer
function} for $L_{q}$. For potentials $q$ decaying sufficiently fast
$f_{\boldsymbol{a}}\left(  x,z\right)  $ coincides with the Jost solution. It
should be mentioned that R. Johnson \cite{j} was the first who introduced the
Weyl functions to Sato's theory.

As we have observed above the invertibility of $T\left(  g\boldsymbol{a}%
\right)  $ is crucial, which is verified with the aid of
\textbf{tau-functions} in this paper. The tau-function was first introduced by
Hirota and its mathematical meaning was found by Sato. In our context it is
defined as the Fredholm determinant of the operator%
\[
g^{-1}T\left(  g\boldsymbol{a}\right)  T\left(  \boldsymbol{a}\right)
^{-1}:H_{N}\left(  D_{+}\right)  \rightarrow H_{N}\left(  D_{+}\right)
\text{,}%
\]
that is%
\[
\tau_{\boldsymbol{a}}\left(  g\right)  =\det\left(  g^{-1}T\left(
g\boldsymbol{a}\right)  T\left(  \boldsymbol{a}\right)  ^{-1}\right)  \text{.}%
\]
However to avoid a technical difficulty one version%
\[
\tau_{\boldsymbol{a}}^{\left(  2\right)  }\left(  g\right)  =\det
\nolimits_{2}\left(  g^{-1}T\left(  g\boldsymbol{a}\right)  T\left(
\boldsymbol{a}\right)  ^{-1}\right)
\]
is employed, whose definition is possible when the operator $g^{-1}T\left(
g\boldsymbol{a}\right)  T\left(  \boldsymbol{a}\right)  ^{-1}-I$ is of
Hilbert-Schmidt. The invertibility of $T\left(  g\boldsymbol{a}\right)  $ is
equivalent to $\tau_{\boldsymbol{a}}^{\left(  2\right)  }\left(  g\right)
\neq0$. Any $g\in\Gamma_{n}$ can be approximated by rational functions $r$
with the same number of zeros and poles in $D_{-}$. For such an $r$ the image
of $r^{-1}T\left(  r\boldsymbol{a}\right)  T\left(  \boldsymbol{a}\right)
^{-1}$ is finite dimensional and $m_{r\boldsymbol{a}}$, $\tau_{\boldsymbol{a}%
}\left(  r\right)  $ are computable by $\left\{  \varphi_{\boldsymbol{a}%
}\text{, }m_{\boldsymbol{a}}\right\}  $. Another key observation is%
\[%
\begin{tabular}
[c]{lll}%
$\tau_{\boldsymbol{a}}^{\left(  2\right)  }\left(  g\right)  \neq0\text{ for
any }g\in\Gamma_{n}$ & $\Longleftrightarrow$ & $\tau_{\boldsymbol{a}}\left(
r\right)  \geq0\text{ for any real rational functions}$\\
& $\Longleftrightarrow$ & $\dfrac{\operatorname{Im}m_{\boldsymbol{a}}\left(
z\right)  }{\operatorname{Im}z}>0$%
\end{tabular}
\ \ \ \ \ \text{,}%
\]
if $\boldsymbol{a}\in\boldsymbol{A}_{L}^{inv}\left(  C\right)  $ satisfies
$\boldsymbol{a}\left(  \lambda\right)  =\overline{\boldsymbol{a}\left(
\overline{\lambda}\right)  }$ on $C$, which yields $g\boldsymbol{a}%
\in\boldsymbol{A}_{L}^{inv}\left(  C\right)  $ for such an $\boldsymbol{a}%
\in\boldsymbol{A}_{L}^{inv}\left(  C\right)  $. Keeping these facts in mind we
define%
\[
\boldsymbol{A}_{L,+}^{inv}\left(  C\right)  =\left\{
\begin{array}
[c]{c}%
\boldsymbol{a}\in\boldsymbol{A}_{L}^{inv}\left(  C\right)  \text{;
}\boldsymbol{a}\left(  \lambda\right)  =\overline{\boldsymbol{a}\left(
\overline{\lambda}\right)  }\text{ on }C\text{,\ }\tau_{\boldsymbol{a}}\left(
r\right)  \geq0\text{ for real rational}\\
\text{function }r\text{ with the same number of zeros and poles in }D_{-}%
\end{array}
\right\}  \text{.}%
\]
One can obtain concrete elements of $\boldsymbol{A}_{L,+}^{inv}\left(
C\right)  $ by defining $\boldsymbol{a}$ directly from $m_{\pm}$. For a given
potential $q$ assume (\ref{e18}) and define $m$ by (\ref{e20}). Then $m$ is
analytic on $\mathbb{C}\backslash\left(  \left[  -\mu_{0},\mu_{0}\right]  \cup
i\mathbb{R}\right)  $ ($\mu_{0}=\sqrt{-\lambda_{0}}$) and satisfies\medskip
\newline(M.1) \ $m(z)=\overline{m(\overline{z})}$\ and%
\[
\left\{
\begin{array}
[c]{lcl}%
\dfrac{\operatorname{Im}m\left(  z\right)  }{\operatorname{Im}z}>0 & \text{on}
& \mathbb{C}\backslash\left(  \mathbb{R}\cup i\mathbb{R}\right)  \smallskip\\
\dfrac{m(x)-m(-x)}{x}>0 & \text{if} & x\in\mathbb{R}\text{ and\ }\left\vert
x\right\vert >\mu_{0}%
\end{array}
\right.  \text{.}%
\]
Assume further\newline(M.2) \ $m$ has an asymptotic behavior:%
\[
m\left(  z\right)  =z+\sum_{1\leq k\leq L-2}m_{k}z^{-k}+O\left(
z^{-L+1}\right)  \text{ \ on }D_{-}\text{.}%
\]
\newline Then one has

\begin{theorem}
\label{t1}If $m$ satisfies (M.1), (M.2) for $L\geq2$, then $\boldsymbol{m}%
\left(  z\right)  \equiv\left(  1,m(z)/z\right)  \in\boldsymbol{A}_{L,+}%
^{inv}\left(  C\right)  $ and the $m$-function $m_{\boldsymbol{m}}$ for
$\boldsymbol{m}$ is $m$.
\end{theorem}

If $q\in C^{L-2}\left(  -\delta,\delta\right)  $, then it is known that the
asymptotics of (M.2) holds in a sector%
\[
\left\vert \arg z\right\vert <\frac{\pi}{2}-\epsilon\text{, \ \ \ }\left\vert
\pi-\arg z\right\vert <\frac{\pi}{2}-\epsilon\text{.}%
\]
However the domain $D_{-}$ is wider even for $n=1$, and its boundary
approaches to the axis $i\mathbb{R}$ if $n\geq3$, therefore it is not trivial
to find $q$ satisfying (M.2) in $D_{-}$. Later in Theorems \ref{t3}, \ref{t4}
$m$ associated with the Weyl functions $m_{\pm}$ will be shown to fulfill
(M.2) if $q$ decays sufficiently fast or oscillates suitably.

Set%
\[
\mathcal{Q}_{L}\left(  C\right)  =\left\{  q\left(  x\right)  =-2\partial
_{x}\kappa_{1}\left(  e_{x}\boldsymbol{a}\right)  \text{; \ }\boldsymbol{a}%
\in\boldsymbol{A}_{L,+}^{inv}\left(  C\right)  \right\}  \text{.}%
\]
Then $m_{\boldsymbol{a}}$ is identified with $m_{\pm}$ of $q$ by (\ref{e20})
for $\boldsymbol{a}\in\boldsymbol{A}_{L,+}^{inv}\left(  C\right)  $, hence the
inverse spectral theory show that $m_{\boldsymbol{a}}$ determines $q$. This
makes it possible to define%
\[
\left(  K\left(  g\right)  q\right)  \left(  x\right)  =-2\partial_{x}%
\kappa_{1}\left(  e_{x}g\boldsymbol{a}\right)  \text{ with }q\left(  x\right)
=-2\partial_{x}\kappa_{1}\left(  e_{x}\boldsymbol{a}\right)
\]
for $\boldsymbol{a}\in\boldsymbol{A}_{L,+}^{inv}\left(  C\right)  $,
$g\in\Gamma_{n}$. One has

\begin{theorem}
\label{t2}Suppose $L\geq\max\left\{  n+1,3\right\}  $. Then $\left\{  K\left(
g\right)  \right\}  _{g\in\Gamma_{n}}$ defines a flow on $\mathcal{Q}%
_{L}\left(  C\right)  $. For a real odd polynomial $h$ of degree $n$ the
function $\left(  K\left(  e^{th}\right)  q\right)  \left(  x\right)  $ is
$C^{1}$ in $t$ and $C^{n}$ in $x$ and satisfies the $\left(  n+1\right)  /2$th
KdV equation. Especially $K\left(  e^{tz^{3}}\right)  q\left(  x\right)  $
satisfies the KdV equation%
\[
\partial_{t}q(t,x)=\dfrac{1}{4}\partial_{x}^{3}q(t,x)-\dfrac{3}{2}%
q(t,x)\partial_{x}q(t,x)
\]
if $q\in\mathcal{Q}_{L}\left(  C\right)  $ for $L\geq4$.
\end{theorem}

We summarize the procedure to obtain $K(g)q$ for a given $q$ as follows.
Define $m$ by (\ref{e20}) and assume the condition (M.2) for $m$. Solve the
equation in $H_{N}\left(  D_{+}\right)  $ for $z\in D_{+}$%
\[
1=e^{xz}g\left(  z\right)  \boldsymbol{f}\left(  z\right)  u\left(
x,z\right)  +\dfrac{1}{2\pi i}\int_{C}\dfrac{e^{x\lambda}g\left(
\lambda\right)  \left(  \boldsymbol{m}\left(  \lambda\right)  -\boldsymbol{f}%
\left(  \lambda\right)  \right)  u\left(  x,\lambda\right)  }{\lambda
-z}d\lambda\text{,}%
\]
where $\boldsymbol{m}\left(  z\right)  =(1,m(z)/z)$ and $\boldsymbol{f}\left(
z\right)  =\left(  1,f(z)\right)  $ with a bounded analytic function $f$ on
$D_{+}$ such that%
\[
m\left(  z\right)  -zf\left(  z\right)  =O\left(  z^{-L+1}\right)  \text{ on
}C\text{,}%
\]
which is possible due to (M.2). The Baker-Akhiezer function is obtained by
$f_{g\boldsymbol{m}}\left(  x,z\right)  =g\left(  z\right)  \boldsymbol{m}%
\left(  z\right)  u\left(  x,z\right)  $, and $\kappa_{1}\left(
e_{x}g\boldsymbol{m}\right)  $ is determined by%
\[
e^{xz}f_{g\boldsymbol{m}}\left(  x,z\right)  =1+\kappa_{1}\left(
e_{x}g\boldsymbol{m}\right)  z^{-1}+o\left(  z^{-1}\right)  \text{.}%
\]
Then we have $\left(  K(g)q\right)  \left(  x\right)  =-2\partial_{x}%
\kappa_{1}\left(  e_{x}g\boldsymbol{m}\right)  $. Especially if $g=1$, one can
recover $q$ from the Weyl functions $m_{\pm}$, which yields another way of the
inverse spectral problem.\medskip

Any concrete example of initial data for the KdV flow is provided by Theorem
\ref{t1}. For a given $m$ we have to verify the condition (M.2). There are two
classes of potentials satisfying (M.2).

If $q^{\left(  j\right)  }\in L^{1}\left(  \mathbb{R}\right)  $ for $j=0$,
$1$,$\cdots$, $L-2$, then (M.2) is valid in $\overline{\mathbb{C}}_{+}%
\equiv\left\{  z\in\mathbb{C}\text{; }\operatorname{Im}z\geq0\right\}  $ for
$L$ (see \cite{r1}), which will be shown in Proposition \ref{p10}.

The extended notion of reflection coefficients is defined by%
\[
R(z,q)=\frac{m_{+}\left(  z\right)  +\overline{m_{-}\left(  z\right)  }}%
{m_{+}\left(  z\right)  +m_{-}\left(  z\right)  }%
\]
where $m_{\pm}$ are the Weyl functions for $q$. The modulus $\left\vert
R(z,q)\right\vert $ coincides with that of the conventional reflection
coefficient on $\mathbb{R}$ if $q$ decays sufficiently fast.

\begin{theorem}
\label{t3}If $R(z,q)$ satisfies%
\begin{equation}
\int_{0}^{\infty}\lambda^{M}\left\vert R(\lambda,q)\right\vert d\lambda
<\infty\text{,} \label{108}%
\end{equation}
then $q\in\mathcal{Q}_{L}\left(  C\right)  $ with $L=M+2-\left(  n+1\right)
/2$ holds.
\end{theorem}

If $R(\lambda,q)=0$ for a.e.$\lambda\geq\lambda_{1}$ for some $\lambda_{1}%
\leq\mathbb{R}$ (which means $q$ is reflectionless on $\left(  \lambda
_{1},\infty\right)  $), the condition (\ref{108}) is satisfied for any
$M\geq1$. This case was already treated in \cite{k3}. The resulting potential
$q$ is known to be meromorphic on $\mathbb{C}$ and uniformly bounded on
$\mathbb{R}$ including all its derivatives.

Since $\left\vert R(\lambda,q)\right\vert $ is invariant under the flow
$K\left(  g\right)  $, that is%
\begin{equation}
\left\vert R(\lambda,q)\right\vert =\left\vert R(\lambda,K\left(  g\right)
q)\right\vert \text{ \ for a.e. }\lambda\in\mathbb{R}\text{,} \label{109}%
\end{equation}
the condition (\ref{108}) is supposed to play a significant role to
investigate the flow $K(g)$ in future. (\ref{109}) will be shown in a separate
paper by using transfer matrices of $K\left(  g\right)  $.

On the other hand \cite{k1} showed for ergodic potential $q_{\omega}\left(
x\right)  $%
\[
\Sigma_{ac}\left(  q_{\omega}\right)  =\left\{  \lambda\in\mathbb{R}\text{;
}\left\vert R(\lambda,q_{\omega})\right\vert =0\right\}  \text{.}%
\]
Therefore in this case (\ref{108}) is equivalent to%
\begin{equation}
\int_{\mathbb{R}_{+}\backslash\Sigma_{ac}\left(  q_{\omega}\right)  }%
\lambda^{M}d\lambda<\infty\text{.} \label{111}%
\end{equation}
In particular for periodic potentials $\Sigma_{ac}\left(  q_{\omega}\right)
=\Sigma\left(  q_{\omega}\right)  $ (the spectrum of $L_{q}$) is valid, hence
(\ref{111}) means that the total length of spectral gaps is small, which can
be estimated by the norms of derivatives of $q$.

For ergodic potentials the condition (\ref{111}) requires the existence of
rich ac spectrum, although it admits singular spectrum. This situation can be
improved by replacing (\ref{108}) by a similar condition on the curve
$\widehat{C}=\left\{  -z^{2}\text{; }z\in C\text{, }\operatorname{Re}%
z>0\right\}  $, which enables us to have

\begin{theorem}
\label{t4}Let $\left\{  q_{\omega}\left(  x\right)  =q\left(  \theta_{x}%
\omega\right)  \right\}  $ be an ergodic process on $\left(  \Omega
,\mathcal{F},P,\left\{  \theta_{x}\right\}  _{x\in\mathbb{R}}\right)  $.
Suppose $q_{\omega}\in C_{b}^{m}\left(  \mathbb{R}\right)  $. Then,
$q_{\omega}\in\mathcal{Q}_{L}\left(  C\right)  $ holds for a.e. $\omega
\in\Omega$ for $L\leq\left(  m-3\left(  n-1\right)  \right)  /6$. In this case
$\left(  K(g)q_{\omega}\right)  \left(  x\right)  =f_{g}\left(  \theta
_{x}\omega\right)  $ for $g\in\Gamma_{n}$ is valid with $f_{g}\left(
\omega\right)  =\left(  K(g)q_{\omega}\right)  \left(  0\right)  $.
\end{theorem}

Any almost periodic potentials can be considered as ergodic potentials and one
can apply Theorem \ref{t3} or Theorem \ref{t4} to have solutions starting from
almost periodic functions. Under the condition of Theorem \ref{t4} one has
$q_{\omega}\in\mathcal{Q}_{L}\left(  C\right)  $ for every $\omega$ not for
a.e. $\omega$.

Rybkin \cite{r3} obtained solutions to the KdV equation with step like initial
data, which is decaying on $\mathbb{R}_{+}$ and arbitrary on $\mathbb{R}_{-}$.
He employed Hirota's tau-function and the Hankel transform on $\mathbb{R}$,
which restricts the class of initial data to step like functions. In our
approach the decaying condition on $\mathbb{R}_{+}$ can be removed, since we
represent the solutions through information of the Weyl functions $m_{\pm}$ on
$\mathbb{C}_{+}$ not on $\mathbb{R}$. Our framework also admits step like
initial data. For instance if $q\in C^{\infty}\left(  \mathbb{R}\right)  $ is
almost periodic on one axis $\mathbb{R}_{+}$ or $\mathbb{R}_{-}$ and decaying
on the opposite axis, since such a potential can be easily verified to satisfy (M.2).

For almost periodic initial data there are several papers. Egorova \cite{e}
treated limit periodic initial data. Damanik-Goldstein \cite{d-g} and
Eichinger-VandenBoom-Yuditskii \cite{evy} considered almost periodic
potentials. Their approaches are different from ours and the associated
Schr\"{o}dinger operators must have only ac spectrum. Tsugawa \cite{tu}
obtained solutions starting from quasi-periodic initial data without assuming
pure ac spectrum, but he could not show the existence of global solutions in time.

One of the advantages of Sato's approach lies on the algebraic nature of the
group $\Gamma_{n}$ acting on the space of symbols. Especially the factor
$q_{\zeta}\left(  z\right)  =\left(  1-\zeta^{-1}z\right)  ^{-1}$ plays a role
of primes in number theory, which will be frequently used in the present
paper.\medskip

Throughout the paper the following notations will be employed:%
\[
\left\{
\begin{array}
[c]{l}%
\mathbb{R}=\text{the set of all real numbers}\\
\mathbb{C}=\text{the set of all complex numbers}\\
\mathbb{Z}=\text{the set of all integers}\\
\mathbb{R}_{+}=\left\{  x\in\mathbb{R}\text{, }x\geq0\right\}  \text{,
\ \ \ \ \ }\mathbb{R}_{-}=\left\{  x\in\mathbb{R}\text{, }x\leq0\right\} \\
\mathbb{C}_{+}=\left\{  z\in\mathbb{C}\text{, }\operatorname{Im}z>0\right\}
\text{, \ }\mathbb{C}_{-}=\left\{  z\in\mathbb{C}\text{, }\operatorname{Im}%
z<0\right\} \\
\mathbb{Z}_{+}=\left\{  n\in\mathbb{Z}\text{, }n\geq0\right\}  \text{,
\ \ \ \ \ \ }\mathbb{Z}_{-}=\left\{  n\in\mathbb{Z}\text{, }n\leq0\right\} \\
\overline{z}\text{ denotes the complex conjugate of }z\text{: }\overline
{x+iy}=x-iy
\end{array}
\right.  \text{.}%
\]

\section{Hardy spaces and Toeplitz operators}

In the previous paper \cite{k2} we employed Segal-Wilson's version of Sato's
theory, in which they constructed KdV flow on a Grassmann manifold in $H\equiv
L^{2}\left(  \left\vert z\right\vert =r\right)  $. In his theory the Fourier
space $H$ is used as symbols space of pseudo-differential operators and the
separation%
\[
H=\left(  L^{2}\text{closure of }\left\{  z^{k}\right\}  _{k\geq0}\right)
\oplus\left(  L^{2}\text{closure of }\left\{  z^{k}\right\}  _{k<0}\right)
\equiv H_{+}\oplus H_{-}%
\]
is essential since the $H_{+}$ component exhibits the part of differential
operators and he had to take out differential operators parts from
pseudo-differential operators. Therefore the projection $\mathfrak{p}_{+}$ to
the Hardy space $H_{+}$ plays an essential role.

Since in this framework only a special class of solutions meromorphic on
$\mathbb{C}$ is possible to treat, the circle $\left\vert z\right\vert =r$
should be replaced by a certain unbounded curve to have more general solutions.

\subsection{Hardy spaces and projections}

Let $C$ be a simple closed smooth curve passing $\infty$ in the Riemann sphere
$\mathbb{C\cup}\left\{  \infty\right\}  $ and oriented anti-clockwisely. We
assume $C=-C$. The curve $C$ separates $\mathbb{C}_{\infty}$ into two domains
$D_{\pm}$, where $D_{+}$ contains the origin $0$. The situation was
illustrated in ( f1 ). The curve $C$ is chosen so that $e^{h\left(  z\right)
}$ remains bounded on $C$ , where $h(z)$ is a given real polynomial of odd degree.

Set%
\[
H\left(  D_{\pm}\right)  =L^{2}\left(  C\right)  \text{-closure of }\left\{
\text{rational functions with no poles in }D_{\pm}\right\}  \text{.}%
\]
Then%
\[
L^{2}\left(  C\right)  =H\left(  D_{+}\right)  \oplus H\left(  D_{-}\right)
\]
holds. For $f\in L^{2}\left(  C\right)  $ define%
\[
\mathfrak{p}_{\pm}f(z)=\pm\dfrac{1}{2\pi i}%
{\displaystyle\int_{C}}
\dfrac{f\left(  \lambda\right)  }{\lambda-z}d\lambda\text{ \ \ for \ }z\in
D_{\pm}\text{,}%
\]
and%
\[
\left(  \Theta f\right)  (z)=\lim_{\epsilon\downarrow0}\dfrac{1}{2\pi i}%
{\displaystyle\int_{C\cap\left\{  \left\vert \lambda-z\right\vert
>\epsilon\right\}  }}
\dfrac{f\left(  \lambda\right)  }{\lambda-z}d\lambda\text{ \ for \ }z\in C.
\]
It is known that $\Theta$ is a bounded operator on $L^{2}\left(  C\right)  $
(see \cite{ca} and \cite{cjs}) and $\mathfrak{p}_{\pm}$ have a finite limit
a.e. when $z$ approaches to an element of $C$. They satisfy%
\[
\left\{
\begin{array}
[c]{l}%
\mathfrak{p}_{+}f(z)=f\left(  z\right)  /2+\left(  \Theta f\right)  (z)\\
\mathfrak{p}_{-}f(z)=f\left(  z\right)  /2-\left(  \Theta f\right)  (z)
\end{array}
\text{ \ for }z\in C\right.  \text{,}%
\]
and $\mathfrak{p}_{\pm}$ are projections from $L^{2}\left(  C\right)  $ onto
$H\left(  D_{\pm}\right)  $ respectively. It should be noted that
$\mathfrak{p}_{\pm}$ are generally not orthogonal projections. If $D_{+}$ is a
disc, $H\left(  D_{+}\right)  $ coincides with the conventional Hardy space.

To treat an analogue of $z^{2}$-action on $L^{2}\left(  C\right)  $ for a
function $u$ on $C$ we define the even part and the odd part of $u$ by%
\[
\left\{
\begin{array}
[c]{c}%
u_{e}\left(  z\right)  =\dfrac{1}{2}\left(  u\left(  z\right)  +u\left(
-z\right)  \right)  \smallskip\\
u_{o}\left(  z\right)  =\dfrac{1}{2}\left(  u\left(  z\right)  -u\left(
-z\right)  \right)
\end{array}
\right.  \text{.}%
\]
The numbers $\pm1$ in front of $z$ come from the solutions to $\omega^{2}=1$.
Any function on $C$ can be represented as $u=u_{e}+u_{o}$ and this yields an
orthogonal decomposition in $L^{2}\left(  C\right)  $. It should be noted also
that%
\begin{equation}
\mathfrak{p}_{+}:L_{e}^{2}\left(  C\right)  \rightarrow H\left(  D_{+}\right)
\cap L_{e}^{2}\left(  C\right)  \text{ and }\mathfrak{p}_{+}:L_{o}^{2}\left(
C\right)  \rightarrow H\left(  D_{+}\right)  \cap L_{o}^{2}\left(  C\right)
\text{,} \label{e1}%
\end{equation}
where $L_{e}^{2}\left(  C\right)  $, $L_{o}^{2}\left(  C\right)  $ denote the
even part and the odd part respectively.

\subsection{Toeplitz operators}

In the previous paper \cite{k2} we considered $z^{2}$-invariant subspaces of
$L^{2}\left(  C\right)  $ when $C$ is a disc with center $0$. If the curve $C$
is unbounded, in place of subspaces of $L^{2}\left(  C\right)  $ we consider a
family of bounded vector functions $\boldsymbol{a}\left(  z\right)  =\left(
a_{1}(z),a_{2}(z)\right)  $ on $C$. The subspace corresponding to
$\boldsymbol{a}\left(  z\right)  $ is%
\[
W_{\boldsymbol{a}}\equiv\left\{  \boldsymbol{a}\left(  z\right)  u\left(
z\right)  \text{; \ }u\in H\left(  D_{+}\right)  \right\}  \subset
L^{2}\left(  C\right)  \text{,}%
\]
where%
\begin{equation}
\boldsymbol{a}\left(  z\right)  u\left(  z\right)  \equiv a_{1}(z)u_{e}%
(z)+a_{2}(z)u_{o}(z)\text{.} \label{2}%
\end{equation}
This space will not appear explicitly in the sequel, but the Toeplitz operator
with symbol $\boldsymbol{a}\left(  z\right)  $ plays an essential role.

Set%
\[
A(C)=\left\{  a\left(  \lambda\right)  \text{; \ }\sup_{\lambda\in
C}\left\vert a(\lambda)\right\vert <\infty\right\}  \text{,}%
\]
and%
\[
\left(  T\left(  a\right)  u\right)  \left(  z\right)  =\mathfrak{p}%
_{+}\left(  au\right)  \left(  z\right)  \text{ \ for }a\in A(C)\text{.}%
\]
Then $T\left(  a\right)  $ defines a bounded operator on $L^{2}\left(
C\right)  $ and is called a Toeplitz operator with symbol $a$.

To investigate the differentiability of solutions to the KdV equation we have
to admit the multiplication operation by rational functions on $H\left(
D_{\pm}\right)  $. To realize such operations some modification of the spaces
$H\left(  D_{\pm}\right)  $ is necessary. For an $N\in\mathbb{Z}_{+}$ and
$b\in D_{-}$ set%
\begin{equation}
H_{N}\left(  D_{+}\right)  =\left(  z-b\right)  ^{N}H\left(  D_{+}\right)  .
\label{3}%
\end{equation}
Clearly $H_{N}\left(  D_{+}\right)  $ does not depend on $b$, and
\[
z^{k}\in H_{N}\left(  D_{+}\right)  \text{ \ for any integer \ }k\leq
N-1\text{.}%
\]
For $u\in z^{N}L^{2}\left(  C\right)  $ we extend the definition of the
projections by%
\begin{equation}
\left\{
\begin{array}
[c]{c}%
\left(  \mathfrak{p}_{+}u\right)  \left(  z\right)  =\lim\limits_{D_{-}\ni
b^{\prime},\operatorname{Re}b^{\prime}\rightarrow\infty}\dfrac{\left(
z-b^{\prime}\right)  ^{N}}{2\pi i}%
{\displaystyle\int_{C}}
\dfrac{u\left(  \lambda\right)  }{\lambda-z}\left(  \lambda-b^{\prime}\right)
^{-N}d\lambda\text{ for }z\in D_{+}\medskip\\
\left(  \mathfrak{p}_{-}u\right)  \left(  z\right)  =\lim\limits_{D_{-}\ni
b^{\prime},\operatorname{Re}b^{\prime}\rightarrow\infty}\dfrac{\left(
z-b^{\prime}\right)  ^{N}}{2\pi i}%
{\displaystyle\int_{C}}
\dfrac{u\left(  \lambda\right)  }{z-\lambda}\left(  \lambda-b^{\prime}\right)
^{-N}d\lambda\text{ for }z\in D_{-}%
\end{array}
\right.  \text{,} \label{e4}%
\end{equation}
if they exist finitely. It should be noted that if they exist for an $N\geq0$,
then they exist also for any $N^{\prime}\geq N$ and take the same values.

The extended $\mathfrak{p}_{\pm}$ are well-defined for a certain $au$.

\begin{lemma}
\label{l1}For $N\in\mathbb{Z}_{+}$ if $a\in A(C)$ satisfies%
\begin{equation}
z^{N}\left(  a(z)-f(z)\right)  \text{ is bounded on }C \label{e5}%
\end{equation}
with a bounded analytic function $f$ on $D_{+}$, then, for $u\in H_{N}\left(
D_{+}\right)  $%
\begin{equation}
\left\{
\begin{array}
[c]{l}%
\mathfrak{p}_{+}\left(  au\right)  \left(  z\right)  =f\left(  z\right)
u\left(  z\right)  +\dfrac{1}{2\pi i}%
{\displaystyle\int_{C}}
\dfrac{\left(  a\left(  \lambda\right)  -f\left(  \lambda\right)  \right)
u\left(  \lambda\right)  }{\lambda-z}d\lambda\in H_{N}\left(  D_{+}\right)
\medskip\\
\mathfrak{p}_{-}\left(  au\right)  \left(  z\right)  =\dfrac{1}{2\pi i}%
{\displaystyle\int_{C}}
\dfrac{\left(  a\left(  \lambda\right)  -f\left(  \lambda\right)  \right)
u\left(  \lambda\right)  }{z-\lambda}d\lambda\in H\left(  D_{-}\right)
\end{array}
\right.  \label{e6}%
\end{equation}
hold. In particular for $u\in H\left(  D_{+}\right)  $ we have $\mathfrak{p}%
_{\pm}\left(  au\right)  \in H\left(  D_{\pm}\right)  $ respectively, and they
satisfy%
\[
\left\{
\begin{array}
[c]{l}%
\mathfrak{p}_{+}u=u\text{ \ if \ }u\in H_{N}\left(  D_{+}\right)  \smallskip\\
\mathfrak{p}_{+}u=0\text{ \ if \ }u\in H\left(  D_{-}\right)
\end{array}
\right.  \text{ and }\left\{
\begin{array}
[c]{l}%
\mathfrak{p}_{-}u=0\text{ \ if \ }u\in H_{N}\left(  D_{+}\right)  \smallskip\\
\mathfrak{p}_{-}u=u\text{ \ if \ }u\in H\left(  D_{-}\right)
\end{array}
\right.  \text{,}%
\]
which implies $H_{N}\left(  D_{+}\right)  \cap H\left(  D_{-}\right)
=\left\{  0\right\}  $.
\end{lemma}

\begin{proof}
If $u\in H_{N}\left(  D_{+}\right)  $ and $b^{\prime}\in D_{-}$,%
\begin{align*}
&  \dfrac{\left(  z-b^{\prime}\right)  ^{N}}{2\pi i}\int_{C}\dfrac{a\left(
\lambda\right)  u\left(  \lambda\right)  }{\lambda-z}\left(  \lambda
-b^{\prime}\right)  ^{-N}d\lambda\\
&  =\dfrac{1}{2\pi i}\int_{C}\dfrac{\left(  a\left(  \lambda\right)  -f\left(
\lambda\right)  \right)  u\left(  \lambda\right)  }{\lambda-z}\left(
\dfrac{z-b^{\prime}}{\lambda-b^{\prime}}\right)  ^{N}d\lambda+f\left(
z\right)  u\left(  z\right)
\end{align*}
holds for $z\in D_{+}$ due to $fu\left(  z-b^{\prime}\right)  ^{-N}\in
H\left(  D_{+}\right)  $. Since $\left(  a-f\right)  u\in L^{2}\left(
C\right)  $%
\[
\lim_{b^{\prime}\rightarrow\infty}\int_{C}\dfrac{\left(  a\left(
\lambda\right)  -f\left(  \lambda\right)  \right)  u\left(  \lambda\right)
}{\lambda-z}\left(  \dfrac{z-b^{\prime}}{\lambda-b^{\prime}}\right)
^{N}d\lambda=\int_{C}\dfrac{\left(  a\left(  \lambda\right)  -f\left(
\lambda\right)  \right)  u\left(  \lambda\right)  }{\lambda-z}d\lambda
\]
is valid, which shows%
\[
\mathfrak{p}_{+}\left(  au\right)  \left(  z\right)  =f\left(  z\right)
u\left(  z\right)  +\dfrac{1}{2\pi i}\int_{C}\dfrac{\left(  a\left(
\lambda\right)  -f\left(  \lambda\right)  \right)  u\left(  \lambda\right)
}{\lambda-z}d\lambda\in H_{N}\left(  D_{+}\right)  \text{.}%
\]
On the other hand, due to $fu\left(  z-b^{\prime}\right)  ^{-N}\in H\left(
D_{+}\right)  $%
\[
\left(  z-b^{\prime}\right)  ^{N}\int_{C}\dfrac{a\left(  \lambda\right)
u\left(  \lambda\right)  }{z-\lambda}\left(  \lambda-b^{\prime}\right)
^{-N}d\lambda=\int_{C}\dfrac{\left(  a\left(  \lambda\right)  -f\left(
\lambda\right)  \right)  u\left(  \lambda\right)  }{z-\lambda}\left(
\dfrac{z-b^{\prime}}{\lambda-b^{\prime}}\right)  ^{N}d\lambda
\]
holds for $z\in D_{-}$, and $\left(  a-f\right)  u\in L^{2}\left(  C\right)  $
implies%
\begin{align*}
\mathfrak{p}_{-}\left(  au\right)  \left(  z\right)   &  =\lim_{b^{\prime
}\rightarrow\infty}\dfrac{\left(  z-b^{\prime}\right)  ^{N}}{2\pi i}\int
_{C}\dfrac{a\left(  \lambda\right)  u\left(  \lambda\right)  }{z-\lambda
}\left(  \lambda-b^{\prime}\right)  ^{-N}d\lambda\\
&  =\dfrac{1}{2\pi i}\int_{C}\dfrac{\left(  a\left(  \lambda\right)  -f\left(
\lambda\right)  \right)  u\left(  \lambda\right)  }{z-\lambda}d\lambda\text{,}%
\end{align*}
which shows (\ref{e6}). If $u\in H\left(  D_{-}\right)  $, then due to $au\in
L^{2}\left(  C\right)  $ we easily have%
\[
\lim_{b^{\prime}\rightarrow\infty}\left(  z-b^{\prime}\right)  ^{N}\int
_{C}\dfrac{a\left(  \lambda\right)  u\left(  \lambda\right)  }{\lambda
-z}\left(  \lambda-b^{\prime}\right)  ^{-N}d\lambda=\int_{C}\dfrac{a\left(
\lambda\right)  u\left(  \lambda\right)  }{\lambda-z}d\lambda
\]
for $z\notin C$. Therefore, $\mathfrak{p}_{\pm}\left(  au\right)  \in H\left(
D_{\pm}\right)  $ respectively. The rest of the proof is clear.\medskip
\end{proof}

Consequently the projections $\mathfrak{p}_{\pm}$ can be extended to%
\begin{equation}
L_{N}^{2}\left(  C\right)  \equiv H_{N}\left(  D_{+}\right)  \oplus H\left(
D_{-}\right)  \text{ \ }\left(  \subset\left\vert \lambda\right\vert ^{N}%
L^{2}\left(  C\right)  \right)  \label{e7}%
\end{equation}
by (\ref{e7}) as projections. The norm in $L_{N}^{2}\left(  C\right)  $ is
defined by%
\[
\left\Vert u\right\Vert _{N}^{2}=\int_{C}\left\vert \left(  \lambda-b\right)
^{-N}\mathfrak{p}_{+}u\left(  \lambda\right)  \right\vert ^{2}\left\vert
d\lambda\right\vert +\int_{C}\left\vert \mathfrak{p}_{-}u\left(
\lambda\right)  \right\vert ^{2}\left\vert d\lambda\right\vert \text{.}%
\]
Moreover this lemma enables us to extend the Toeplitz operator $T\left(
a\right)  $ as a bounded operator on $H_{N}\left(  D_{+}\right)  $ for a
bounded function $a$ satisfying (\ref{e5}). Subsequently a subset
$A_{L}\left(  C\right)  $ of $A\left(  C\right)  $ for $L\in\mathbb{Z}_{+}$ is
introduced as follows:%
\begin{equation}
A_{L}\left(  C\right)  =\left\{
\begin{array}
[c]{l}%
a\in A\left(  C\right)  \text{; there exists an analytic function }f\text{ on
}D_{+}\\
\text{such that }\sup\limits_{z\in D_{+}}\left\vert f\left(  z\right)
\right\vert <\infty\text{, }\sup\limits_{\lambda\in C}\left\vert \lambda
^{L}\left(  a\left(  \lambda\right)  -f\left(  \lambda\right)  \right)
\right\vert <\infty
\end{array}
\right\}  \text{.} \label{e8}%
\end{equation}
Lemma \ref{l1} enables us to define the Toeplitz operator on $H_{N}\left(
D_{+}\right)  $ by%
\[
T_{N}(a)u=\mathfrak{p}_{+}\left(  au\right)  \in H_{N}\left(  D_{+}\right)
\text{.}%
\]
Let $L\geq N^{\prime}\geq N$. Then $\left\{  T_{N}\left(  a\right)  \right\}
_{N\geq0}$ has the property that if $a\in A_{L}\left(  C\right)  $, then
$\left.  T_{N^{\prime}}(a)\right\vert _{H_{N}\left(  D_{+}\right)  }=T_{N}%
(a)$. Therefore we use the notation%
\[
T(a)=T_{N}(a)\text{.}%
\]
The vector version of $A_{L}\left(  C\right)  $ and $T\left(  a\right)  $ are
defined by%
\begin{equation}
\left\{
\begin{array}
[c]{l}%
\boldsymbol{A}_{L}\left(  C\right)  =\left\{  \boldsymbol{a}=\left(
a_{1},a_{2}\right)  \text{; \ }a_{1}\text{, }a_{2}\in A_{L}\left(  C\right)
\right\}  \smallskip\\
\left(  T\left(  \boldsymbol{a}\right)  u\right)  \left(  z\right)  =T\left(
a_{1}\right)  u_{e}\left(  z\right)  +T\left(  a_{2}\right)  u_{o}\left(
z\right)
\end{array}
\right.  \text{.} \label{e9}%
\end{equation}
Set%
\begin{equation}
\boldsymbol{A}_{L}^{inv}\left(  C\right)  =\left\{  \boldsymbol{a}%
\in\boldsymbol{A}_{L}\left(  C\right)  \text{; }T(\boldsymbol{a})\text{ is
invertible on }H_{L}\left(  D_{+}\right)  \right\}  \text{.} \label{e10}%
\end{equation}

It should be noted that $\boldsymbol{A}_{L}^{inv}\left(  C\right)
\supset\boldsymbol{A}_{L^{\prime}}^{inv}\left(  C\right)  $ holds if
$L^{\prime}\geq L$.

\subsection{Characteristic functions and $m$-functions for $\boldsymbol{a}%
\in\boldsymbol{A}_{L}^{inv}\left(  C\right)  $}

In this section we define several quantities which will be necessary later
when $T(\boldsymbol{a})^{-1}$ exists.

For $\boldsymbol{a}\in\boldsymbol{A}_{2}^{inv}\left(  C\right)  $ one can
define two functions of $H\left(  D_{-}\right)  $ which characterize
$W_{\boldsymbol{a}}$ and are closely related to the tau-function introduced
later. For $\boldsymbol{a}\in\boldsymbol{A}_{2}^{inv}\left(  C\right)  $ set%
\begin{equation}
\left\{
\begin{array}
[c]{c}%
u\left(  z\right)  =\left(  T\left(  \boldsymbol{a}\right)  ^{-1}1\right)
\left(  z\right)  \in H_{1}\left(  D_{+}\right)  \smallskip\\
v(z)=\left(  T\left(  \boldsymbol{a}\right)  ^{-1}z\right)  \left(  z\right)
\in H_{2}\left(  D_{+}\right)
\end{array}
\right.  \text{,} \label{e2}%
\end{equation}
which is possible due to $1\in H_{1}\left(  D_{+}\right)  $, $z\in
H_{2}\left(  D_{+}\right)  $, and%
\begin{equation}
\left\{
\begin{array}
[c]{l}%
\varphi_{\boldsymbol{a}}\left(  z\right)  =\mathfrak{p}_{-}\left(
\boldsymbol{a}u\right)  \left(  z\right)  =a_{1}(z)u_{e}(z)+a_{2}%
(z)u_{o}(z)-1\in H\left(  D_{-}\right)  \smallskip\\
\psi_{\boldsymbol{a}}\left(  z\right)  =\mathfrak{p}_{-}\left(  \boldsymbol{a}%
v\right)  \left(  z\right)  =a_{1}(z)v_{e}(z)+a_{2}(z)v_{o}(z)-z\in H\left(
D_{-}\right)  \smallskip\\
\Delta_{\boldsymbol{a}}\left(  z\right)  =\dfrac{\left(  1+\varphi
_{\boldsymbol{a}}\left(  -z\right)  \right)  \left(  \psi_{\boldsymbol{a}%
}\left(  z\right)  +z\right)  -\left(  1+\varphi_{\boldsymbol{a}}\left(
z\right)  \right)  \left(  \psi_{\boldsymbol{a}}\left(  -z\right)  -z\right)
}{2z}%
\end{array}
\right.  \text{.} \label{53}%
\end{equation}

\begin{lemma}
\label{l12}If $\boldsymbol{a}\in\boldsymbol{A}_{2}^{inv}\left(  C\right)  $,
$\left\{  \varphi_{\boldsymbol{a}},\psi_{\boldsymbol{a}}\right\}  $ satisfies
the following properties.\newline(i) $\ \ \Delta_{\boldsymbol{a}}\left(
b\right)  \neq0$ on $D_{-}$ and%
\begin{equation}
T\left(  \boldsymbol{a}\right)  ^{-1}\dfrac{1}{z+b}=\dfrac{\left(
\varphi_{\boldsymbol{a}}\left(  b\right)  +1\right)  v-\left(  \psi
_{\boldsymbol{a}}\left(  b\right)  +b\right)  u}{\Delta_{\boldsymbol{a}%
}\left(  b\right)  \left(  z^{2}-b^{2}\right)  }\text{.} \label{62}%
\end{equation}
\newline(ii) $\ \left\{  \varphi_{\boldsymbol{a}},\psi_{\boldsymbol{a}%
}\right\}  $ determines $\boldsymbol{a}$.\newline(iii) There exist $\kappa
_{1}\left(  \boldsymbol{a}\right)  \in\mathbb{C}$ and $\phi_{\boldsymbol{a}%
}\in H\left(  D_{-}\right)  $ such that%
\begin{equation}
\varphi_{\boldsymbol{a}}\left(  z\right)  =\kappa_{1}\left(  \boldsymbol{a}%
\right)  z^{-1}+\phi_{\boldsymbol{a}}\left(  z\right)  z^{-1}\text{.}
\label{e11}%
\end{equation}

\end{lemma}

\begin{proof}
Here the suffix $\boldsymbol{a}$ is omitted. (\ref{53}) implies that for $b\in
D_{-}$ we have decompositions%
\[
\left\{
\begin{array}
[c]{l}%
\dfrac{a_{1}(z)u_{e}(z)+a_{2}(z)u_{o}(z)}{z^{2}-b^{2}}\\
=\dfrac{1}{2b}\left(  \dfrac{\varphi\left(  b\right)  +1}{z-b}-\dfrac
{\varphi\left(  -b\right)  +1}{z+b}\right)  +\dfrac{1}{2b}\left(
\dfrac{\varphi\left(  z\right)  -\varphi\left(  b\right)  }{z-b}%
-\dfrac{\varphi\left(  z\right)  -\varphi\left(  -b\right)  }{z+b}\right) \\
\dfrac{a_{1}(z)v_{e}(z)+a_{2}(z)v_{o}(z)}{z^{2}-b^{2}}\\
=\dfrac{1}{2b}\left(  \dfrac{\psi\left(  b\right)  +b}{z-b}-\dfrac{\psi\left(
-b\right)  -b}{z+b}\right)  +\dfrac{1}{2b}\left(  \dfrac{\psi\left(  z\right)
-\psi\left(  b\right)  }{z-b}-\dfrac{\psi\left(  z\right)  -\psi\left(
-b\right)  }{z+b}\right)
\end{array}
\right.
\]
into elements of $H_{2}\left(  D_{+}\right)  $ and $H\left(  D_{-}\right)  $,
hence%
\begin{equation}
\left\{
\begin{array}
[c]{c}%
T\left(  \boldsymbol{a}\right)  \dfrac{u}{z^{2}-b^{2}}=\dfrac{1}{2b}\left(
\dfrac{\varphi\left(  b\right)  +1}{z-b}-\dfrac{\varphi\left(  -b\right)
+1}{z+b}\right) \\
T\left(  \boldsymbol{a}\right)  \dfrac{v}{z^{2}-b^{2}}=\dfrac{1}{2b}\left(
\dfrac{\psi\left(  b\right)  +b}{z-b}-\dfrac{\psi\left(  -b\right)  -b}%
{z+b}\right)
\end{array}
\right.  \text{,} \label{51}%
\end{equation}
which yields%
\begin{equation}
\left\{
\begin{array}
[c]{l}%
T\left(  \boldsymbol{a}\right)  \dfrac{\left(  \varphi\left(  b\right)
+1\right)  v-\left(  \psi\left(  b\right)  +b\right)  u}{z^{2}-b^{2}}%
=\dfrac{\Delta\left(  b\right)  }{z+b}\smallskip\\
T\left(  \boldsymbol{a}\right)  \dfrac{\left(  \varphi\left(  -b\right)
+1\right)  v-\left(  \psi\left(  -b\right)  -b\right)  u}{z^{2}-b^{2}}%
=\dfrac{\Delta\left(  b\right)  }{z-b}%
\end{array}
\right.  \text{.} \label{52}%
\end{equation}
If $\Delta\left(  b\right)  =0$, then the invertibility of $T\left(
\boldsymbol{a}\right)  $ implies%
\[
\left\{
\begin{array}
[c]{l}%
\left(  \psi\left(  b\right)  +b\right)  u-\left(  \varphi\left(  b\right)
+1\right)  v=0\smallskip\\
\left(  \psi\left(  -b\right)  -b\right)  u-\left(  \varphi\left(  -b\right)
+1\right)  v=0
\end{array}
\right.  \text{.}%
\]
Applying $T\left(  \boldsymbol{a}\right)  $ we have%
\[
\left\{
\begin{array}
[c]{l}%
\left(  \psi\left(  b\right)  +b\right)  -\left(  \varphi\left(  b\right)
+1\right)  z=0\smallskip\\
\left(  \psi\left(  -b\right)  -b\right)  -\left(  \varphi\left(  -b\right)
+1\right)  z=0
\end{array}
\right.  \text{,}%
\]
which means%
\[
\psi\left(  b\right)  +b=\varphi\left(  b\right)  +1=\psi\left(  -b\right)
-b=\varphi\left(  -b\right)  +1=0\text{,}%
\]
hence (\ref{51}) implies $u=v=0$. This contradicts $T\left(  \boldsymbol{a}%
\right)  u=1$, $T\left(  \boldsymbol{a}\right)  v=z$, and we have
$\Delta\left(  b\right)  \neq0$. (\ref{62}) can be deduced from (\ref{52}).

Then, (\ref{52}) shows that for any rational function $f(z)$ with simple poles
only in $D_{-}$ and of order $O\left(  z^{-1}\right)  $ there exist two even
rational functions $r_{1}(z)$, $r_{2}(z)$ with the same property such that%
\[
T\left(  \boldsymbol{a}\right)  \left(  r_{1}u+r_{2}v\right)  (z)=f(z)
\]
holds, hence%
\[
r_{1}(z)u(z)+r_{2}(z)v(z)=\left(  T\left(  \boldsymbol{a}\right)
^{-1}f\right)  \left(  z\right)  \text{.}%
\]
Since such rational functions are dense in $H_{2}\left(  D_{+}\right)  $,
approximating $T\left(  \boldsymbol{a}\right)  1\boldsymbol{,}$ $T\left(
\boldsymbol{a}\right)  z\in H_{2}\left(  D_{+}\right)  $ by rational functions
$f_{n}$ and the continuity of $T\left(  \boldsymbol{a}\right)  ^{-1}$ show
that $T\left(  \boldsymbol{a}\right)  ^{-1}f_{n}$ converges to $1$, $z$
compact uniformly on $D_{+}$. Let
\[
\mathcal{Z}=\left\{  z\in D_{+}\text{; }\det\left(
\begin{array}
[c]{cc}%
u_{e}\left(  z\right)  & v_{e}\left(  z\right) \\
u_{o}\left(  z\right)  & v_{o}\left(  z\right)
\end{array}
\right)  =0\right\}  \text{.}%
\]
Then, the linear independence of $\left\{  u,v\right\}  $ implies the
discreteness of $\mathcal{Z}$, and hence for $z\in D_{+}\backslash\mathcal{Z}$
the associated $r_{n,1}(z)$, $r_{n,2}(z)$ converge to $\left\{  r_{1}%
^{(1)}(z),r_{2}^{(1)}(z)\right\}  $ and $\left\{  r_{1}^{(2)}(z),r_{2}%
^{(2)}(z)\right\}  $ depending on $T\left(  \boldsymbol{a}\right)  ^{-1}%
f_{n}\rightarrow1$ or $T\left(  \boldsymbol{a}\right)  ^{-1}f_{n}\rightarrow
z$. Since $r_{n,1}$, $r_{n,2}$ have no poles on $D_{+}$, their limits
$r_{j}^{(i)}\left(  z\right)  $ have no singularity on $D_{+}$ either, hence%
\begin{equation}
\left\{
\begin{array}
[c]{c}%
r_{1}^{(1)}(z)u(z)+r_{2}^{(1)}(z)v(z)=1\smallskip\\
r_{1}^{(2)}(z)u(z)+r_{2}^{(2)}(z)v(z)=z
\end{array}
\right.  \label{54}%
\end{equation}
holds for $z\in D_{+}$. Since $\left\{  r_{j}^{(i)}\right\}  $ are constructed
by $\left\{  \varphi,\psi\right\}  $, one see that $\left\{  r_{j}%
^{(i)}\right\}  $ depend only on $\left\{  \varphi,\psi\right\}  $, hence so
does $\left\{  u,v\right\}  $. On the other hand (\ref{53}) shows that
$a_{1}\left(  z\right)  $, $a_{2}\left(  z\right)  $ are recovered from
$\left\{  \varphi,\psi\right\}  $. We remark that%
\[
\det\left(
\begin{array}
[c]{cc}%
u_{e}\left(  z\right)  & v_{e}\left(  z\right) \\
u_{o}\left(  z\right)  & v_{o}\left(  z\right)
\end{array}
\right)  \neq0\text{ \ on \ }D_{+}\text{.}%
\]

(iii) is proved by applying (\ref{13}) of Lemma \ref{l1} to $\boldsymbol{a}%
\in\mathcal{A}_{2}^{inv}\left(  C\right)  $. Namely we have%
\begin{align*}
\varphi\left(  z\right)   &  =\mathfrak{p}_{-}\left(  a_{1}(z)u_{e}%
(z)+a_{2}(z)u_{o}(z)\right) \\
&  =\dfrac{1}{2\pi i}\int_{C}\dfrac{\left(  a_{1}(\lambda)-f_{1}%
(\lambda)\right)  u_{e}(\lambda)+\left(  a_{2}(\lambda)-f_{2}(\lambda)\right)
u_{o}(\lambda)}{z-\lambda}d\lambda\text{.}%
\end{align*}
Since $a_{j}(\lambda)-f_{j}(\lambda)=O\left(  \lambda^{-2}\right)  $ for
$j=1$, $2$ and $u_{e}$, $u_{o}\in H_{1}\left(  D_{+}\right)  $,%
\[
\lambda\left(  a_{1}(\lambda)-f_{1}(\lambda)\right)  u_{e}(\lambda)\text{,
\ }\lambda\left(  a_{2}(\lambda)-f_{2}(\lambda)\right)  u_{o}(\lambda)\in
L^{2}\left(  C\right)
\]
hold, hence%
\begin{align*}
&  \dfrac{1}{2\pi i}\int_{C}\dfrac{\left(  a_{1}(\lambda)-f_{1}(\lambda
)\right)  u_{e}(\lambda)+\left(  a_{2}(\lambda)-f_{2}(\lambda)\right)
u_{o}(\lambda)}{z-\lambda}d\lambda\\
&  =\kappa_{1}\left(  \boldsymbol{a}\right)  z^{-1}+\phi_{\boldsymbol{a}%
}\left(  z\right)  z^{-1}%
\end{align*}
is valid with%
\[
\left\{
\begin{array}
[c]{l}%
\kappa_{1}\left(  \boldsymbol{a}\right)  =\dfrac{1}{2\pi i}%
{\displaystyle\int_{C}}
\left(  \left(  a_{1}(\lambda)-f_{1}(\lambda)\right)  u_{e}(\lambda)+\left(
a_{2}(\lambda)-f_{2}(\lambda)\right)  u_{o}(\lambda)\right)  d\lambda
\smallskip\\
\phi_{\boldsymbol{a}}\left(  z\right)  =\dfrac{1}{2\pi i}%
{\displaystyle\int_{C}}
\dfrac{\lambda\left(  \left(  a_{1}(\lambda)-f_{1}(\lambda)\right)
u_{e}(\lambda)+\left(  a_{2}(\lambda)-f_{2}(\lambda)\right)  u_{o}%
(\lambda)\right)  }{z-\lambda}d\lambda
\end{array}
\right.  \text{.}%
\]
\medskip
\end{proof}

Owing to this Lemma we call $\left\{  \varphi_{\boldsymbol{a}},\psi
_{\boldsymbol{a}}\right\}  $ as the \textbf{characteristic} \textbf{functions}
of $\boldsymbol{a}$ or $W_{\boldsymbol{a}}$. This Lemma also implies the
possibility of a kind of Riemann-Hilbert factorization of $\left(
a_{1}(z),a_{2}(z)\right)  $, namely (\ref{53}), (\ref{54}) yield a
representation%
\begin{align*}
&  \left(
\begin{array}
[c]{cc}%
a_{1,e}(z) & a_{1,o}(z)\\
a_{2,o}(z) & a_{2,e}(z)
\end{array}
\right) \\
&  =\left(
\begin{array}
[c]{cc}%
1 & 0\\
0 & z
\end{array}
\right)  ^{-1}\left(
\begin{array}
[c]{cc}%
r_{1}^{(1)}(z) & r_{1}^{(1)}(z)\\
r_{2}^{(2)}(z) & r_{2}^{(2)}(z)
\end{array}
\right)  \left(
\begin{array}
[c]{cc}%
1+\varphi_{e}\left(  z\right)  & \varphi_{o}\left(  z\right) \\
\psi_{e}\left(  z\right)  & z+\psi_{o}\left(  z\right)
\end{array}
\right)  \text{,}%
\end{align*}
where the second term is analytic on $D_{+}$ and the third term is analytic on
$D_{-}$. The possibility of this factorization is very close to a sufficient
condition for the invertibility of $T(\boldsymbol{a})$.

The $m$-\textbf{function} for $W_{\boldsymbol{a}}$ is defined by%
\begin{equation}
m_{\boldsymbol{a}}\left(  z\right)  =\dfrac{z+\psi_{\boldsymbol{a}}\left(
z\right)  }{1+\varphi_{\boldsymbol{a}}\left(  z\right)  }+\kappa_{1}\left(
\boldsymbol{a}\right)  \text{.} \label{56}%
\end{equation}
$\kappa_{1}\left(  \boldsymbol{a}\right)  $ is added so that we have an
asymptotic behavior%
\begin{equation}
m_{\boldsymbol{a}}\left(  z\right)  =z+o\left(  1\right)  \text{ as
}z\rightarrow\infty\text{ in }\left\vert \arg\left(  \pm z\right)  \right\vert
<\pi/2-\epsilon\label{57}%
\end{equation}
(i) of (\ref{l12}) implies $1+\varphi_{\boldsymbol{a}}\left(  z\right)  $ is
not identically $0$, hence $m_{\boldsymbol{a}}$ is meromorphic on $D_{-}$.
Later we will see that $m_{\boldsymbol{a}}$ determines the potential $q$ and
is equal to the Weyl function under a certain condition on $\boldsymbol{a}%
$.\medskip

A good subset of $\boldsymbol{A}_{L}^{inv}\left(  C\right)  $ is given by
$\boldsymbol{M}_{L}\left(  C\right)  $:%
\begin{align}
&  \boldsymbol{M}_{L}\left(  C\right) \nonumber\\
&  =\left\{
\begin{array}
[c]{l}%
\boldsymbol{m}\left(  z\right)  =\left(  m_{1}\left(  z\right)  ,m_{2}\left(
z\right)  \right)  \text{; }\boldsymbol{m}\text{ is analytic on }%
\mathbb{C}\backslash\left(  \left[  -\mu_{0},\mu_{0}\right]  \cup
i\mathbb{R}\right) \\
\text{with }\mu_{0}=\sqrt{-\lambda_{0}}\text{ and satisfies (i), (ii) below:}%
\end{array}
\right\}  \label{9}%
\end{align}
(i) \ $\boldsymbol{m}(z)=\boldsymbol{1}+\sum\limits_{1\leq k<L}\boldsymbol{m}%
_{k}z^{-k}+O\left(  z^{-L}\right)  $ on $D_{-}$\ with $\boldsymbol{1}=\left(
1,1\right)  $, $\boldsymbol{m}_{k}\in\mathbb{C}^{2}$.\newline(ii)
$m_{1}(z)m_{2}(-z)+m_{1}(-z)m_{2}(z)\neq0$ on $\mathbb{C}\backslash\left(
\left[  -\mu_{0},\mu_{0}\right]  \cup i\mathbb{R}\right)  $\medskip\newline
For $\boldsymbol{m}$, $\boldsymbol{n}\in\boldsymbol{M}_{L}\left(  C\right)  $
define new elements by%
\[
\left\{
\begin{array}
[c]{l}%
\left(  \boldsymbol{m}\cdot\boldsymbol{n}\right)  (z)=\left(  m_{1}%
(z)n_{1,e}(z)+m_{2}(z)n_{1,o}(z),m_{1}(z)n_{2,o}(z)+m_{2}(z)n_{2,e}(z)\right)
\\
\widehat{\boldsymbol{m}}(z)=\left(  \dfrac{2\left(  m_{2,e}(z)-m_{1,o}%
(z)\right)  }{m_{1}(z)m_{2}(-z)+m_{1}(-z)m_{2}(z)},\dfrac{2\left(
m_{1,e}(z)-m_{2,o}(z)\right)  }{m_{1}(z)m_{2}(-z)+m_{1}(-z)m_{2}(z)}\right)
\end{array}
\right.  \text{.}%
\]
We have

\begin{lemma}
\label{l2}$\boldsymbol{M}_{L}\left(  C\right)  $ satisfies the group property:%
\begin{equation}
\left\{
\begin{array}
[c]{l}%
\boldsymbol{m}\cdot\boldsymbol{n}\text{, \ \ }\widehat{\boldsymbol{m}}%
\in\boldsymbol{M}_{L}\left(  C\right) \\
\boldsymbol{m}\cdot\widehat{\boldsymbol{m}}=\widehat{\boldsymbol{m}}%
\cdot\boldsymbol{m}=\boldsymbol{1}%
\end{array}
\right.  \text{.} \label{63}%
\end{equation}
Moreover it holds that%
\begin{equation}
T\left(  \boldsymbol{m}\cdot\boldsymbol{n}\right)  =T\left(  \boldsymbol{m}%
\right)  T\left(  \boldsymbol{n}\right)  \label{27}%
\end{equation}
for $\boldsymbol{m}$, $\boldsymbol{n}\in\boldsymbol{M}_{L}\left(  C\right)  $.
Consequently $\boldsymbol{M}_{L}\left(  C\right)  \subset\boldsymbol{A}%
_{L}^{inv}\left(  C\right)  $ is valid.
\end{lemma}

\begin{proof}
First note $\boldsymbol{M}_{L}\left(  C\right)  \subset\boldsymbol{A}%
_{L}\left(  C\right)  $. This is because for any $b\in D_{-}$ there exist
$\widetilde{\boldsymbol{m}}_{k}\in\mathbb{C}^{2}$ such that%
\[
\sum\limits_{1\leq k<L}\boldsymbol{m}_{k}z^{-k}=\sum\limits_{1\leq
k<L}\widetilde{\boldsymbol{m}}_{k}\left(  z-b\right)  ^{-k}+O\left(
z^{-L}\right)  \text{.}%
\]
The group property (\ref{63}) is clear. To show (\ref{27}) note $mH\left(
D_{-}\right)  \subset H\left(  D_{-}\right)  $ for any bounded analytic
function $m$ on $D_{-}$. Therefore, if scalers $m_{1}$, $m_{2}$ satisfy the
condition (i) of (\ref{9}), then for $u\in H_{N}\left(  D_{+}\right)  $%
\[
\mathfrak{p}_{+}\left(  m_{1}m_{2}u\right)  =\mathfrak{p}_{+}\left(
m_{1}\mathfrak{p}_{+}m_{2}u\right)  +\mathfrak{p}_{+}\left(  m_{1}%
\mathfrak{p}_{-}m_{2}u\right)  =\mathfrak{p}_{+}\left(  m_{1}\mathfrak{p}%
_{+}m_{2}u\right)
\]
holds. With this property in mind we have%
\begin{align*}
&  T\left(  \boldsymbol{m}\cdot\boldsymbol{n}\right)  u\\
&  =\mathfrak{p}_{+}\left(  \left(  m_{1}n_{1,e}+m_{2}n_{1,o}\right)
u_{e}+\left(  m_{1}n_{2,o}+m_{2}n_{2,e}\right)  u_{o}\right) \\
&  =\mathfrak{p}_{+}\left(  m_{1}\mathfrak{p}_{+}\left(  n_{1,e}u_{e}%
+n_{2,o}u_{o}\right)  \right)  +\mathfrak{p}_{+}\left(  m_{2}\mathfrak{p}%
_{+}\left(  n_{1,o}u_{e}+n_{2,e}u_{o}\right)  \right) \\
&  =T\left(  m_{1}\right)  \left(  T\left(  \boldsymbol{n}\right)  u\right)
_{e}+T\left(  m_{2}\right)  \left(  T\left(  \boldsymbol{n}\right)  u\right)
_{o}\\
&  =T\left(  \boldsymbol{m}\right)  T\left(  \boldsymbol{n}\right)  u\text{,}%
\end{align*}
which shows (\ref{27}).\medskip
\end{proof}

If $\boldsymbol{a}\left(  z\right)  =\boldsymbol{m}\left(  z\right)  =\left(
m_{1}\left(  z\right)  ,m_{2}\left(  z\right)  \right)  \in\boldsymbol{M}%
_{2}\left(  C\right)  $, then due to $\boldsymbol{m}\left(  z\right)
=\boldsymbol{1}+\boldsymbol{m}_{1}z^{-1}+O\left(  z^{-2}\right)  $ we have%
\[
T(\boldsymbol{m})1=\mathfrak{p}_{+}m_{1}(z)=1\text{, \ \ }T(\boldsymbol{m}%
)z=\mathfrak{p}_{+}m_{2}(z)z=z+m_{12}%
\]
with $\boldsymbol{m}_{1}=\left(  m_{11},m_{12}\right)  $, hence $u(z)=1$,
$v(z)=z-m_{12}$. Therefore%
\[
\left\{
\begin{array}
[c]{l}%
\varphi_{\boldsymbol{a}}\left(  z\right)  =m_{1}(z)-1\smallskip\\
\psi_{\boldsymbol{a}}\left(  z\right)  =-m_{12}m_{1}(z)+zm_{2}(z)-z
\end{array}
\right.
\]
follows, which yields%
\begin{equation}
m_{\boldsymbol{a}}(z)=\frac{-m_{12}m_{1}(z)+zm_{2}(z)}{m_{1}(z)}+m_{11}%
=\frac{zm_{2}(z)}{m_{1}(z)}+m_{11}-m_{12}\text{.} \label{35}%
\end{equation}

\section{Group action on $\boldsymbol{A}_{L}^{inv}\left(  C\right)  $}

The KdV flow is described by a group action on $\boldsymbol{A}_{L}%
^{inv}\left(  C\right)  $. For $m\in\mathbb{Z}_{-}$ and odd $n\in
\mathbb{Z}_{+}$ let $\Gamma_{n}^{\left(  m\right)  }$ be%
\begin{equation}
\left\{
\begin{array}
[c]{l}%
\Gamma_{n}^{\left(  m\right)  }=\left\{
\begin{array}
[c]{c}%
g=re^{h}\text{; }r\text{ is a rational function of order }m\text{ which}\\
\text{do not have poles nor zeros on }\left[  -\mu_{0},\mu_{0}\right]  \cup
i\mathbb{R}\text{,}\\
\text{and }h\text{ is a real odd polynomials of degree}\leq n
\end{array}
\right\} \\
\Gamma_{n}=\left\{  g=e^{h}\text{; }h\text{ is a real odd polynomials of
degree}\leq n\right\}  \subset\Gamma_{n}^{\left(  0\right)  }%
\end{array}
\right.  \text{,} \label{28}%
\end{equation}
where $\mu_{0}=\sqrt{-\lambda_{0}}$ and the order $m$ of a rational function
$r$ is defined by%
\[
m=\deg p-\deg q\text{ when }r=p/q\text{ with polynomials }p\text{, }q\text{.}%
\]
When we consider $g=re^{h}\in\Gamma_{n}^{\left(  m\right)  }$, the curve $C$
is taken so that $e^{h}$ remains bounded on $D_{+}$. Therefore, $C$ is
parametrized as%
\begin{equation}
C=\left\{
\begin{array}
[c]{c}%
\pm\omega\left(  y\right)  +iy\text{; \ }y\in\mathbb{R}\text{, }\omega\left(
y\right)  >0\text{, }\omega\left(  y\right)  =\omega\left(  -y\right)  \text{,
}\omega\text{ is}\\
\text{smooth and satisfies }\omega\left(  y\right)  =O\left(  y^{-\left(
n-1\right)  }\right)  \text{ as }y\rightarrow\infty
\end{array}
\right\}  \text{.} \label{93}%
\end{equation}
In most of the cases $\Gamma_{n}^{\left(  0\right)  }$ will be treated.
However, in some important cases we would like to consider $q_{\zeta}g$,
$q_{\zeta_{1}}q_{\zeta_{2}}g$ with $g\in\Gamma_{n}^{\left(  0\right)  }$ to
make arguments transparent, hence the numbers for $m$ which are frequently
appears are $0$, $-1$, $-2$. Note that any rational function $r$ can be
represented as a product of finite numbers of $q_{\zeta}\left(  z\right)
=\left(  1-\zeta^{-1}z\right)  ^{-1}$, $q_{\zeta}\left(  z\right)  ^{-1}%
$.\newline

For $\boldsymbol{a}\in\boldsymbol{A}_{L}\left(  C\right)  $, $g\in\Gamma
_{n}^{\left(  m\right)  }$ a natural product $g\boldsymbol{a}$ is bounded on
$C$ due to $m\leq0$, hence $g\boldsymbol{a}\in\boldsymbol{A}_{L}\left(
C\right)  .$ For $a\in A_{L}\left(  C\right)  $ and $g\in\Gamma_{n}^{\left(
m\right)  }$ when $L\geq N+m\geq0$, for $u\in H_{N}\left(  D_{+}\right)  $ an
identity%
\begin{equation}
T\left(  ga\right)  u=T\left(  a\right)  gu\in H_{N+m}\left(  D_{+}\right)
\label{11}%
\end{equation}
and%
\[
T\left(  ga\right)  :H_{N}\left(  D_{+}\right)  \rightarrow H_{N+m}\left(
D_{+}\right)  \subset H_{N}\left(  D_{+}\right)  \text{.}%
\]
hold. It should be noted that generally one cannot expect (\ref{11}) for
$\boldsymbol{a}\in\boldsymbol{A}_{L}\left(  C\right)  $ and $g\in\Gamma
_{n}^{\left(  m\right)  }$ unless $g$ is even.

The invertibility of $T(g\boldsymbol{a})$ is crucial in this paper and this
will be shown by using the tau-function, which is defined by the determinant
of the difference between $T(g\boldsymbol{a})$ and $T(\boldsymbol{a})$, namely%
\[
g^{-1}T(g\boldsymbol{a})T(\boldsymbol{a})^{-1}\text{.}%
\]
The tau-function describes the $\Gamma_{n}^{\left(  m\right)  }$ action very
well. To define the determinant we have to show the relevant operators are of
Hilbert-Schmidt type.

For $a\in A_{L}\left(  C\right)  $ let $f$ be an analytic function on $D_{+}$
such that%
\[
\sup_{z\in D_{+}}\left\vert f\left(  z\right)  \right\vert <\infty\text{ \ and
\ }\sup_{z\in C}\left\vert z^{L}\left(  a\left(  z\right)  -f\left(  z\right)
\right)  \right\vert <\infty\text{.}%
\]
For $g_{1}\in\Gamma_{n}^{\left(  0\right)  }$, $g_{2}\in\Gamma_{n}^{\left(
m\right)  }$ and fixed $b\in D_{-}$ define%
\begin{equation}
\left\{
\begin{array}
[c]{l}%
S_{\boldsymbol{a}}u\left(  z\right)  =\dfrac{1}{2\pi i}%
{\displaystyle\int_{C}}
\dfrac{\widetilde{\boldsymbol{a}}\left(  \lambda\right)  }{z-\lambda}u\left(
\lambda\right)  d\lambda\text{ \ for }u\in H_{N}\left(  D_{+}\right)
\smallskip\\
H_{g_{2}}u=\mathfrak{p}_{+}\left(  g_{2}u\right)  \text{ \ for }u\in H\left(
D_{-}\right)  \smallskip\\
R_{\boldsymbol{a}}\left(  g_{1},g_{2}\right)  u=\mathfrak{p}_{+}\left(
g_{2}\mathfrak{p}_{-}g_{1}\boldsymbol{a}u\right)  \text{ \ for }u\in
H_{N}\left(  D_{+}\right)  \smallskip
\end{array}
\right.  \text{.} \label{43}%
\end{equation}
The domains and images for the above maps are as follows:%
\[
\left\{
\begin{array}
[c]{ll}%
T\left(  g_{2}g_{1}\boldsymbol{a}\right)  & H_{N}\left(  D_{+}\right)
\rightarrow H_{N+m}\left(  D_{+}\right) \\
T\left(  g_{1}\boldsymbol{a}\right)  & H_{N}\left(  D_{+}\right)  \rightarrow
H_{N}\left(  D_{+}\right) \\
g_{2} & H_{N}\left(  D_{+}\right)  \rightarrow H_{N+m}\left(  D_{+}\right)
\end{array}%
\begin{array}
[c]{ll}%
S_{g_{1}\boldsymbol{a}} & H_{N}\left(  D_{+}\right)  \rightarrow H\left(
D_{-}\right) \\
H_{g_{2}} & H\left(  D_{-}\right)  \rightarrow H_{N+m}\left(  D_{+}\right) \\
&
\end{array}
\right.  \text{.}%
\]
$H_{g_{2}}$ is akin to a Hankel operator if $D_{+}$ is the unit disc. Recall
that the norms in $H_{N}\left(  D_{+}\right)  $ to $H\left(  D_{-}\right)  $
are respectively%
\[
\sqrt{%
{\displaystyle\int_{C}}
\left\vert u\left(  \lambda\right)  \right\vert ^{2}\left\vert \lambda
\right\vert ^{-2N}\left\vert d\lambda\right\vert }\text{, \ \ \ }\sqrt{%
{\displaystyle\int_{C}}
\left\vert u\left(  \lambda\right)  \right\vert ^{2}\left\vert d\lambda
\right\vert }\text{.}%
\]

\begin{lemma}
\label{l3}Assume $L\geq N\geq0$. Then we have\newline(i) $S_{\boldsymbol{a}}$
defines a Hilbert-Schmidt class operator from $H_{N}\left(  D_{+}\right)  $ to
$H\left(  D_{-}\right)  $ if%
\begin{equation}%
{\displaystyle\int_{C^{2}}}
\left\vert \dfrac{z^{N}\widetilde{a}_{j}\left(  z\right)  -\lambda
^{N}\widetilde{a}_{j}\left(  \lambda\right)  }{z-\lambda}\right\vert
^{2}\left\vert dz\right\vert \left\vert d\lambda\right\vert <\infty\text{
\ for }j=1\text{, }2\text{.} \label{39}%
\end{equation}
\newline(ii) Suppose $N+m\geq0$. Then $H_{g_{2}}$ is of Hilbert-Schmidt class
from $H\left(  D_{-}\right)  $ to $H_{N+m}\left(  D_{+}\right)  $ if%
\begin{equation}%
{\displaystyle\int_{C^{2}}}
\left\vert \frac{g_{2}\left(  \lambda\right)  -g_{2}\left(  z\right)
}{\lambda-z}\right\vert ^{2}\left\vert z\right\vert ^{-2\left(  N+m\right)
}\left\vert dz\right\vert \left\vert d\lambda\right\vert <\infty\text{.}
\label{48}%
\end{equation}
(iii) Identities%
\begin{equation}
\left\{
\begin{array}
[c]{l}%
T\left(  g_{2}g_{1}\boldsymbol{a}\right)  =g_{2}T\left(  g_{1}\boldsymbol{a}%
\right)  +R_{\boldsymbol{a}}\left(  g_{1},g_{2}\right)  \smallskip\\
R_{\boldsymbol{a}}\left(  g_{1},g_{2}\right)  =H_{g_{2}}S_{g_{1}%
\boldsymbol{a}}%
\end{array}
\right.  \label{40}%
\end{equation}
hold, and $R_{a}\left(  g_{1},g_{2}\right)  $ defines a trace class operator
from $H_{N}\left(  D_{+}\right)  $ to $H_{N}\left(  D_{+}\right)  $ under the
conditions (\ref{39}), (\ref{48}).
\end{lemma}

\begin{proof}
For $u\in H_{N}\left(  D_{+}\right)  $ it holds that%
\begin{align*}
T\left(  g_{2}g_{1}\boldsymbol{a}\right)  u  &  =\mathfrak{p}_{+}\left(
g_{2}\mathfrak{p}_{+}g_{1}\boldsymbol{a}u\right)  +\mathfrak{p}_{+}\left(
g_{2}\mathfrak{p}_{-}g_{1}\boldsymbol{a}u\right) \\
&  =g_{2}T\left(  g_{1}\boldsymbol{a}\right)  u+\mathfrak{p}_{+}\left(
g_{2}\mathfrak{p}_{-}g_{1}\boldsymbol{a}u\right) \\
&  =g_{2}T\left(  g_{1}\boldsymbol{a}\right)  u+H_{g_{2}}S_{g_{1}%
\boldsymbol{a}}u\text{.}%
\end{align*}
If $L\geq N$, we have%
\[
\left(  \mathfrak{p}_{-}\left(  g_{1}\boldsymbol{a}u\right)  \right)  \left(
z\right)  =\frac{1}{2\pi i}\int_{C}\frac{g_{1}\left(  \lambda\right)
\widetilde{\boldsymbol{a}}\left(  \lambda\right)  }{z-\lambda}u\left(
\lambda\right)  d\lambda
\]
for $u\in H_{N}\left(  D_{+}\right)  $, $z\in D_{-}$. Note here an identity%
\begin{align*}
\left(  S_{g_{1}\boldsymbol{a}}\right)  u\left(  z\right)   &  =\frac{1}{2\pi
i}\int_{C}\frac{g_{1}\left(  \lambda\right)  \left(  \lambda-b\right)
^{N}\widetilde{\boldsymbol{a}}\left(  \lambda\right)  }{z-\lambda}\left(
\lambda-b\right)  ^{-N}u\left(  \lambda\right)  d\lambda\\
&  =\frac{z^{-M}}{2\pi i}\int_{C}\frac{g_{1}\left(  \lambda\right)  \left(
\boldsymbol{s}\left(  \lambda\right)  -\boldsymbol{s}\left(  z\right)
\right)  }{z-\lambda}\left(  \lambda-b\right)  ^{-N}\lambda^{2N}u\left(
\lambda\right)  \lambda^{-2N}d\lambda
\end{align*}
for $z\in C$ with $\boldsymbol{s}\left(  \lambda\right)  =\left(
\lambda-b\right)  ^{N}\widetilde{\boldsymbol{a}}\left(  \lambda\right)  $ due
to $u\in H_{N}\left(  D_{+}\right)  $. Since we can regard $S_{g_{1}%
\boldsymbol{a}}$ as a map from $z^{N}L^{2}\left(  C\right)  $ to $L^{2}\left(
C\right)  $, we see that $S_{g_{1}\boldsymbol{a}}$ is of Hilbert-Schmidt class
from $H_{N}\left(  D_{+}\right)  $ to $H\left(  D_{-}\right)  $ if%
\[
\int_{C^{2}}\left\vert g_{1}\left(  \lambda\right)  \frac{\boldsymbol{s}%
\left(  \lambda\right)  -\boldsymbol{s}\left(  z\right)  }{z-\lambda}\left(
\lambda-b\right)  ^{-N}\lambda^{2N}\right\vert ^{2}\left\vert \lambda
\right\vert ^{-2N}\left\vert dz\right\vert \left\vert d\lambda\right\vert
<\infty\text{,}%
\]
which is equivalent to (\ref{39}) if we replace $\boldsymbol{s}\left(
\lambda\right)  $ by $\lambda^{N}\widetilde{a}_{j}$ here.

On the other hand the assumption $m\leq0$ implies $\sup_{z\in D_{+}}\left\vert
g_{2}\left(  z\right)  \right\vert <\infty$, hence $H_{g_{2}}$ defines an
operator from $H\left(  D_{-}\right)  $ to $H\left(  D_{+}\right)  $. We find
a condition for $H_{g_{2}}$ also to be of Hilbert-Schmidt class. Note%
\[
H_{g_{2}}u\left(  z\right)  =\frac{1}{2\pi i}\int_{C}\frac{g_{2}\left(
\lambda\right)  }{\lambda-z}u\left(  \lambda\right)  d\lambda=\frac{1}{2\pi
i}\int_{C}\frac{g_{2}\left(  \lambda\right)  -g_{2}\left(  z\right)  }%
{\lambda-z}u\left(  \lambda\right)  d\lambda
\]
for $u\in H\left(  D_{-}\right)  $, $z\in D_{+}$. Then, $H_{g_{2}}$ is of
Hilbert-Schmidt class from $H\left(  D_{-}\right)  $ to $H_{N+m}\left(
D_{+}\right)  $ if%
\[%
{\displaystyle\int_{C^{2}}}
\left\vert \frac{g_{2}\left(  \lambda\right)  -g_{2}\left(  z\right)
}{\lambda-z}\right\vert ^{2}\left\vert z\right\vert ^{-2\left(  N+m\right)
}\left\vert dz\right\vert \left\vert d\lambda\right\vert
\]
which is (\ref{48}).\medskip
\end{proof}

For later purpose we find a sufficient condition on $g\in\Gamma_{n}^{\left(
m\right)  }$ under which (\ref{48}) is satisfied. From now on we assume
without loss of generality the curve $C$ fulfills%
\begin{equation}
\sup_{z\in C}%
{\displaystyle\int_{\left\vert z-\lambda\right\vert \leq1,\lambda\in C}}
\left\vert d\lambda\right\vert <\infty\text{,} \label{32}%
\end{equation}
and there exists a neighborhood $U$ of the closure of $D_{+}$ and $\epsilon>0$
such that%
\begin{equation}
z\text{, }\lambda\in C\text{, }\left\vert z-\lambda\right\vert \leq
\epsilon\Longrightarrow\left(  z-\lambda\right)  t+\lambda\in U\text{ for
}t\in\left[  0,1\right]  \text{.} \label{34}%
\end{equation}
For $g\in\Gamma_{n}^{\left(  m\right)  }$ let $c>0$ be a constant such that%
\begin{equation}
\left\{
\begin{array}
[c]{l}%
c^{-1}\left\vert z\right\vert ^{m}\leq\left\vert g\left(  z\right)
\right\vert \leq c\left\vert z\right\vert ^{m}\\
\left\vert g^{\prime}\left(  z\right)  \right\vert \leq c\left\vert
z\right\vert ^{m+n-1}%
\end{array}
\right.  \text{ hold for }z\in U\text{.} \label{8}%
\end{equation}
For $N\in\mathbb{Z}_{+}$ set%
\[
\Delta=%
{\displaystyle\int_{C^{2}}}
\left\vert \dfrac{g\left(  z\right)  -g\left(  \lambda\right)  }{z-\lambda
}\right\vert ^{2}\left\vert z\right\vert ^{-2\left(  N+m\right)  }\left\vert
dz\right\vert \left\vert d\lambda\right\vert \text{,}%
\]
which is the square of the Hilbert-Schmidt norm of the operator
\[
H_{g}:H\left(  D_{-}\right)  \rightarrow H_{N+m}\left(  D_{+}\right)  .
\]

\begin{lemma}
\label{l10}If $N\geq\max\left\{  n,1-m\right\}  $ hold, there exists a
constant $c_{0}$ depending only on $c$ such that $\Delta\leq c_{0}$.
\end{lemma}

\begin{proof}
Let $\epsilon$ be $0<\epsilon<1$. We first show that there exists a constant
$c_{1}$ depending on the constant $c$ of (\ref{8}) such that%
\begin{equation}
\left\vert \dfrac{g\left(  z\right)  -g\left(  \lambda\right)  }{z-\lambda
}\right\vert \leq c_{1}\left\{
\begin{array}
[c]{ll}%
\left\vert z\right\vert ^{m+n-1} & \text{if }\left\vert z-\lambda\right\vert
\leq\epsilon\left\vert \lambda\right\vert \\
\left\vert z\right\vert ^{-1}\left(  \left\vert z\right\vert ^{m}+\left\vert
\lambda\right\vert ^{m}\right)  & \text{if }\left\vert z-\lambda\right\vert
>\epsilon\left\vert \lambda\right\vert
\end{array}
\right.  \label{64}%
\end{equation}
holds for $z$, $\lambda\in C$. Since%
\[
\dfrac{g\left(  z\right)  -g\left(  \lambda\right)  }{z-\lambda}=\int_{0}%
^{1}g^{\prime}\left(  \left(  \lambda-z\right)  t+z\right)  dt
\]
for $z$, $\lambda\in C$, the properties (\ref{34}), (\ref{8}) show that there
exists a constant $c_{1}$ such that%
\[
\left\vert \dfrac{g\left(  z\right)  -g\left(  \lambda\right)  }{z-\lambda
}\right\vert \leq c_{1}\left\vert z\right\vert ^{m+n-1}I_{\left\vert
z-\lambda\right\vert \leq\epsilon\left\vert \lambda\right\vert }\text{.}%
\]
The other estimate is clear and we have (\ref{64}).

Note $\Delta\leq\Delta_{1}+\Delta_{2}$ with%
\[
\left\{
\begin{array}
[c]{l}%
\Delta_{1}=c_{1}^{2}%
{\displaystyle\int_{\left\vert z-\lambda\right\vert \leq\epsilon\left\vert
\lambda\right\vert }}
\left\vert z\right\vert ^{2\left(  m+n-1\right)  }\left\vert z\right\vert
^{-2\left(  N+m\right)  }\left\vert dz\right\vert \left\vert d\lambda
\right\vert \\
\Delta_{2}=c_{1}^{2}%
{\displaystyle\int_{\left\vert z-\lambda\right\vert >\epsilon\left\vert
\lambda\right\vert }}
\left\vert \lambda\right\vert ^{-2}\left(  \left\vert z\right\vert
^{m}+\left\vert \lambda\right\vert ^{m}\right)  ^{2}\left\vert z\right\vert
^{-2\left(  N+m\right)  }\left\vert dz\right\vert \left\vert d\lambda
\right\vert
\end{array}
\right.  \text{.}%
\]
Since the exponent of the integrand of $\Delta_{1}$ is equal to%
\[
2\left(  m+n-1\right)  -2\left(  N+m\right)  \text{,}%
\]
$\Delta_{1}<\infty$ is valid if $N>n-1/2$. On the other hand, there exists a
constant $c_{2}$ such that%
\[
\Delta_{2}\leq c_{2}%
{\displaystyle\int_{\left\vert z-\lambda\right\vert >\epsilon\left\vert
\lambda\right\vert }}
\left\vert \lambda\right\vert ^{-2}\left(  \left\vert z\right\vert
^{2m}+\left\vert \lambda\right\vert ^{2m}\right)  \left\vert z\right\vert
^{-2\left(  N+m\right)  }\left\vert dz\right\vert \left\vert d\lambda
\right\vert \text{.}%
\]
The right side is finite if%
\[
-2N<-1\text{, }-2(N+m)<-1
\]
which is equivalent to $N+m\geq1$. The above constants $c_{1}$, $c_{2}$ can be
chosen depending on $c$, hence so does $c_{0}$.\medskip
\end{proof}

The dependence of the constant $c_{0}$ on the constant $c$ will be used in the
proof of the continuity of the tau-function later.

\section{Derivation of Schr\"{o}dinger operator and KdV equation}

Schr\"{o}dinger operators and solutions to the KdV equation can be obtained
from $T\left(  e^{xz}\boldsymbol{a}\right)  $, $T\left(  e^{xz+tz^{3}%
}\boldsymbol{a}\right)  $ under their invertibility. This section is devoted
to the rigorous derivation of these equations.

\subsection{Differentiability}

The KdV flow is constructed by one-parameter group $g_{t}\left(  z\right)
=e^{th\left(  z\right)  }$ with odd polynomial $h$, and for the construction
the differentiability of $T\left(  g_{t}\boldsymbol{a}\right)  $ with respect
to $t$ will be necessary. In this section we extend the definition $T\left(
a\right)  $. For a polynomial $h$ of degree $n$%
\[
hu\in H_{N+n}\left(  D_{+}\right)  \text{ \ \ if \ }u\in H_{N}\left(
D_{+}\right)  \text{,}%
\]
so for $a\in A_{L}\left(  C\right)  $ define%
\[
T\left(  ha\right)  u=\mathfrak{p}_{+}\left(  hau\right)  =T\left(  a\right)
hu\in H_{N+n}\left(  D_{+}\right)  \text{,}%
\]
which is possible if $L\geq N+n$. For $\boldsymbol{a}=\left(  a_{1}%
,a_{2}\right)  \in\boldsymbol{A}_{L}\left(  C\right)  $ we define%
\[
T\left(  h\boldsymbol{a}\right)  u=T\left(  ha_{1}\right)  u_{e}+T\left(
ha_{2}\right)  u_{o}%
\]

\begin{lemma}
\label{l5}Let $\boldsymbol{a}\in\boldsymbol{A}_{L}\left(  C\right)  $ and
$g_{t}\left(  z\right)  =e^{th\left(  z\right)  }\in\Gamma_{n}$. Assume
$g_{t}\boldsymbol{a}\in\boldsymbol{A}_{N+n}^{inv}\left(  C\right)  $ for any
$t\in\mathbb{R}$. Then if $L\geq N+n$, for any $u\in H_{N}\left(
D_{+}\right)  $%
\[
\left\{
\begin{array}
[c]{l}%
\partial_{t}T\left(  g_{t}\boldsymbol{a}\right)  u=T\left(  hg_{t}%
\boldsymbol{a}\right)  u\in H_{N+n}\left(  D_{+}\right)  \smallskip\\
\partial_{t}T\left(  g_{t}\boldsymbol{a}\right)  ^{-1}u=-T\left(
g_{t}\boldsymbol{a}\right)  ^{-1}T\left(  hg_{t}\boldsymbol{a}\right)
T\left(  g_{t}\boldsymbol{a}\right)  ^{-1}u\in H_{N+n}\left(  D_{+}\right)
\end{array}
\right.
\]
holds. Any higher derivative $\partial_{t}^{k}T\left(  g_{t}\boldsymbol{a}%
\right)  ^{-1}u$ exists if $L\geq N+kn$.
\end{lemma}

\begin{proof}
Let $N_{1}=N+n$. Recall $T\left(  \boldsymbol{a}\right)  u=T\left(
a_{1}\right)  u_{e}+T\left(  a_{2}\right)  u_{o}$ if $\boldsymbol{a}=\left(
a_{1},a_{2}\right)  $ (see (\ref{16})). The first identity follows easily from%
\[
\dfrac{T\left(  g_{t}\boldsymbol{a}\right)  u-T\left(  g_{s}\boldsymbol{a}%
\right)  u}{t-s}=\dfrac{1}{t-s}\int_{s}^{t}\mathfrak{p}_{+}\left(  g_{\tau
}h\boldsymbol{a}u\right)  d\tau=\dfrac{1}{t-s}\int_{s}^{t}T\left(  hg_{\tau
}\boldsymbol{a}\right)  ud\tau\text{.}%
\]
To show the second identity first we verify the continuity of $T\left(
g_{t}\boldsymbol{a}\right)  ^{-1}u$ in $H_{N_{1}}\left(  D_{+}\right)  $ with
respect to $t$. Applying (ii) of Lemma \ref{l3} with $g_{1}=1$, $g_{2}=g$ and
replacing $N$ by $N_{1}$, we have%
\[
T\left(  g_{t}\boldsymbol{a}\right)  =g_{t}T\left(  \boldsymbol{a}\right)
+R_{\boldsymbol{a}}\left(  1,g_{t}\right)
\]
with%
\[
R_{\boldsymbol{a}}\left(  1,g_{t}\right)  =H_{g_{t}}S_{\boldsymbol{a}}\text{.}%
\]
Therefore%
\[
T\left(  g_{t}\boldsymbol{a}\right)  ^{-1}=\left(  I+T\left(  \boldsymbol{a}%
\right)  ^{-1}g_{t}^{-1}R_{\boldsymbol{a}}\left(  1,g_{t}\right)  \right)
^{-1}T\left(  \boldsymbol{a}\right)  ^{-1}g_{t}^{-1}%
\]
holds. We show $g_{t}^{-1}R_{\boldsymbol{a}}\left(  1,g_{t}\right)  $ is
continuous in the Hilbert-Schmidt norm on $H_{N_{1}}\left(  D_{+}\right)  $,
which is reduced to that of $g_{t}^{-1}H_{g_{t}}$ as an operator from
$H\left(  D_{-}\right)  $ to $H_{N_{1}}\left(  D_{+}\right)  $. The HS-norm of
$g_{t}^{-1}H_{g_{t}}$ is%
\begin{equation}
\Delta=%
{\displaystyle\int_{C^{2}}}
\left\vert \frac{g_{t}\left(  \lambda\right)  g_{t}\left(  z\right)
^{-1}-g_{s}\left(  \lambda\right)  g_{s}\left(  z\right)  ^{-1}}{\lambda
-z}\right\vert ^{2}\left\vert z\right\vert ^{-2N_{1}}\left\vert dz\right\vert
\left\vert d\lambda\right\vert \text{.} \label{112}%
\end{equation}
The proof is carried out similarly to that of (\ref{48}). Observe%
\begin{align*}
\left\vert \dfrac{g_{t}\left(  \lambda\right)  g_{t}\left(  z\right)
^{-1}-g_{s}\left(  \lambda\right)  g_{s}\left(  z\right)  ^{-1}}{\lambda
-z}\right\vert  &  \leq\left\vert \dfrac{h\left(  \lambda\right)  -h\left(
z\right)  }{\lambda-z}\right\vert \int_{s}^{t}\left\vert e^{\tau\left(
h\left(  \lambda\right)  -h\left(  z\right)  \right)  }\right\vert d\tau\\
&  \leq c\left(  \left\vert z\right\vert ^{n-1}+\left\vert \lambda\right\vert
^{n-1}\right)  \left(  t-s\right)
\end{align*}
for $z$, $\lambda\in C$. Then separating the integral (\ref{112}) on
$\left\vert \lambda-z\right\vert \leq\epsilon\left\vert \lambda\right\vert $
and $\left\vert \lambda-z\right\vert >\epsilon\left\vert \lambda\right\vert $,
we have%
\begin{align*}
\Delta &  \leq c_{1}\left(  t-s\right)  ^{2}\int_{\left\vert \lambda
-z\right\vert \leq\epsilon\left\vert \lambda\right\vert }\left(  \left\vert
z\right\vert ^{2\left(  n-1\right)  }+\left\vert \lambda\right\vert ^{2\left(
n-1\right)  }\right)  \left\vert z\right\vert ^{-2N_{1}}\left\vert
dz\right\vert \left\vert d\lambda\right\vert \\
&  +c_{1}\int_{\left\vert \lambda-z\right\vert >\epsilon\left\vert
\lambda\right\vert }\left\vert \lambda\right\vert ^{-2}\left\vert g_{t}\left(
\lambda\right)  g_{t}\left(  z\right)  ^{-1}-g_{s}\left(  \lambda\right)
g_{s}\left(  z\right)  ^{-1}\right\vert ^{2}\left\vert z\right\vert ^{-2N_{1}%
}\left\vert dz\right\vert \left\vert d\lambda\right\vert \text{.}%
\end{align*}
The first term is dominated by $c_{2}\left(  t-s\right)  ^{2}$ if%
\[
2\left(  n-1\right)  -2N_{1}<-1\Longrightarrow N_{1}\geq n\text{,}%
\]
which is satisfied if $N\geq0$. The second term tends to $0$ as $s\rightarrow
t$ if $-2N_{1}<-1$, which is always valid if $N_{1}\geq1$. Therefore we have
the continuity of $g_{t}^{-1}R_{\boldsymbol{a}}\left(  1,g_{t}\right)  $ in
the HS-norm on $H_{N_{1}}\left(  D_{+}\right)  $, which implies $T\left(
g_{t}\boldsymbol{a}\right)  ^{-1}u$ is continuous in $t$ for any fixed $u\in
H_{N_{1}}\left(  D_{+}\right)  $ if $L\geq N+n$. Consequently noting the
identity%
\[
\epsilon^{-1}\left(  T\left(  g_{t+\epsilon}\boldsymbol{a}\right)
^{-1}-T\left(  g_{t}\boldsymbol{a}\right)  ^{-1}\right)  =T\left(
g_{t+\epsilon}\boldsymbol{a}\right)  ^{-1}\epsilon^{-1}\left(  T\left(
g_{t}\boldsymbol{a}\right)  -T\left(  g_{t+\epsilon}\boldsymbol{a}\right)
\right)  T\left(  g_{t}\boldsymbol{a}\right)  ^{-1}\text{,}%
\]
we have the Lemma. The existence of higher derivatives can be shown similarly.
\end{proof}

\subsection{Derivation of Schr\"{o}dinger operator}

First we derive a Schr\"{o}dinger operator from $\boldsymbol{a}\in
\boldsymbol{A}_{1}\left(  C\right)  $ and $g=e_{x}$ with%
\[
e_{x}\left(  z\right)  =e^{xz}\text{.}%
\]
The curve $C$ is chosen so that $e_{x}\left(  z\right)  $ remains bounded for
any fixed $x\in\mathbb{R}$, namely%
\begin{equation}
C=\left\{
\begin{array}
[c]{c}%
\pm\omega\left(  y\right)  +iy\text{, }y\in\mathbb{R}\text{; }\omega\left(
y\right)  \text{ is a positive even smooth}\\
\text{function on }\mathbb{R}\text{ such that }\omega\left(  y\right)
=O\left(  1\right)  \text{ as }\left\vert y\right\vert \rightarrow\infty
\end{array}
\right\}  \text{.} \label{24}%
\end{equation}
Recall%
\[
\boldsymbol{a}\left(  z\right)  f\left(  z\right)  =a_{1}(z)f_{e}%
(z)+a_{2}(z)f_{o}(z)
\]
for a vector function $\boldsymbol{a}\left(  z\right)  =\left(  a_{1}%
(z),a_{2}(z)\right)  $ and a function $f(z)$ on $C$. For $L\geq3$,
$\boldsymbol{a}\in\boldsymbol{A}_{L}\left(  C\right)  $ assume $e_{x}%
\boldsymbol{a}\in\boldsymbol{A}_{L}^{inv}\left(  C\right)  $ for any
$x\in\mathbb{R}$. Let $u_{x}\in H_{1}\left(  D_{+}\right)  $ be%
\[
u_{x}=T\left(  e_{x}\boldsymbol{a}\right)  ^{-1}1\in H_{1}\left(
D_{+}\right)
\]
and set%
\[
w_{x}=\mathfrak{p}_{-}\left(  e_{x}\boldsymbol{a}u_{x}\right)  \in H\left(
D_{-}\right)  \text{.}%
\]
Then, for a bounded analytic vector $\boldsymbol{f}\left(  z\right)  $ on
$D_{+}$ satisfying $\widetilde{\boldsymbol{a}}\left(  \lambda\right)
=\boldsymbol{a}\left(  \lambda\right)  -\boldsymbol{f}\left(  \lambda\right)
=O\left(  \lambda^{-L}\right)  $ on $C$%
\[
w_{x}\left(  z\right)  =\dfrac{1}{2\pi i}\int_{C}\dfrac{e^{x\lambda}%
\widetilde{\boldsymbol{a}}\left(  \lambda\right)  }{z-\lambda}u_{x}\left(
\lambda\right)  d\lambda
\]
holds. Since Lemma \ref{l5} implies $\partial_{x}^{j}u_{x}\in H_{j+1}\left(
D_{+}\right)  $ for $j\leq L-1$,%
\begin{equation}
\left\{
\begin{array}
[c]{l}%
\partial_{x}w_{x}\left(  z\right)  =\dfrac{1}{2\pi i}%
{\displaystyle\int_{C}}
e^{x\lambda}\dfrac{\lambda\widetilde{\boldsymbol{a}}\left(  \lambda\right)
u_{x}\left(  \lambda\right)  +\widetilde{\boldsymbol{a}}\left(  \lambda
\right)  \partial_{x}u_{x}\left(  \lambda\right)  }{z-\lambda}d\lambda\\
\partial_{x}^{2}w_{x}\left(  z\right)  =\dfrac{1}{2\pi i}%
{\displaystyle\int_{C}}
e^{x\lambda}\dfrac{\lambda^{2}\widetilde{\boldsymbol{a}}\left(  \lambda
\right)  u_{x}\left(  \lambda\right)  +2\lambda\widetilde{\boldsymbol{a}%
}\left(  \lambda\right)  \partial_{x}u_{x}\left(  \lambda\right)
+\widetilde{\boldsymbol{a}}\left(  \lambda\right)  \partial_{x}^{2}%
u_{x}\left(  \lambda\right)  }{z-\lambda}d\lambda
\end{array}
\right.  \text{.} \label{33}%
\end{equation}
Since $u_{x}\in H_{1}\left(  D_{+}\right)  $, the expansion%
\[
\left(  z-\lambda\right)  ^{-1}=\sum_{1\leq k\leq M}z^{-k}\lambda^{k-1}%
+z^{-M}\lambda^{M}\left(  z-\lambda\right)  ^{-1}%
\]
shows%
\begin{align}
w_{x}\left(  z\right)   &  =\sum_{1\leq k\leq L-1}z^{-k}\dfrac{1}{2\pi i}%
\int_{C}\lambda^{k-1}e^{x\lambda}\widetilde{\boldsymbol{a}}\left(
\lambda\right)  u_{x}\left(  \lambda\right)  d\lambda\nonumber\\
&  \text{ \ \ \ \ \ }+z^{-L+1}\dfrac{1}{2\pi i}\int_{C}\dfrac{\lambda
^{L-1}e^{x\lambda}\widetilde{\boldsymbol{a}}\left(  \lambda\right)
}{z-\lambda}u_{x}\left(  \lambda\right)  d\lambda\nonumber\\
&  \equiv\sum_{1\leq k\leq L-1}z^{-k}s_{k}(x)+z^{-L+1}\widetilde{w}_{x}(z)
\label{10}%
\end{align}
\ with $\widetilde{w}_{x}\in H\left(  D_{-}\right)  $. Since $\partial
_{x}u_{x}\in H_{2}\left(  D_{+}\right)  $, $\partial_{x}^{2}u_{x}\in
H_{3}\left(  D_{+}\right)  $, (\ref{33}) shows similarly%
\begin{equation}
\left\{
\begin{array}
[c]{c}%
\partial_{x}w_{x}\left(  z\right)  =\sum_{1\leq k\leq L-2}z^{-k}s_{k}^{\prime
}(x)+z^{-L+2}\widetilde{w}_{x}^{\left(  1\right)  }(z)\\
\partial_{x}^{2}w_{x}\left(  z\right)  =\sum_{1\leq k\leq L-3}z^{-k}%
s_{k}^{\prime\prime}(x)+z^{-L+3}\widetilde{w}_{x}^{\left(  2\right)  }(z)
\end{array}
\in H\left(  D_{-}\right)  \right.  \text{.} \label{18}%
\end{equation}
with $\widetilde{w}_{x}^{\left(  1\right)  }$, $\widetilde{w}_{x}^{\left(
2\right)  }\in H\left(  D_{-}\right)  $. The notations in (\ref{53}) and Lemma
\ref{l12} imply%
\begin{equation}
\left\{
\begin{array}
[c]{l}%
w_{x}\left(  z\right)  =\varphi_{e_{x}\boldsymbol{a}}\left(  z\right) \\
s_{1}(x)=\kappa_{1}\left(  e_{x}\boldsymbol{a}\right)
\end{array}
\right.  \text{.} \label{4}%
\end{equation}
Set%
\begin{equation}
f_{\boldsymbol{a}}\left(  x,z\right)  =\boldsymbol{a}\left(  z\right)
u_{x}\left(  z\right)  =e^{-xz}\left(  1+\varphi_{e_{x}\boldsymbol{a}}\left(
z\right)  \right)  \text{.} \label{13}%
\end{equation}

\begin{proposition}
\label{p1}Let $L\geq3$ and assume $e_{x}\boldsymbol{a}\in\boldsymbol{A}%
_{L}^{inv}\left(  C\right)  $ for any $x\in\mathbb{R}$. Set%
\[
q(x)=-2s_{1}^{\prime}(x)=-2\partial_{x}\kappa_{1}\left(  e_{x}\boldsymbol{a}%
\right)  \text{.}%
\]
Then, a Schr\"{o}dinger equation%
\begin{equation}
-\partial_{x}^{2}f_{\boldsymbol{a}}\left(  x,z\right)  +q(x)f_{\boldsymbol{a}%
}\left(  x,z\right)  =-z^{2}f_{\boldsymbol{a}}\left(  x,z\right)  \label{22}%
\end{equation}
holds, and $\left\{  s_{k}(x)\right\}  _{2\leq k\leq L-2}$ in (\ref{10}) is
determined by a recurrence relation%
\begin{equation}
s_{k}^{\prime\prime}+2s_{1}^{\prime}s_{k}-2s_{k+1}^{\prime}=0\text{,\ \ }%
\left(  1\leq k\leq L-3\right)  \label{21}%
\end{equation}
for given $s_{1}\left(  x\right)  $, $\left\{  s_{k}(0)\right\}  _{2\leq k\leq
L-2}$.
\end{proposition}

\begin{proof}
The identity $w_{x}=e_{x}\boldsymbol{a}u_{x}-1$ yields%
\[
\left\{
\begin{array}
[c]{l}%
\partial_{x}w_{x}=e_{x}\left(  z\boldsymbol{a}u_{x}+\boldsymbol{a}\partial
_{x}u_{x}\right) \\
\partial_{x}^{2}w_{x}=e_{x}\left(  z^{2}\boldsymbol{a}u_{x}+\boldsymbol{a}%
\partial_{x}^{2}u_{x}+2z\boldsymbol{a}\partial_{x}u_{x}\right)
\end{array}
\right.  \text{,}%
\]
which implies%
\begin{equation}
\partial_{x}^{2}w_{x}-2z\partial_{x}w_{x}=e_{x}\boldsymbol{a}\left(
\partial_{x}^{2}u_{x}-z^{2}u_{x}\right)  \text{.} \label{23}%
\end{equation}
Here we have used the identity%
\[
z^{2}\boldsymbol{a}\left(  z\right)  u\left(  z\right)  =\boldsymbol{a}\left(
z\right)  z^{2}u\left(  z\right)  \text{.}%
\]
Our strategy is to modify (\ref{23}) so that the left hand side is an element
of $H\left(  D_{-}\right)  $ and the right hand side is an element of $\left(
e_{x}\boldsymbol{a}\right)  H_{3}\left(  D_{+}\right)  $. From (\ref{18})%
\[
z\partial_{x}w_{x}=s_{1}^{\prime}(x)+\sum_{2\leq k\leq L-2}z^{-k+1}%
s_{k}^{\prime}(x)+z^{-L+3}\widetilde{w}_{x}^{\left(  1\right)  }%
(z)=s_{1}^{\prime}\left(  x\right)  +v_{x}%
\]
follows with $v_{x}\in H\left(  D_{-}\right)  $, which combined with
(\ref{23}) yields%
\begin{align}
&  \partial_{x}^{2}w_{x}+2s_{1}^{\prime}w_{x}-2v_{x}\nonumber\\
&  =e_{x}\boldsymbol{a}\left(  \partial_{x}^{2}u_{x}-z^{2}u_{x}\right)
+2z\partial_{x}u_{x}+2s_{1}^{\prime}\left(  e_{x}\boldsymbol{a}u_{x}-1\right)
-2v_{x}\nonumber\\
&  =e_{x}\boldsymbol{a}\left(  \partial_{x}^{2}u_{x}-z^{2}u_{x}+2s_{1}%
^{\prime}u_{x}\right)  \text{.} \label{30}%
\end{align}
Note%
\[
\partial_{x}^{2}w_{x}+2s_{1}^{\prime}w_{x}-2v_{x}\in H\left(  D_{-}\right)
\text{, \ }\partial_{x}^{2}u_{x}-z^{2}u_{x}+2s_{1}^{\prime}u_{x}\in
H_{3}\left(  D_{+}\right)  \text{.}%
\]
Then, applying $\mathfrak{p}_{+}$ to (\ref{30}) we have%
\[
0=T\left(  e_{x}\boldsymbol{a}\right)  \left(  \partial_{x}^{2}u_{x}%
-z^{2}u_{x}+2s_{1}^{\prime}u_{x}\right)  \text{,}%
\]
and the invertibility of $T\left(  e_{x}\boldsymbol{a}\right)  $ on
$H_{3}\left(  D_{+}\right)  $ yields (\ref{22}).%
\[
\partial_{x}^{2}u_{x}-z^{2}u_{x}+2s_{1}^{\prime}u_{x}=0\text{.}%
\]
Since%
\begin{align*}
0  &  =\partial_{x}^{2}w_{x}+2s_{1}^{\prime}w_{x}-2v_{x}\\
&  =\sum_{1\leq k\leq L-3}\left(  s_{k}^{\prime\prime}+2s_{1}^{\prime}%
s_{k}-2s_{k+1}^{\prime}\right)  z^{-k}+z^{-L+3}\widetilde{v}_{x}%
\end{align*}
with $\widetilde{v}_{x}\in H\left(  D_{-}\right)  $, we have (\ref{21}%
).\medskip
\end{proof}

\begin{remark}
\label{r1}$q(x)$, $f_{\boldsymbol{a}}\left(  x,z\right)  $ themselves are
well-defined as continuous functions if $L\geq2$, and $\partial_{x}%
f_{\boldsymbol{a}}\left(  x,z\right)  $ exists as a continuous function.
Although the Schr\"{o}dinger equation (\ref{22}) is satisfied if $L\geq3$, it
seems that $\partial_{x}^{2}f_{\boldsymbol{a}}\left(  x,z\right)  $ exists
even if $L=2$ due to (\ref{22}). However we have no rigorous proof.
\end{remark}

Now the $m$-function $m_{\boldsymbol{a}}$ has another representation by
$f_{\boldsymbol{a}}\left(  x,z\right)  $.

\begin{corollary}
\label{c1}It holds that%
\begin{equation}
\left.  \dfrac{\partial_{x}f_{\boldsymbol{a}}\left(  x,z\right)
}{f_{\boldsymbol{a}}\left(  x,z\right)  }\right\vert _{x=0}=-m_{\boldsymbol{a}%
}\left(  z\right)  \text{.} \label{46}%
\end{equation}

\end{corollary}

\begin{proof}
Set $b(z)=\left.  \partial_{x}u_{x}\left(  z\right)  \right\vert _{x=0}$,
$\widetilde{b}(z)=\left.  \partial_{x}f_{\boldsymbol{a}}\left(  x,z\right)
\right\vert _{x=0}$. Then%
\[
b(z)\in H_{2}\left(  D_{+}\right)  \text{ \ and \ }\boldsymbol{a}\left(
z\right)  b(z)=\widetilde{b}(z)
\]
hold. Since $e^{xz}f_{\boldsymbol{a}}\left(  x,z\right)  =1+w_{x}\left(
z\right)  $, (\ref{18}) shows%
\begin{align*}
\widetilde{b}(z)  &  =-z\left(  1+w_{0}\left(  z\right)  \right)  +\sum_{1\leq
k\leq L-2}z^{-k}s_{k}^{\prime}(0)+z^{-L+2}\widetilde{w}_{0}^{\left(  1\right)
}(z)\\
&  =-z-s_{1}(0)+w\left(  z\right)
\end{align*}
with $w\in H\left(  D_{-}\right)  $. Applying $\mathfrak{p}_{+}$ to the
identity%
\[
\boldsymbol{a}\left(  z\right)  b(z)=-z-s_{1}(0)+w\left(  z\right)
\]
we have%
\[
T\left(  \boldsymbol{a}\right)  b(z)=-z-s_{1}(0)\text{,}%
\]
hence%
\[
b(z)=-T\left(  \boldsymbol{a}\right)  ^{-1}z-s_{1}(0)T\left(  \boldsymbol{a}%
\right)  ^{-1}1\text{,}%
\]
which implies%
\[
\widetilde{b}(z)=\boldsymbol{a}\left(  z\right)  b(z)=-\left(  z+\psi
_{\boldsymbol{a}}\left(  z\right)  \right)  -s_{1}(0)\left(  1+\varphi
_{\boldsymbol{a}}\left(  z\right)  \right)  \text{.}%
\]
Consequently we have%
\begin{align*}
\left.  \dfrac{\partial_{x}f_{\boldsymbol{a}}\left(  x,z\right)
}{f_{\boldsymbol{a}}\left(  x,z\right)  }\right\vert _{x=0}  &  =\dfrac
{-\left(  z+\psi_{\boldsymbol{a}}\left(  z\right)  \right)  -s_{1}(0)\left(
1+\varphi_{\boldsymbol{a}}\left(  z\right)  \right)  }{1+\varphi
_{\boldsymbol{a}}\left(  z\right)  }\\
&  =-m_{\boldsymbol{a}}\left(  z\right)  +\kappa_{1}\left(  \boldsymbol{a}%
\right)  -s_{1}(0)=-m_{\boldsymbol{a}}\left(  z\right)
\end{align*}
due to $s_{1}(0)=\kappa_{1}\left(  \boldsymbol{a}\right)  $.
\end{proof}

\subsection{Derivation of KdV equation}

Our next task is to derive the KdV equation from $g=e_{t,x}$ with%

\[
e_{t,x}(z)=e^{xz+tz^{3}}\text{.}%
\]
In this case the curve $C$ is determined by requesting $e_{t,x}(z)$ to be
bounded on $C$, hence%
\[
C=\left\{
\begin{array}
[c]{c}%
\pm\omega\left(  y\right)  +iy\text{, }y\in\mathbb{R}\text{; }\omega\text{\ is
smooth positive even on }\mathbb{R}\text{,}\\
\text{and }\omega\text{ satisfies }\omega\left(  y\right)  =O\left(
y^{-2}\right)  \text{ as }y\rightarrow\infty\text{.}%
\end{array}
\right\}  \text{.}%
\]
For $L\geq4$ let $\boldsymbol{a}\in\boldsymbol{A}_{L}\left(  C\right)  $ and
assume $e_{t,x}\boldsymbol{a}\in\boldsymbol{A}_{L}^{inv}\left(  C\right)  $
for any $t$, $x$ $\in\mathbb{R}$. Let $u_{t,x}=T\left(  e_{t,x}\boldsymbol{a}%
\right)  ^{-1}1\in H_{1}\left(  D_{+}\right)  $ and $w_{t,x}=\mathfrak{p}%
_{-}\left(  e_{t,x}\boldsymbol{a}f_{t,x}\right)  \in H\left(  D_{-}\right)  $.
Then, Lemma \ref{l5} implies%
\[
w_{t,x}\left(  z\right)  =\sum_{1\leq k\leq L-1}z^{-k}s_{k}(t,x)+z^{-L+1}%
\widetilde{w}_{t,x}(z)\text{ \ with }s_{k}\in\mathbb{C}\text{, }\widetilde
{w}_{t,x}\in H\left(  D_{-}\right)
\]
for $t$, $x\in\mathbb{R}$.

\begin{proposition}
\label{p2}Let $L\geq4$ and assume $e_{t,x}\boldsymbol{a}\in\boldsymbol{A}%
_{L}^{inv}\left(  C\right)  $ for any $t$, $x\in\mathbb{R}$. Set
\[
q(t,x)=-2\partial_{x}s_{1}\left(  t,x\right)  =-2\partial_{x}\kappa_{1}\left(
e_{t,x}\boldsymbol{a}\right)  \text{.}%
\]
Then $q$ satisfies the KdV equation%
\begin{equation}
\partial_{t}q(t,x)=\dfrac{1}{4}\partial_{x}^{3}q(t,x)-\dfrac{3}{2}%
q(t,x)\partial_{x}q(t,x)\text{.} \label{26}%
\end{equation}

\end{proposition}

\begin{proof}
In this case $N=1$, $n=3$. Our strategy for the proof is similar to that of
the last one. Since $w_{t,x}=e_{t,x}\boldsymbol{a}u_{t,x}-1\in H\left(
D_{-}\right)  $%
\[
\left\{
\begin{array}
[c]{l}%
\partial_{x}w_{t,x}=e_{t,x}\left(  z\boldsymbol{a}u_{t,x}+\boldsymbol{a}%
\partial_{x}u_{t,x}\right)  \smallskip\\
\partial_{t}w_{t,x}=e_{t,x}\left(  z^{3}\boldsymbol{a}u_{t,x}+\boldsymbol{a}%
\partial_{t}u_{t,x}\right)
\end{array}
\right.  \text{,}%
\]
which leads us to%
\[
\partial_{t}w_{t,x}=z^{2}\partial_{x}w_{t,x}+e_{t,x}\boldsymbol{a}\left(
\partial_{t}u_{t,x}-z^{2}\partial_{x}u_{t,x}\right)  \text{.}%
\]
Since%
\[
z^{2}\partial_{x}w_{t,x}\equiv zs_{1}^{\prime}+s_{2}^{\prime}+s_{3}^{\prime
}z^{-1}\text{ }\operatorname{mod}\text{ }z^{-1}H\left(  D_{-}\right)  \text{,}%
\]
substituting $1=e_{t,x}\boldsymbol{a}u_{t,x}-w_{t,x}$ and%
\begin{align*}
z  &  =ze_{t,x}\boldsymbol{a}u_{t,x}-zw_{t,x}\\
&  =\partial_{x}w_{t,x}-e_{t,x}\boldsymbol{a}\partial_{x}u_{t,x}-zw_{t,x}\\
&  \equiv s_{1}^{\prime}z^{-1}-s_{1}-s_{2}z^{-1}-e_{t,x}\boldsymbol{a}%
\partial_{x}u_{t,x}\\
&  \equiv\left(  s_{1}^{\prime}-s_{2}\right)  z^{-1}-s_{1}\left(
e_{t,x}\boldsymbol{a}u_{t,x}-w_{t,x}\right)  -e_{t,x}\boldsymbol{a}%
\partial_{x}u_{t,x}\\
&  \equiv\left(  s_{1}^{\prime}-s_{2}+s_{1}^{2}\right)  z^{-1}-e_{t,x}%
\boldsymbol{a}\left(  s_{1}u_{t,x}+\partial_{x}u_{t,x}\right)
\end{align*}
($\equiv$ means $\operatorname{mod}$ $z^{-1}H\left(  D_{-}\right)  $) into the
above identity yields%
\begin{align*}
&  z^{2}\partial_{x}w_{t,x}\\
&  \equiv s_{1}^{\prime}\left(  \left(  s_{1}^{\prime}-s_{2}+s_{1}^{2}\right)
z^{-1}-e_{t,x}\boldsymbol{a}\left(  s_{1}u_{t,x}+\partial_{x}u_{t,x}\right)
\right)  +s_{2}^{\prime}+s_{3}^{\prime}z^{-1}\\
&  =s_{2}^{\prime}\left(  e_{t,x}\boldsymbol{a}f_{t,x}-w_{t,x}\right)
+\left(  s_{1}^{\prime}\left(  s_{1}^{\prime}-s_{2}+s_{1}^{2}\right)
+s_{3}^{\prime}\right)  z^{-1}-e_{t,x}\boldsymbol{a}\left(  s_{1}%
u_{t,x}+\partial_{x}u_{t,x}\right) \\
&  =-s_{2}^{\prime}w_{t,x}+\left(  s_{1}^{\prime}\left(  s_{1}^{\prime}%
-s_{2}+s_{1}^{2}\right)  +s_{3}^{\prime}\right)  z^{-1}+e_{t,x}\boldsymbol{a}%
\left(  \left(  s_{2}^{\prime}-s_{1}\right)  u_{t,x}-\partial_{x}%
u_{t,x}\right)
\end{align*}
Therefore%
\begin{align*}
\partial_{t}w_{t,x}  &  \equiv-s_{2}^{\prime}w_{t,x}+\left(  s_{1}^{\prime
}\left(  s_{1}^{\prime}-s_{2}+s_{1}^{2}\right)  +s_{3}^{\prime}\right)
z^{-1}\\
&  +e_{t,x}\boldsymbol{a}\left(  \left(  s_{2}^{\prime}-s_{1}\right)
u_{t,x}+\partial_{t}u_{t,x}-z^{2}\partial_{x}u_{t,x}-\partial_{x}%
u_{t,x}\right)
\end{align*}
holds, and we have%
\begin{equation}
\left\{
\begin{array}
[c]{l}%
\left(  s_{2}^{\prime}-s_{1}\right)  u_{t,x}+\partial_{t}u_{t,x}-z^{2}%
\partial_{x}u_{t,x}-\partial_{x}u_{t,x}=0\smallskip\\
\partial_{t}w_{t,x}+s_{2}^{\prime}w_{t,x}-\left(  s_{1}^{\prime}\left(
s_{1}^{\prime}-s_{2}+s_{1}^{2}\right)  +s_{3}^{\prime}\right)  z^{-1}\equiv0
\end{array}
\right.  \text{.} \label{29}%
\end{equation}
due to the invertibility of $T\left(  e_{t,x}\boldsymbol{a}\right)  $ on
$H_{3}\left(  D_{+}\right)  $. Since the coefficient of $z^{-1}$ for the
second identity vanishes, it follows that%
\[
\partial_{t}s_{1}+s_{2}^{\prime}s_{1}-\left(  s_{1}^{\prime}\left(
s_{1}^{\prime}-s_{2}+s_{1}^{2}\right)  +s_{3}^{\prime}\right)  =0\text{.}%
\]
Here the identities for $k=1$, $2$ and $e^{tz^{3}}\boldsymbol{a}$ in
(\ref{21}) of Proposition \ref{p1}
\[
\left\{
\begin{array}
[c]{c}%
s_{1}^{\prime\prime}+2s_{1}^{\prime}s_{1}-2s_{2}^{\prime}=0\\
s_{2}^{\prime\prime}+2s_{1}^{\prime}s_{2}-2s_{3}^{\prime}=0
\end{array}
\right.
\]
allow us to have%
\begin{align*}
\partial_{t}s_{1}  &  =-s_{2}^{\prime}s_{1}+\left(  s_{1}^{\prime}\left(
s_{1}^{\prime}-s_{2}+s_{1}^{2}\right)  +s_{3}^{\prime}\right) \\
&  =-s_{1}\left(  s_{1}^{\prime\prime}/2+s_{1}^{\prime}s_{1}\right)  +\left(
s_{1}^{\prime}\right)  ^{2}+s_{1}^{2}s_{1}^{\prime}-s_{1}^{\prime}s_{2}%
+s_{2}^{\prime\prime}/2+s_{1}^{\prime}s_{2}\\
&  =-s_{1}\left(  s_{1}^{\prime\prime}/2+s_{1}^{\prime}s_{1}\right)  +\left(
s_{1}^{\prime}\right)  ^{2}+s_{1}^{2}s_{1}^{\prime}+\left(  s_{1}%
^{\prime\prime}/2+s_{1}^{\prime}s_{1}\right)  ^{\prime}/2\\
&  =s_{1}^{\prime\prime\prime}/4+3\left(  s_{1}^{\prime}\right)
^{2}/2\text{,}%
\end{align*}
which is (\ref{26}) by substituting $q(t,x)=-2\partial_{x}s_{1}\left(
t,x\right)  $.\medskip
\end{proof}

The above calculation does not reveal explicitly the reason why the KdV
equation appears. There is a hidden algebraic structure behind discovered by
Sato \cite{sa}.

The invertibility of $T\left(  e_{t,x}\boldsymbol{a}\right)  $ plays an
essential role in the above calculation. The non-existence of $T\left(
e_{t,x}\boldsymbol{a}\right)  ^{-1}$ at some point $\left(  t,x\right)  $
means that the solution $q(t,x)$ has a singularity at $\left(  t,x\right)  $.

\section{Tau-function}

Hirota introduced an object $\tau_{\boldsymbol{a}}$ called tau-function whose
mathematical meaning was discovered by Sato later. In this paper the
tau-function will be used to prove the invertibility of $T\left(
g\boldsymbol{a}\right)  $.

The tau-function defined in \cite{k2} is written in the present context as%
\begin{equation}
\tau_{\boldsymbol{a}}\left(  g\right)  =\det\left(  g^{-1}T\left(
g\boldsymbol{a}\right)  T(\boldsymbol{a})^{-1}\right)  \text{.} \label{31}%
\end{equation}
The operator $g^{-1}T\left(  g\boldsymbol{a}\right)  T(\boldsymbol{a})^{-1}$
is a map on $H_{N}\left(  D_{+}\right)  $ and the determinant is well-defined
if the operator $g^{-1}T\left(  g\boldsymbol{a}\right)  T(\boldsymbol{a}%
)^{-1}-I$ is of trace class on $H_{N}\left(  D_{+}\right)  $. The identity
(\ref{40}) implies%
\[
g^{-1}T\left(  g\boldsymbol{a}\right)  T\left(  \boldsymbol{a}\right)
^{-1}-I=g^{-1}R_{\boldsymbol{a}}\left(  1,g\right)  T\left(  \boldsymbol{a}%
\right)  ^{-1}\text{,}%
\]
hence it is sufficient for this that $R_{\boldsymbol{a}}\left(  1,g\right)  $
is of trace class. Since $R_{\boldsymbol{a}}\left(  1,g\right)  $ is a product
of $H_{g}$ and $S_{\boldsymbol{a}}$, there are two cases where
$R_{\boldsymbol{a}}\left(  1,g\right)  $ is of trace class.\medskip\newline(i)
$\ H_{g}$ is of trace class.\newline(ii) $H_{g}$ and $S_{\boldsymbol{a}}$ are
of Hilbert-Schmidt class.\medskip\newline(i) is the case for $g\in\Gamma
_{0}^{\left(  0\right)  }$ (rational function of order $0$) due to Lemma
\ref{l13}. However for other cases one has to impose an extra condition on
$\boldsymbol{a}$. To avoid this inconvenience we use the modified determinant
$\det_{2}$, namely%
\[
\det\nolimits_{2}\left(  I+A\right)  =\det(I+A)e^{-\text{\textrm{tr}}%
A}\text{.}%
\]
It is known that this determinant can be extended to any operator $A$ of
Hilbert-Schmidt (HS in short) class. Since $I+A$ is invertible if and only if
$\det\nolimits_{2}\left(  I+A\right)  \neq0$, this $\det\nolimits_{2}$ is
sufficient to verify the existence of $T\left(  g\boldsymbol{a}\right)  ^{-1}%
$. Set%
\begin{equation}
\tau_{\boldsymbol{a}}^{\left(  2\right)  }\left(  g\right)  =\det
\nolimits_{2}\left(  g^{-1}T\left(  g\boldsymbol{a}\right)  T(\boldsymbol{a}%
)^{-1}\right)  \text{ for }\boldsymbol{a}\in\boldsymbol{A}_{L}^{inv}\left(
C\right)  \text{, }g\in\Gamma_{n}^{\left(  m\right)  }\text{.} \label{5}%
\end{equation}
$\tau_{\boldsymbol{a}}^{\left(  2\right)  }\left(  g\right)  $ can be defined
if $R_{\boldsymbol{a}}\left(  1,g\right)  =H_{g}S_{a}$ is of HS class, which
is valid if $H_{g}$ is of HS-class as a map from $H\left(  D_{-}\right)  $ to
$H_{N+m}\left(  D_{+}\right)  $, and Lemma \ref{l10} implies that $H_{g}$ is
of HS class if%
\begin{equation}
N\geq\max\left\{  n,1-m\right\}  \text{,} \label{6}%
\end{equation}
which implies%
\begin{equation}
L\geq\max\left\{  n,1-m\right\}  \text{.} \label{19}%
\end{equation}
Conversely the existence of $N$ satisfying (\ref{6}) follows from (\ref{19}).

In the definitions of $\tau_{\boldsymbol{a}}\left(  g\right)  $,
$\tau_{\boldsymbol{a}}^{\left(  2\right)  }\left(  g\right)  $ the operator
$g^{-1}T\left(  g\boldsymbol{a}\right)  T(\boldsymbol{a})^{-1}$ is a map on
$H_{N}\left(  D_{+}\right)  $, hence $\tau_{\boldsymbol{a}}\left(  g\right)
$, $\tau_{\boldsymbol{a}}^{\left(  2\right)  }\left(  g\right)  $ may depend
on $N$. However we have

\begin{lemma}
\label{l8}For $\boldsymbol{a}\in\boldsymbol{A}_{L}\left(  C\right)  $,
$g\in\Gamma_{n}^{\left(  m\right)  }$ assume that $N$, $N^{\prime}%
\in\mathbb{Z}_{+}$ satisfy (\ref{6}) and the two operators%
\[
\left\{
\begin{array}
[c]{l}%
T\left(  \boldsymbol{a}\right)  :H_{N}\left(  D_{+}\right)  \rightarrow
H_{N}\left(  D_{+}\right)  \smallskip\\
T\left(  \boldsymbol{a}\right)  :H_{N^{\prime}}\left(  D_{+}\right)
\rightarrow H_{N^{\prime}}\left(  D_{+}\right)
\end{array}
\right.
\]
are bijective. Then the determinants and the modified determinants on
$H_{N}\left(  D_{+}\right)  $ and $H_{N^{\prime}}\left(  D_{+}\right)  $ in
(\ref{31}) and (\ref{5}) are equal, hence $\tau_{\boldsymbol{a}}\left(
g\right)  $, $\tau_{\boldsymbol{a}}^{\left(  2\right)  }\left(  g\right)  $ do
not depend on $N$.
\end{lemma}

To show this lemma we use a metric free nature of determinant. For the
necessary facts of determinant refer to \cite{si}.

\begin{lemma}
\label{l9}Let $H$, $H_{1}$ be two Hilbert spaces and $H_{1}$ be a subspace of
$H$ as vector spaces. Assume $H_{1}$ is dense. Suppose a linear operator $A$
on $H$ is of Hilbert-Schmidt class and satisfies%
\[
AH_{1}\subset H_{1}\text{ \ and }A_{H_{1}}\equiv\left.  A\right\vert _{H_{1}%
}\text{ is of Hilbert-Schmidt class in }H_{1}\text{.}%
\]
Then%
\[
\det\nolimits_{2}\left(  I_{H}+A\right)  =\det\nolimits_{2}\left(  I_{H_{1}%
}+A_{H_{1}}\right)
\]
holds. If $A$ is of trace class, then $\det\left(  I_{H}+A\right)
=\det\left(  I_{H_{1}}+A_{H_{1}}\right)  $ holds as well.
\end{lemma}

\begin{proof}
Let $\left\{  e_{n}\right\}  _{n\geq1}$ be a complete orthonormal basis of
$H_{1}$ and $\left\{  f_{n}\right\}  _{n\geq1}$ be the orthonormal vectors in
$H$ generated from $\left\{  e_{n}\right\}  _{n\geq1}$ by the Gram-Schmidt
process. Since $H_{1}$ is dense in $H$, $\left\{  f_{n}\right\}  _{n\geq1}$
turns to be complete in $H$. Let $V_{n}$ be the $n$-dimensional subspace
generated by $\left\{  e_{k}\right\}  _{1\leq k\leq n}$ \ and $P_{n}$, $Q_{n}$
be the orthogonal projections to $V_{n}$ from $H_{1}$, $H$ respectively. Then
it is known that%
\begin{equation}
\left\{
\begin{array}
[c]{l}%
\det\nolimits_{2}\left(  I_{H_{1}}+P_{n}A_{H_{1}}P_{n}\right)  \rightarrow
\det\nolimits_{2}(I_{H_{1}}+A_{H_{1}})\smallskip\\
\det\nolimits_{2}\left(  I_{H}+Q_{n}AQ_{n}\right)  \rightarrow\det
\nolimits_{2}(I_{H}+A)
\end{array}
\right.  \text{ as }n\rightarrow\infty\text{.} \label{45}%
\end{equation}
Let $B$, $E$ be $n\times n$ matrices whose entries are%
\[
f_{i}=\sum_{1\leq j\leq n}b_{ij}e_{j}\rightarrow B=\left(  b_{ij}\right)
\text{, and\ }E=\left(  \left(  e_{i},e_{j}\right)  _{H}\right)  \text{.}%
\]
If we denote%
\[
A_{n}^{H_{1}}=\left(  \left(  Ae_{i},e_{j}\right)  _{H_{1}}\right)  _{1\leq
i,j\geq n}\text{, \ \ \ }A_{n}^{H}=\left(  \left(  Af_{i},f_{j}\right)
_{H}\right)  _{1\leq i,j\leq n}\text{,}%
\]
then the identity $A_{n}^{H}=BA_{n}^{H_{1}}EB^{\ast}$ yields%
\begin{align*}
\det\nolimits_{2}\left(  I_{H}+Q_{n}AQ_{n}\right)   &  =\det\nolimits_{2}%
\left(  I_{n}+A_{n}^{H}\right) \\
&  =\det\nolimits_{2}\left(  I_{n}+BA_{n}^{H_{1}}EB^{\ast}\right) \\
&  =\det\nolimits_{2}\left(  I_{n}+BA_{n}^{H_{1}}B^{-1}\right) \\
&  =\det\nolimits_{2}\left(  I_{n}+A_{n}^{H_{1}}\right)  =\det\nolimits_{2}%
\left(  I_{H_{1}}+P_{n}A_{H_{1}}P_{n}\right)
\end{align*}
due to $BEB^{\ast}=I_{n}$. This together with (\ref{45}) completes the proof
for $\det_{2}$. The proof for $\det$ is similar.\medskip
\end{proof}

\textbf{Proof of Lemma \ref{l8}}. Note that for $N^{\prime}\geq N$ an identity%
\begin{equation}
\left.  g^{-1}R_{\boldsymbol{a}}^{\left(  N^{\prime}\right)  }\left(
1,g\right)  T_{N^{\prime}}\left(  \boldsymbol{a}\right)  ^{-1}\right\vert
_{H_{N}\left(  D_{+}\right)  }=g^{-1}R_{\boldsymbol{a}}^{\left(  N\right)
}\left(  1,g\right)  T_{N}\left(  \boldsymbol{a}\right)  ^{-1} \label{38}%
\end{equation}
holds on $H_{N}\left(  D_{+}\right)  $. On the other hand, for $u\in H\left(
D_{+}\right)  $ the identity%
\[
\left(  b-z\right)  ^{N}u\left(  z\right)  =\lim_{\epsilon\rightarrow0}\left(
b-z\right)  ^{N}\left(  1-\epsilon z\right)  ^{-N}u\left(  z\right)
\]
implies $H\left(  D_{+}\right)  $ is dense in $H_{N}\left(  D_{+}\right)  $.
Then, applying Lemma \ref{l9} to
\[
H=H_{N^{\prime}}\left(  D_{+}\right)  \text{, }H_{1}=H_{N}\left(
D_{+}\right)  \text{, }A=g^{-1}R_{\boldsymbol{a}}^{\left(  N^{\prime}\right)
}\left(  1,g\right)  T_{N^{\prime}}\left(  \boldsymbol{a}\right)  ^{-1}%
\]
yields the lemma.\medskip

This Lemma yields flexibility in choosing $N$, namely $N$ can be arbitrary if
(\ref{6}) is satisfied for given $L$, $m$, $n$. Therefore $\tau
_{\boldsymbol{a}}^{\left(  2\right)  }\left(  g\right)  $ can be defined for
$\boldsymbol{a}\in\boldsymbol{A}_{L}^{inv}\left(  C\right)  $, $g\in\Gamma
_{n}^{\left(  m\right)  }$ under the condition (\ref{19}). However it should
be noted that for any rational function $r\in\Gamma_{0}^{\left(  m\right)  }$
Lemma \ref{l13} shows $r^{-1}T\left(  r\boldsymbol{a}\right)  T(\boldsymbol{a}%
)^{-1}-I$ is of finite rank on any space $H_{N}\left(  D_{+}\right)  $ with
$N$ such that $-m\leq N\leq L$. Therefore, $\tau_{\boldsymbol{a}}^{\left(
2\right)  }\left(  r\right)  $, $\tau_{\boldsymbol{a}}\left(  r\right)  $ can
be defined for $r\in\Gamma_{0}^{\left(  m\right)  }$, $\boldsymbol{a}%
\in\boldsymbol{A}_{L}^{inv}\left(  C\right)  $ if $L\geq-m$.

\subsection{Cocycle property of tau-function}

The tau-function is a key material to study the KdV flow and in this section
we give fundamental properties for the tau-function.

Note%
\[
\boldsymbol{a}\in\boldsymbol{A}_{L}\left(  C\right)  \text{, \ }g\in\Gamma
_{n}^{\left(  m\right)  }\Longrightarrow g\boldsymbol{a}\in\boldsymbol{A}%
_{L}\left(  C\right)  \text{,}%
\]
since $g$ is analytic and bounded on $D_{+}$. Assume further%
\[
\boldsymbol{a}\in\boldsymbol{A}_{L}^{inv}\left(  C\right)  \text{, \ }g_{1}%
\in\Gamma_{n}^{\left(  0\right)  }\text{, \ }g_{2}\in\Gamma_{n}^{\left(
m\right)  }\text{.}%
\]
We consider three tau-functions $\tau_{\boldsymbol{a}}^{\left(  2\right)
}\left(  g_{1}g_{2}\right)  $, $\tau_{\boldsymbol{a}}^{\left(  2\right)
}\left(  g_{1}\right)  $, $\tau_{g_{1}\boldsymbol{a}}^{\left(  2\right)
}\left(  g_{2}\right)  $ simultaneously, which is possible if $L\in
\mathbb{Z}_{+}$ satisfies (\ref{9}) and $g_{1}\boldsymbol{a}\in\boldsymbol{A}%
_{L}^{inv}\left(  C\right)  $.

For simplicity of notations set%
\begin{align}
E_{\boldsymbol{a}}\left(  g_{1},g_{2}\right)   &  =\text{\textrm{tr}}\left(
\left(  \left(  g_{1}g_{2}\right)  ^{-1}T\left(  g_{2}g_{1}\boldsymbol{a}%
\right)  T(g_{1}\boldsymbol{a})^{-1}g_{1}-I\right)  \left(  g_{1}^{-1}%
T(g_{1}\boldsymbol{a})T(\boldsymbol{a})^{-1}-I\right)  \right) \nonumber\\
&  =\text{\textrm{tr}}\left(  \left(  g_{1}g_{2}\right)  ^{-1}%
R_{\boldsymbol{a}}\left(  g_{1},g_{2}\right)  T(g_{1}\boldsymbol{a}%
)^{-1}R_{\boldsymbol{a}}\left(  1,g_{1}\right)  T(\boldsymbol{a})^{-1}\right)
\text{.} \label{7}%
\end{align}

\begin{lemma}
\label{l7}Assume $L$ satisfies (\ref{19}) and let $N$ be one $N$ of (\ref{6}).
Then we have the followings.\newline(i) \ The map $T\left(  g_{1}%
\boldsymbol{a}\right)  $ is bijective on $H_{N}\left(  D_{+}\right)  $ if and
only if $\tau_{\boldsymbol{a}}^{\left(  2\right)  }\left(  g_{1}\right)
\neq0$. Similarly the map%
\[
T\left(  g_{2}g_{1}\boldsymbol{a}\right)  :H_{N}\left(  D_{+}\right)
\rightarrow H_{N+m}\left(  D_{+}\right)
\]
is bijective if and only if $\tau_{\boldsymbol{a}}^{\left(  2\right)  }\left(
g_{1}g_{2}\right)  \neq0$.\newline(ii) If $\tau_{\boldsymbol{a}}^{\left(
2\right)  }\left(  g_{1}\right)  \neq0$, then it holds that%
\begin{equation}
\tau_{\boldsymbol{a}}^{\left(  2\right)  }\left(  g_{1}g_{2}\right)
=\tau_{\boldsymbol{a}}^{\left(  2\right)  }\left(  g_{1}\right)  \tau
_{g_{1}\boldsymbol{a}}^{\left(  2\right)  }\left(  g_{2}\right)  \exp\left(
-E_{\boldsymbol{a}}\left(  g_{1},g_{2}\right)  \right)  \text{.} \label{17}%
\end{equation}
Additionally if $r_{1}\in\Gamma_{0}^{\left(  0\right)  }$, $r_{2}\in\Gamma
_{0}^{\left(  m\right)  }$ and $\tau_{\boldsymbol{a}}\left(  r_{1}\right)
\neq0$, then%
\begin{equation}
\tau_{\boldsymbol{a}}\left(  r_{1}r_{2}\right)  =\tau_{\boldsymbol{a}}\left(
r_{1}\right)  \tau_{r_{1}\boldsymbol{a}}\left(  r_{2}\right)  \text{.}
\label{12}%
\end{equation}
(iii) Suppose $g_{1}$ satisfies $g_{1}\left(  z\right)  =g_{1}\left(
-z\right)  $. Then, it holds that $T\left(  g_{2}g_{1}\boldsymbol{a}\right)
=T\left(  g_{2}\boldsymbol{a}\right)  g_{1}$ and%
\[
\tau_{g_{1}\boldsymbol{a}}^{\left(  2\right)  }\left(  g_{2}\right)
=\tau_{\boldsymbol{a}}^{\left(  2\right)  }\left(  g_{2}\right)  \text{.}%
\]
Similarly, for rational functions $r_{1}\in\Gamma_{0}^{\left(  0\right)  }$,
$r_{2}\in\Gamma_{0}^{\left(  m\right)  }$ satisfying $r_{1}\left(  z\right)
=r_{1}\left(  -z\right)  $ we have%
\[
\tau_{r_{1}\boldsymbol{a}}\left(  r_{2}\right)  =\tau_{\boldsymbol{a}}\left(
r_{2}\right)  \text{, \ \ }\tau_{\boldsymbol{a}}\left(  r_{1}r_{2}\right)
=\tau_{\boldsymbol{a}}\left(  r_{1}\right)  \tau_{\boldsymbol{a}}\left(
r_{2}\right)  \text{.}%
\]

\end{lemma}

\begin{proof}
The identity (\ref{7}) implies%
\[
g_{2}^{-1}T\left(  g_{2}g_{1}\boldsymbol{a}\right)  T\left(  g_{1}%
\boldsymbol{a}\right)  ^{-1}=I+g_{2}^{-1}R_{\boldsymbol{a}}\left(  g_{1}%
,g_{2}\right)  T\left(  g_{1}\boldsymbol{a}\right)  ^{-1}%
\]
with the Hilbert-Schmidt class operator $R_{\boldsymbol{a}}\left(  g_{1}%
,g_{2}\right)  $. Then general theory of Fredholm determinant shows that the
operator $g_{2}^{-1}T\left(  g_{2}g_{1}\boldsymbol{a}\right)  T\left(
g_{1}\boldsymbol{a}\right)  ^{-1}$ is bijective if and only if $\det
_{2}\left(  g_{2}^{-1}T\left(  g_{2}g_{1}\boldsymbol{a}\right)  T\left(
g_{1}\boldsymbol{a}\right)  ^{-1}\right)  \neq0$, which implies the
bijectivity of $T\left(  g_{2}g_{1}\boldsymbol{a}\right)  $. The bijectivity
of $T\left(  g_{1}\boldsymbol{a}\right)  $ follows by letting $g_{2}=1$.

The definition of the tau-function says%
\[
\tau_{\boldsymbol{a}}^{\left(  2\right)  }\left(  g_{1}g_{2}\right)
=\det\nolimits_{2}\left(  (g_{1}g_{2})^{-1}T\left(  g_{1}g_{2}\boldsymbol{a}%
\right)  T(\boldsymbol{a})^{-1}\right)  \text{.}%
\]
On the other hand it holds that%
\begin{align}
&  (g_{1}g_{2})^{-1}T\left(  g_{1}g_{2}\boldsymbol{a}\right)  T(\boldsymbol{a}%
)^{-1}\nonumber\\
&  =\left(  g_{1}^{-1}\left(  g_{2}^{-1}T\left(  g_{2}g_{1}\boldsymbol{a}%
\right)  T(g_{1}\boldsymbol{a})^{-1}\right)  g_{1}\right)  \left(  g_{1}%
^{-1}T(g_{1}\boldsymbol{a})T(\boldsymbol{a})^{-1}\right)  \text{.} \label{25}%
\end{align}
Note the identity%
\[
\left\{
\begin{array}
[c]{l}%
\det\nolimits_{2}\left(  G^{-1}\left(  I+A\right)  G\right)  =\det
\nolimits_{2}\left(  I+A\right) \\
\det\nolimits_{2}\left(  \left(  I+A\right)  \left(  I+B\right)  \right)
=\det\nolimits_{2}\left(  I+A\right)  \det\nolimits_{2}\left(  I+B\right)
e^{-\text{\textrm{tr}}\left(  AB\right)  }%
\end{array}
\right.
\]
for Hilbert-Schmidt operators $A$, $B$ and a bounded operator $G$ having
bounded inverse. Then taking determinant in (\ref{25}) yields%
\begin{align*}
\tau_{\boldsymbol{a}}^{\left(  2\right)  }\left(  g_{1}g_{2}\right)   &
=\det\nolimits_{2}\left(  g_{1}^{-1}g_{2}^{-1}T\left(  g_{1}g_{2}%
\boldsymbol{a}\right)  T(g_{2}\boldsymbol{a})^{-1}g_{1}\right)  \det
\nolimits_{2}\left(  g_{1}^{-1}T(g_{1}\boldsymbol{a})T(\boldsymbol{a}%
)^{-1}\right) \\
&  \times\exp\left(  -\text{\textrm{tr}}\left(  g_{1}^{-1}g_{2}^{-1}%
R_{\boldsymbol{a}}\left(  g_{1},g_{2}\right)  T\left(  g_{1}\boldsymbol{a}%
\right)  ^{-1}g_{1}g_{1}^{-1}R_{\boldsymbol{a}}\left(  1,g_{1}\right)
T\left(  \boldsymbol{a}\right)  ^{-1}\right)  \right) \\
&  =\tau_{\boldsymbol{a}}^{\left(  2\right)  }\left(  g_{1}\right)
\tau_{g_{1}\boldsymbol{a}}^{\left(  2\right)  }\left(  g_{2}\right)
\exp\left(  -E_{\boldsymbol{a}}\left(  g_{1},g_{2}\right)  \right)
\end{align*}
For $r\in\Gamma_{0}^{\left(  0\right)  }$ Lemma \ref{l13} implies that
$r^{-1}T\left(  r\boldsymbol{a}\right)  T(\boldsymbol{a})^{-1}$ is of finite
rank, hence $\tau_{\boldsymbol{a}}\left(  r\right)  $ can be defined. Taking
the determinant in (\ref{25}) we easily have (\ref{12}).

Suppose $g_{1}(z)=g_{1}(-z)$. For $u\in H_{N}\left(  D_{+}\right)  $,
$\boldsymbol{a}=\left(  a_{1},a_{2}\right)  $ we have%
\begin{align*}
T\left(  g_{2}g_{1}\boldsymbol{a}\right)  u  &  =\mathfrak{p}_{+}\left(
g_{1}g_{2}a_{1}\right)  u_{e}+\mathfrak{p}_{+}\left(  g_{1}g_{2}a_{2}\right)
u_{o}\\
&  =\mathfrak{p}_{+}\left(  g_{2}a_{1}\right)  g_{1}u_{e}+\mathfrak{p}%
_{+}\left(  g_{2}a_{2}\right)  g_{1}u_{o}\\
&  =\mathfrak{p}_{+}\left(  g_{2}a_{1}\right)  \left(  g_{1}u\right)
_{e}+\mathfrak{p}_{+}\left(  g_{2}a_{2}\right)  \left(  g_{1}u\right)  _{o}\\
&  =T\left(  g_{2}\boldsymbol{a}\right)  \left(  g_{1}u\right)  \text{,}%
\end{align*}
which implies $T\left(  g_{2}g_{1}\boldsymbol{a}\right)  =T\left(
g_{2}\boldsymbol{a}\right)  g_{1}$ on $H_{N}\left(  D_{+}\right)  $, hence%
\[
T\left(  g_{2}g_{1}\boldsymbol{a}\right)  T(g_{1}\boldsymbol{a})^{-1}=T\left(
g_{2}\boldsymbol{a}\right)  g_{1}g_{1}^{-1}T(\boldsymbol{a})^{-1}=T\left(
g_{2}\boldsymbol{a}\right)  T(\boldsymbol{a})^{-1}%
\]
is valid. This shows%
\[
\tau_{g_{1}\boldsymbol{a}}^{\left(  2\right)  }\left(  g_{2}\right)
=\det\nolimits_{2}\left(  g_{2}^{-1}T\left(  g_{2}g_{1}\boldsymbol{a}\right)
T(g_{1}\boldsymbol{a})^{-1}\right)  =\det\nolimits_{2}\left(  g_{2}%
^{-1}T\left(  g_{2}\boldsymbol{a}\right)  T(\boldsymbol{a})^{-1}\right)
=\tau_{\boldsymbol{a}}^{\left(  2\right)  }\left(  g_{2}\right)  \text{.}%
\]
\medskip
\end{proof}

\subsection{Continuity of tau-functions}

Since the determinant $\det_{2}$ is estimated by the HS-norm, the continuity
of $\tau_{\boldsymbol{a}}^{\left(  2\right)  }\left(  g\right)  $ with respect
to $g\in\Gamma_{n}^{\left(  m\right)  }$ follows from that of $g^{-1}%
R_{\boldsymbol{a}}\left(  1,g\right)  $ with respect to the HS-norm, which is
reduced to the continuity of $H_{g}=\mathfrak{p}_{+}\left(  g\cdot\right)  $
on $H\left(  D_{-}\right)  $ with respect to the HS-norm due to%
\begin{equation}
R_{\boldsymbol{a}}\left(  1,g\right)  =H_{g}S_{\boldsymbol{a}}\text{ \ (see
Lemma \ref{l3}),} \label{55}%
\end{equation}
where%
\[
\left\{
\begin{array}
[c]{ll}%
R_{\boldsymbol{a}}\left(  1,g\right)  & :H_{N}\left(  D_{+}\right)
\rightarrow H_{N}\left(  D_{+}\right) \\
H_{g} & :H\left(  D_{-}\right)  \rightarrow H_{N}\left(  D_{+}\right) \\
S_{\boldsymbol{a}} & :H_{N}\left(  D_{+}\right)  \rightarrow H\left(
D_{-}\right)
\end{array}
\right.  \text{.}%
\]
The condition for $L$, $n$ is (\ref{19}), namely%
\[
L\geq\max\left\{  n,1-m\right\}
\]
and $N$ is arbitrary if it satisfies (\ref{6}). Denote%
\begin{align*}
&  d_{N}\left(  g_{1},g_{2}\right) \\
&  =\left(
\begin{array}
[c]{c}%
{\displaystyle\int_{C^{2}}}
\left\vert \dfrac{g_{1}\left(  z\right)  ^{-1}g_{1}\left(  \lambda\right)
-g_{2}\left(  z\right)  ^{-1}g_{2}\left(  \lambda\right)  }{z-\lambda
}\right\vert ^{2}\left\vert z\right\vert ^{-2N}\left\vert dz\right\vert
\left\vert d\lambda\right\vert
\end{array}
\right)  ^{1/2}%
\end{align*}
if $g_{1},g_{2}\in\Gamma_{n}^{\left(  m\right)  }$, we have

\begin{lemma}
\label{l14}Let $\boldsymbol{a}\in\boldsymbol{A}_{L}^{inv}\left(  C\right)  $
with $L$ satisfying (\ref{19}). Assume $g_{1}$, $g_{2}\in\Gamma_{n}^{\left(
m\right)  }$ and $d_{N}\left(  g_{1},1\right)  \leq c_{1}$. Then there exists
a constant $c_{\boldsymbol{a}}$ depending on $c_{1}$, $\boldsymbol{a}$, $N$
such that%
\begin{equation}
\left\vert \tau_{\boldsymbol{a}}^{\left(  2\right)  }\left(  g_{1}\right)
-\tau_{\boldsymbol{a}}^{\left(  2\right)  }\left(  g_{2}\right)  \right\vert
\leq c_{\boldsymbol{a}}d_{N}\left(  g_{1},g_{2}\right)  \text{.} \label{50}%
\end{equation}

\end{lemma}

\begin{proof}
Recall the definition $\tau_{\boldsymbol{a}}^{\left(  2\right)  }\left(
g\right)  =\det_{2}\left(  I+g^{-1}R_{\boldsymbol{a}}\left(  1,g\right)
T\left(  \boldsymbol{a}\right)  ^{-1}\right)  $ for $g\in\Gamma_{n}^{\left(
m\right)  }$, and%
\[
g^{-1}R_{\boldsymbol{a}}\left(  1,g\right)  T\left(  \boldsymbol{a}\right)
^{-1}=g^{-1}H_{g}S_{\boldsymbol{a}}T\left(  \boldsymbol{a}\right)
^{-1}\text{.}%
\]
The HS-norm of this operator is dominated by $d_{N}\left(  g,1\right)  $,
hence if $d_{N}\left(  g,1\right)  <\infty$, then $\tau_{\boldsymbol{a}%
}^{\left(  2\right)  }\left(  g\right)  $ is defined finitely. Generally if
$\left\Vert A\right\Vert _{HS}$, $\left\Vert B\right\Vert _{HS}\leq c_{1}$
there exists a constant $c_{2}$ depending only on $c_{1}$ such that%
\[
\left\vert \det\nolimits_{2}\left(  I+A\right)  -\det\nolimits_{2}\left(
I+B\right)  \right\vert \leq c_{2}\left\Vert A-B\right\Vert _{HS}\text{.}%
\]
Therefore $\tau_{\boldsymbol{a}}^{\left(  2\right)  }\left(  g_{1}\right)
-\tau_{\boldsymbol{a}}^{\left(  2\right)  }\left(  g_{2}\right)  $ can be
estimated by those of the HS-norms of $g_{1}^{-1}H_{g_{1}}-g_{2}^{-1}H_{g_{2}%
}$ on the space $H_{N}\left(  D_{+}\right)  $, which is just equal to
$d_{N}\left(  g_{1},g_{2}\right)  $.\medskip
\end{proof}

For later purpose we give a sufficient condition for the convergence of
$\tau_{\boldsymbol{a}}^{\left(  2\right)  }\left(  g_{k}\right)  $.

\begin{lemma}
\label{l11}Assume the following properties for $g_{k}$, $g\in\Gamma
_{n}^{\left(  0\right)  }$:%
\begin{equation}
\left\{
\begin{array}
[c]{l}%
\text{(i) there exist }c_{1}\text{, }c_{2}>0\text{ such that for }z\in U\\
\text{ \ \ \ \ \ \ \ \ \ }c_{1}\leq\left\vert g_{k}\left(  z\right)
\right\vert \leq c_{2}\text{, \ \ }\left\vert g_{k}^{\prime}\left(  z\right)
\right\vert \leq c_{2}\left\vert z\right\vert ^{n-1}\\
\text{(ii) }g_{k}\left(  z\right)  \rightarrow g\left(  z\right)  \text{ as
}k\rightarrow\infty\text{ for any }z\in C
\end{array}
\right.  \text{,} \label{61}%
\end{equation}
where $U$ is a neighborhood of the closure of $D_{+}$ satisfying (\ref{34}).
Then, for $r\in\Gamma_{0}^{\left(  m\right)  }$ and $\boldsymbol{a}%
\in\boldsymbol{A}_{L}\left(  C\right)  $ with
\[
L\geq\max\left\{  n,1-m\right\}  \text{,}%
\]
it holds that%
\[
\tau_{\boldsymbol{a}}^{\left(  2\right)  }\left(  rg_{k}\right)
\rightarrow\tau_{\boldsymbol{a}}^{\left(  2\right)  }\left(  rg\right)
\text{.}%
\]

\end{lemma}

\begin{proof}
Choose an integer $N\geq0$ such that $L\geq N\geq\max\left\{  n,1-m\right\}
$. Set%
\[
\Delta_{k}\left(  z,\lambda\right)  =\dfrac{g\left(  z\right)  ^{-1}g\left(
\lambda\right)  -g_{k}\left(  z\right)  ^{-1}g_{k}\left(  \lambda\right)
}{z-\lambda}\text{.}%
\]
Since%
\[
d_{L,N}\left(  rg_{k},rg\right)  ^{2}=%
{\displaystyle\int_{C^{2}}}
\left\vert \Delta_{k}\left(  z,\lambda\right)  \right\vert ^{2}\left\vert
z\right\vert ^{-2N}\left\vert dz\right\vert \left\vert d\lambda\right\vert
\text{,}%
\]
it is sufficient for $d_{L,N}\left(  rg_{k},rg\right)  \rightarrow0$ as
$k\rightarrow\infty$ to show that there exists a function $f$ integrable with
respect to $\left\vert z\right\vert ^{-2N}\left\vert dz\right\vert \left\vert
d\lambda\right\vert $ such that%
\[
\left\vert \Delta_{k}\left(  z,\lambda\right)  \right\vert ^{2}\leq f\left(
z,\lambda\right)  \text{.}%
\]
Note%
\begin{align*}
\Delta_{k}\left(  z,\lambda\right)   &  =r\left(  z\right)  ^{-1}g\left(
z\right)  ^{-1}\dfrac{r\left(  \lambda\right)  g\left(  \lambda\right)
-r\left(  z\right)  g\left(  z\right)  }{z-\lambda}\\
&  \text{ \ \ \ \ \ }-r\left(  z\right)  ^{-1}g_{k}\left(  z\right)
^{-1}\dfrac{r\left(  \lambda\right)  g_{k}\left(  \lambda\right)  -r\left(
z\right)  g_{k}\left(  z\right)  }{z-\lambda}\text{.}%
\end{align*}
Then (\ref{64}) implies%
\[
\left\vert \dfrac{r\left(  z\right)  g_{k}\left(  z\right)  -r\left(
\lambda\right)  g_{k}\left(  \lambda\right)  }{z-\lambda}\right\vert \leq
c_{1}\left\{
\begin{array}
[c]{ll}%
\left\vert z\right\vert ^{m+n-1} & \text{if }\left\vert z-\lambda\right\vert
\leq\epsilon\left\vert \lambda\right\vert \\
\left\vert z\right\vert ^{-1}\left(  \left\vert z\right\vert ^{m}+\left\vert
\lambda\right\vert ^{m}\right)  & \text{if }\left\vert z-\lambda\right\vert
>\epsilon\left\vert \lambda\right\vert
\end{array}
\right.  \text{.}%
\]
Set%
\[
f_{1}\left(  z,\lambda\right)  =\left\vert z\right\vert ^{-m}\left(
\left\vert z\right\vert ^{m+n-1}I_{\left\vert z-\lambda\right\vert
\leq\epsilon\left\vert \lambda\right\vert }+\left\vert \lambda\right\vert
^{-1}\left(  \left\vert z\right\vert ^{m}+\left\vert \lambda\right\vert
^{m}\right)  I_{\left\vert z-\lambda\right\vert >\epsilon\left\vert
\lambda\right\vert }\right)  \text{.}%
\]
Since $\left\vert r\left(  z\right)  ^{-1}g_{k}\left(  z\right)
^{-1}\right\vert \leq c_{2}\left\vert z\right\vert ^{-m}$, it is sufficient to
show the integrability of $f_{1}\left(  z,\lambda\right)  ^{2}$. The rest of
proof proceeds just as the proof of Lemma \ref{l10}. The exponent of the first
term is $2\left(  n-1\right)  -2N$, which is less than if $N\geq n$. The
integral of the second term is dominated by%
\[
\int_{C^{2}}\left\vert \lambda\right\vert ^{-2}\left\vert z\right\vert
^{-2m}\left(  \left\vert z\right\vert ^{2m}+\left\vert \lambda\right\vert
^{2m}\right)  \left\vert z\right\vert ^{-2N}\left\vert dz\right\vert
\left\vert d\lambda\right\vert \text{,}%
\]
which is finite if%
\[
-2N<-1\text{, }-2m-2N<-1\Longrightarrow N\geq1\text{, }1-m\text{.}%
\]

\end{proof}

\section{Non-negativity condition of $\boldsymbol{A}_{N}^{inv}\left(
C\right)  $}

Generally potentials arising from $\boldsymbol{a}\in\boldsymbol{A}_{N}%
^{inv}\left(  C\right)  $ are complex valued, so to obtain real potentials
some sort of realness for $\boldsymbol{a}$ and $g$ is required.
$\boldsymbol{a}\in\boldsymbol{A}_{N}\left(  C\right)  $, $g\in\Gamma
_{n}^{\left(  m\right)  }$ are called \textbf{real} if they satisfy%
\begin{equation}
\boldsymbol{a}\left(  \lambda\right)  =\overline{\boldsymbol{a}\left(
\overline{\lambda}\right)  }\text{ \ for }\lambda\in C\text{, \ \ }%
g(z)=\overline{g(\overline{z})}\text{ \ for }z\in\mathbb{C}\text{.} \label{49}%
\end{equation}
If $\boldsymbol{a}$ and $g$ are real in this sense, then clearly we have%
\[
\left\{
\begin{array}
[c]{l}%
\varphi_{\boldsymbol{a}}\left(  z\right)  =\overline{\varphi_{\boldsymbol{a}%
}\left(  \overline{z}\right)  }\text{, \ }\psi_{\boldsymbol{a}}\left(
z\right)  =\overline{\psi_{\boldsymbol{a}}\left(  \overline{z}\right)
}\text{, \ }m_{\boldsymbol{a}}\left(  z\right)  =\overline{m_{\boldsymbol{a}%
}\left(  \overline{z}\right)  }\\
\tau_{\boldsymbol{a}}(g)\text{, }\tau_{\boldsymbol{a}}^{\left(  2\right)
}(g)\in\mathbb{R}%
\end{array}
\right.  \text{,}%
\]
and the associated potential takes real values.

Define a subclass of $\boldsymbol{A}_{L}^{inv}$:%
\begin{align*}
\boldsymbol{A}_{L,+}^{inv}\left(  C\right)   &  =\left\{  \boldsymbol{a}%
\in\boldsymbol{A}_{L}^{inv}\left(  C\right)  \text{; \ }\tau_{\boldsymbol{a}%
}^{\left(  2\right)  }\left(  r\right)  \geq0\text{ for any real rational
}r\in\Gamma_{0}^{\left(  0\right)  }\right\} \\
&  =\left\{  \boldsymbol{a}\in\boldsymbol{A}_{L}^{inv}\left(  C\right)
\text{; \ }\tau_{\boldsymbol{a}}\left(  r\right)  \geq0\text{ for any real
rational }r\in\Gamma_{0}^{\left(  0\right)  }\right\}  \text{.}%
\end{align*}
The second identity follows from the identity%
\[
\tau_{\boldsymbol{a}}^{\left(  2\right)  }\left(  r\right)  =\tau
_{\boldsymbol{a}}\left(  r\right)  \exp\left(  -\text{\textrm{tr}}\left(
r^{-1}T\left(  r\boldsymbol{a}\right)  T\left(  \boldsymbol{a}\right)
^{-1}-I\right)  \right)  \text{.}%
\]
$\tau_{\boldsymbol{a}}\left(  r\right)  $ is well-defined for any rational
function $r\in\Gamma_{0}^{\left(  0\right)  }$ since the relevant operator is
of finite rank. Our strategy to show $\tau_{\boldsymbol{a}}^{\left(  2\right)
}\left(  g\right)  >0$ for real $g\in\Gamma_{n}^{\left(  0\right)  }$ is as
follows:\medskip\newline(i) \ \ Show $\tau_{\boldsymbol{a}}^{\left(  2\right)
}\left(  r\right)  >0$ for any real $r\in\Gamma_{0}^{\left(  0\right)  }$ and
$\boldsymbol{a}\in\boldsymbol{A}_{L,+}^{inv}\left(  C\right)  $.\newline(ii)
\ Approximate a general real $g\in\Gamma_{n}^{\left(  0\right)  }$ by a
sequence of real $r_{k}\in\Gamma_{0}^{\left(  0\right)  }$.\newline(iii) Use
the continuity of $\tau_{\boldsymbol{a}}^{\left(  2\right)  }\left(
\cdot\right)  $ to have $\tau_{\boldsymbol{a}}^{\left(  2\right)  }\left(
gr_{k}^{-1}\right)  >0$ for sufficiently large $k$ and show $gr_{k}%
^{-1}\boldsymbol{a}\in\boldsymbol{A}_{L,+}^{inv}\left(  C\right)  $%
.\newline(iv) Apply the cocycle property of $\tau_{\boldsymbol{a}}^{\left(
2\right)  }\left(  \cdot\right)  $ to have $\tau_{\boldsymbol{a}}^{\left(
2\right)  }\left(  g\right)  >0$, namely%
\[
\tau_{\boldsymbol{a}}^{\left(  2\right)  }\left(  g\right)  =\tau
_{\boldsymbol{a}}^{\left(  2\right)  }\left(  gr_{k}^{-1}r_{k}\right)
=\tau_{\boldsymbol{a}}^{\left(  2\right)  }\left(  gr_{k}^{-1}\right)
\tau_{gr_{k}^{-1}\boldsymbol{a}}^{\left(  2\right)  }\left(  r_{k}\right)
\exp\left(  -E_{\boldsymbol{a}}\left(  gr_{k}^{-1},r_{k}\right)  \right)
>0\text{.}%
\]
This programme will be realized in the next section. At the same time a close
connection of the $m$-function with the Weyl function for Schr\"{o}dinger
operators will be revealed.

\subsection{Non-degeneracy of Tau-functions for $\boldsymbol{a}\in
\boldsymbol{A}_{L,+}^{inv}\left(  C\right)  $}

To investigate properties of $m_{\boldsymbol{a}}$ and $\tau_{\boldsymbol{a}%
}^{\left(  2\right)  }$ for $\boldsymbol{a}\in\boldsymbol{A}_{L,+}%
^{inv}\left(  C\right)  $ we prepare several lemmas. In this section the curve
$C$ is parametrized by%
\[
C=\left\{  \pm\omega\left(  y\right)  +iy\text{; \ }y\in\mathbb{R}\right\}
\]
with a smooth function $\omega\left(  y\right)  >0$ satisfying $\omega\left(
y\right)  =\omega\left(  -y\right)  $ and $\omega\left(  y\right)  =O\left(
y^{-\left(  n-1\right)  }\right)  $.

\begin{lemma}
\label{l6}Let $a(x)$, $b(x)$ be real valued analytic functions on an interval
$I\subset\mathbb{R}$ satisfying%
\[
\left\{
\begin{array}
[c]{l}%
\dfrac{a(y)b(x)-a(x)b(y)}{x-y}\geq0\\
a(x)^{2}+b(x)^{2}>0\\
a\left(  x\right)  \geq0
\end{array}
\text{ for any }x\text{, }y\in I\right.
\]
Then we have either$\ a(x)=0$ identically or $a(x)>0$ for any $x\in I$ .
\end{lemma}

\begin{proof}
Since the first assumption implies%
\[
a(y)b(x)-a(x)b(y)\geq0\text{ \ for any }x>y\text{,}%
\]
if $a(x)$ has a zero at $x_{0}\in I$, setting $x=x_{0}$ or $y=x_{0}$ we have%
\begin{equation}
\left\{
\begin{array}
[c]{l}%
a(y)b(x_{0})\geq0\text{ \ for any }y<x_{0}\\
a(x)b(x_{0})\leq0\text{ \ for any }x>x_{0}%
\end{array}
\right.  \text{.} \label{76}%
\end{equation}
The second assumption implies $b(x_{0})\neq0$, hence (\ref{76}) together with
the property $a(x)\geq0$ shows $a(x)=0$ on $\left(  -\infty,x_{0}\right)  \cap
I$ or $\left(  x_{0},\infty\right)  \cap I$. Then the analyticity of $a$
yields the vanishing of $a(x)$ on $I$.\medskip
\end{proof}

In what follows $\tau_{\boldsymbol{a}}\left(  r\right)  $ for $r\in\Gamma
_{0}^{\left(  0\right)  }$ will be used instead of $\tau_{\boldsymbol{a}%
}^{\left(  2\right)  }\left(  r\right)  $. Recall that $\tau_{\boldsymbol{a}%
}\left(  r\right)  \geq0$ holds for any real $r\in\Gamma_{0}^{\left(
0\right)  }$ if $\boldsymbol{a}\in\boldsymbol{A}_{L,+}^{inv}\left(  C\right)
$. In what follows we have to assume $\boldsymbol{a}\in\boldsymbol{A}%
_{L,+}^{inv}\left(  C\right)  $ with $L\geq2$ since in the proofs we use
$\varphi_{\boldsymbol{a}}$, $m_{\boldsymbol{a}}$. Let%
\[
p_{\zeta}\left(  z\right)  =1+\zeta^{-1}z\text{, \ }q_{\zeta}\left(  z\right)
=\left(  1-\zeta^{-1}z\right)  ^{-1}\text{ \ for }\zeta\in D_{-}\text{.}%
\]

\begin{lemma}
\label{l19}Let $\boldsymbol{a}\in\boldsymbol{A}_{L,+}^{inv}\left(  C\right)  $
with $L\geq2$. Then the followings are valid.\newline(i) $\ \operatorname{Im}%
m_{\boldsymbol{a}}\left(  z\right)  /\operatorname{Im}z>0$ on $D_{-}%
\backslash\mathbb{R}$ and $m_{\boldsymbol{a}}\left(  z\right)  $ is analytic
on $D_{-}$.\newline(ii) $\ 1+\varphi_{\boldsymbol{a}}\left(  z\right)  \neq0$
on $D_{-}$.\newline(iii) $\tau_{\boldsymbol{a}}\left(  q_{x_{1}}p_{x_{2}%
}\right)  =\left(  1+\varphi_{\boldsymbol{a}}\left(  x_{1}\right)  \right)
\left(  1+\varphi_{\boldsymbol{a}}\left(  x_{2}\right)  \right)
\dfrac{m_{\boldsymbol{a}}\left(  x_{1}\right)  -m_{\boldsymbol{a}}\left(
x_{2}\right)  }{\Delta_{\boldsymbol{a}}\left(  x_{2}\right)  \left(
x_{1}-x_{2}\right)  }>0$ for any $x_{1}$, $x_{2}\in D_{-}\cap\mathbb{R}$ if
$x_{1}\neq x_{2}$.
\end{lemma}

\begin{proof}
In the formula in Lemma \ref{l13} setting $\zeta_{1}=\zeta$, $\zeta
_{2}=\overline{\zeta}$, $\zeta_{1}^{\prime}=\eta$, $\zeta_{2}^{\prime
}=\overline{\eta}\in D_{-}$, we see that $q_{\zeta}q_{\overline{\zeta}}%
p_{\eta}p_{\overline{\eta}}$ is a real element of $\Gamma_{0}^{\left(
0\right)  }$, hence $\tau_{\boldsymbol{a}}\left(  q_{\zeta}q_{\overline{\zeta
}}p_{\eta}p_{\overline{\eta}}\right)  \geq0$ if $\boldsymbol{a}\in
\boldsymbol{A}_{L,+}^{inv}\left(  C\right)  $. Lemma \ref{l13} implies that%
\[
\lim_{\eta\rightarrow\infty}\tau_{\boldsymbol{a}}\left(  q_{\zeta}%
q_{\overline{\zeta}}p_{\eta}p_{\overline{\eta}}\right)  =\tau_{\boldsymbol{a}%
}\left(  q_{\zeta}q_{\overline{\zeta}}\right)  \text{,}%
\]
hence we have $\tau_{\boldsymbol{a}}\left(  q_{\zeta}q_{\overline{\zeta}%
}\right)  \geq0$ for all $\zeta\in D_{-}$. On the other hand from Lemma
\ref{l13} we have%
\[
\tau_{\boldsymbol{a}}\left(  q_{\zeta}q_{\overline{\zeta}}\right)  =\left\vert
\varphi_{\boldsymbol{a}}\left(  \zeta\right)  +1\right\vert ^{2}%
\dfrac{\operatorname{Im}m_{\boldsymbol{a}}\left(  \zeta\right)  }%
{\operatorname{Im}\zeta}\text{,}%
\]
which implies $\operatorname{Im}m_{\boldsymbol{a}}\left(  \zeta\right)
/\operatorname{Im}\zeta\geq0$ for $\zeta\in\left\{  \zeta\in D_{-}\text{;
}\varphi_{\boldsymbol{a}}\left(  \zeta\right)  +1\neq0\right\}  \equiv
\mathcal{Z}_{\varphi}$. Note $\varphi_{\boldsymbol{a}}\left(  \zeta\right)
\rightarrow0$ as $\zeta\rightarrow\infty$, hence the set $\mathcal{Z}%
_{\varphi}$ is discrete. Suppose $m_{\boldsymbol{a}}\left(  \zeta_{0}\right)
=0$ for some $\zeta_{0}\in\mathcal{Z}_{\varphi}\cap\mathbb{C}_{+}$. Then the
maximum principle for the harmonic function $-\operatorname{Im}%
m_{\boldsymbol{a}}\left(  \zeta\right)  $ shows $\operatorname{Im}%
m_{\boldsymbol{a}}\left(  \zeta\right)  =0$ identically there, which
contradicts $m_{\boldsymbol{a}}\left(  \zeta\right)  =\zeta+O\left(
\zeta^{-1}\right)  $ as $\zeta\rightarrow\infty$, hence $>0$ on $\mathcal{Z}%
_{\varphi}\cap\mathbb{C}_{+}$. One has the same property also in
$\mathbb{C}_{-}$. On the other hand $m_{\boldsymbol{a}}\left(  \zeta\right)  $
has poles at $\zeta_{0}\in D_{-}\backslash\mathcal{Z}_{\varphi}$ since
$\varphi_{\boldsymbol{a}}\left(  \zeta\right)  +1$ and $\psi_{\boldsymbol{a}%
}\left(  \zeta\right)  +\zeta$ do not vanish simultaneously due to
$\Delta_{\boldsymbol{a}}\left(  \zeta\right)  \neq0$. However this is
impossible if we apply the same argument to $m_{\boldsymbol{a}}\left(
\zeta\right)  ^{-1}$, which implies%
\[
m_{\boldsymbol{a}}\left(  z\right)  \ \text{is analytic and }\frac
{\operatorname{Im}m_{\boldsymbol{a}}\left(  z\right)  }{\operatorname{Im}%
z}>0\text{ holds on }D_{-}\backslash\mathbb{R}\text{.}%
\]
This in particular means $1+\varphi_{\boldsymbol{a}}\left(  z\right)  \neq0$
on $D_{-}\backslash\mathbb{R}$. The remaining problem is the existence or
non-existence of poles of $m_{\boldsymbol{a}}\left(  z\right)  $ on $D_{-}%
\cap\mathbb{R}$. We rewrite the identity in Lemma \ref{l13} as%
\[
\tau_{\boldsymbol{a}}\left(  q_{\zeta_{1}}p_{\zeta_{2}}\right)  =\dfrac
{\left(  \psi_{\boldsymbol{a}}\left(  \zeta_{2}\right)  +\zeta_{2}\right)
\left(  \varphi_{\boldsymbol{a}}\left(  \zeta_{1}\right)  +1\right)  -\left(
\varphi_{\boldsymbol{a}}\left(  \zeta_{2}\right)  +1\right)  \left(
\psi_{\boldsymbol{a}}\left(  \zeta_{1}\right)  +\zeta_{1}\right)  }%
{\Delta_{\boldsymbol{a}}\left(  \zeta_{2}\right)  \left(  \zeta_{1}-\zeta
_{2}\right)  }\text{,}%
\]
which is valid for any $\zeta_{1}$, $\zeta_{2}\in D_{-}$. Set $\zeta_{1}%
=x_{1}$, $\zeta_{2}=x_{2}\in D_{-}\cap\mathbb{R}$ and%
\[
a(x)=\varphi_{\boldsymbol{a}}\left(  x\right)  +1\text{, \ }b(x)=\psi
_{\boldsymbol{a}}\left(  x\right)  +x\text{.}%
\]
Then%
\[
\tau_{\boldsymbol{a}}\left(  q_{x_{1}}p_{x_{2}}\right)  =\dfrac{a\left(
x_{1}\right)  b\left(  x_{2}\right)  -a\left(  x_{2}\right)  b\left(
x_{1}\right)  }{\Delta_{\boldsymbol{a}}\left(  x_{2}\right)  \left(
x_{1}-x_{2}\right)  }%
\]
holds, and the property $\Delta_{\boldsymbol{a}}\left(  x\right)  \neq0$,
$\Delta_{\boldsymbol{a}}\left(  x\right)  \rightarrow-1$ as $\left\vert
x\right\vert \rightarrow\infty$ implies $\Delta_{\boldsymbol{a}}\left(
x\right)  <0$ on $D_{-}\cap\mathbb{R}$, which together with $\tau
_{\boldsymbol{a}}\left(  q_{x_{1}}p_{x_{2}}\right)  \geq0$ (due to
$\boldsymbol{a}\in\boldsymbol{A}_{L,+}^{inv}\left(  C\right)  $) shows%
\[
\dfrac{a\left(  x_{1}\right)  b\left(  x_{2}\right)  -a\left(  x_{2}\right)
b\left(  x_{1}\right)  }{x_{1}-x_{2}}\leq0\text{.}%
\]
Moreover, in the inequality $\tau_{\boldsymbol{a}}\left(  q_{x_{1}}p_{x_{2}%
}\right)  \geq0$ letting $x_{2}\rightarrow\infty$, we have $\tau
_{\boldsymbol{a}}\left(  q_{x_{1}}p_{x_{2}}\right)  \rightarrow\tau
_{\boldsymbol{a}}\left(  q_{x_{1}}\right)  $ just by the same argument as
above and%
\[
a\left(  x_{1}\right)  =\varphi_{\boldsymbol{a}}\left(  x_{1}\right)
+1=\tau_{\boldsymbol{a}}\left(  q_{x_{1}}\right)  \geq0
\]
is valid. The condition $a\left(  x\right)  ^{2}+b(x)^{2}>0$ follows from
$\Delta_{\boldsymbol{a}}\left(  x\right)  \neq0$, hence one can apply Lemma
\ref{l6} to have $a(x)>0$, since $a\left(  x\right)  \rightarrow1$ as
$\left\vert x\right\vert \rightarrow\infty$. We have shown (i) and (ii). To
show (iii) first assume $x_{1}\neq x_{2}$ and suppose $\tau_{\boldsymbol{a}%
}\left(  q_{x_{1}}p_{x_{2}}\right)  =0$. Then
\[
0=\dfrac{a\left(  x_{1}\right)  b\left(  x_{2}\right)  -a\left(  x_{2}\right)
b\left(  x_{1}\right)  }{x_{1}-x_{2}}=a\left(  x_{1}\right)  a\left(
x_{2}\right)  \dfrac{c\left(  x_{2}\right)  -c\left(  x_{1}\right)  }%
{x_{1}-x_{2}}%
\]
with $c(x)=b(x)/a(x)$. Observing $c(x)$ is analytic and $c^{\prime}(x)\geq0$,
this identity implies $c(x)$ is identically constant, which contradicts
$c(x)=x+o\left(  1\right)  $ as $\left\vert x\right\vert \rightarrow\infty$,
hence we have $\tau_{\boldsymbol{a}}\left(  q_{x_{1}}p_{x_{2}}\right)  >0$ if
$x_{1}\neq x_{2}$.\medskip
\end{proof}

\begin{lemma}
\label{l17}Let $\boldsymbol{a}\in\boldsymbol{A}_{L,+}^{inv}\left(  C\right)  $
with $L\geq2$ and $r\in\Gamma_{0}^{\left(  0\right)  }$ be real. Then
$\tau_{\boldsymbol{a}}\left(  r\right)  >0$ holds, which in particular means
$r\boldsymbol{a}\in\boldsymbol{A}_{L,+}^{inv}\left(  C\right)  $.
\end{lemma}

\begin{proof}
First note that for $r_{1}$, $r_{2}\in\Gamma_{0}^{\left(  0\right)  }$ and
$\boldsymbol{a}\in\boldsymbol{A}_{L,+}^{inv}$%
\[
\tau_{\boldsymbol{a}}\left(  r_{1}\right)  >0\text{, }\tau_{\boldsymbol{a}%
}\left(  r_{2}\right)  >0\Longrightarrow\tau_{\boldsymbol{a}}\left(
r_{1}r_{2}\right)  >0\text{.}%
\]
This is because $r_{1}\boldsymbol{a}\in\boldsymbol{A}_{L,+}^{inv}\left(
C\right)  $ since the cocycle property implies%
\[
\tau_{r_{1}\boldsymbol{a}}\left(  r\right)  =\frac{\tau_{\boldsymbol{a}%
}\left(  rr_{1}\right)  }{\tau_{\boldsymbol{a}}\left(  r_{1}\right)  }%
\geq0\text{ \ for any real }r\in\Gamma_{0}^{\left(  0\right)  }\text{,}%
\]
and
\[
\tau_{\boldsymbol{a}}\left(  r_{1}r_{2}\right)  =\tau_{\boldsymbol{a}}\left(
r_{1}\right)  \tau_{r_{1}\boldsymbol{a}}\left(  r_{2}\right)  >0\text{.}%
\]
Generally real $r\in\Gamma_{0}^{\left(  0\right)  }$ is a product of%
\begin{equation}
\left\{
\begin{array}
[c]{l}%
\text{(1) }q_{\zeta}q_{\overline{\zeta}}p_{\eta}p_{\overline{\eta}}\text{ with
}\eta\text{, }\zeta\in D_{-}\backslash\mathbb{R}\\
\text{(2) }q_{s}p_{t}\text{ with }s\text{, }t\in D_{-}\cap\mathbb{R}\\
\text{(3) }q_{\zeta}q_{\overline{\zeta}}p_{s}p_{t}\text{ with }\zeta\in
D_{-}\backslash\mathbb{R}\text{, }s\text{, }t\in D_{-}\cap\mathbb{R}\\
\text{(4) }q_{s}q_{t}p_{\eta}p_{\overline{\eta}}\text{ with }\eta\in
D_{-}\backslash\mathbb{R}\text{, }s\text{, }t\in D_{-}\cap\mathbb{R}%
\end{array}
\right.  \text{.} \label{86}%
\end{equation}
Therefore, if $\tau_{\boldsymbol{a}}\left(  r\right)  >0$ is proved for any
$r$ of these 4 cases, we have $\tau_{\boldsymbol{a}}\left(  r\right)  >0$ for
any real $r\in\Gamma_{0}^{\left(  0\right)  }$.

We begin from the case (1) and let $r\left(  z\right)  =\left(  q_{\zeta
}q_{\overline{\zeta}}p_{\eta}p_{\overline{\eta}}\right)  \left(  z\right)  $.
Then (\ref{60}) of Lemma \ref{l13} implies%
\begin{align}
&  \tau_{\boldsymbol{a}}\left(  r\right)  =\dfrac{\left\vert \eta
+\zeta\right\vert ^{2}\left\vert \eta+\overline{\zeta}\right\vert
^{2}\left\vert \varphi_{\boldsymbol{a}}\left(  \eta\right)  +1\right\vert
^{2}\left\vert \varphi_{\boldsymbol{a}}\left(  \zeta\right)  +1\right\vert
^{2}}{\left(  4\operatorname{Im}\eta\operatorname{Im}\zeta\right)  \left\vert
\Delta_{\boldsymbol{a}}\left(  \eta\right)  \right\vert ^{2}}\nonumber\\
&  \text{ \ \ \ \ \ \ \ \ \ \ \ \ \ \ \ \ \ \ \ \ \ \ }\times\left(
\left\vert \dfrac{m_{\boldsymbol{a}}\left(  \eta\right)  -\overline
{m_{\boldsymbol{a}}\left(  \zeta\right)  }}{\eta^{2}-\overline{\zeta}^{2}%
}\right\vert ^{2}-\left\vert \dfrac{m_{\boldsymbol{a}}\left(  \eta\right)
-m_{\boldsymbol{a}}\left(  \zeta\right)  }{\eta^{2}-\zeta^{2}}\right\vert
^{2}\right)  \label{84}%
\end{align}
due to realness of $\boldsymbol{a}$. Note $\varphi_{\boldsymbol{a}}\left(
z\right)  +1\neq0$ for any $z\in D_{-}$ due to Lemma \ref{l19}. Owing to the
symmetry of $r$ with respect to $\zeta$, $\eta$, one can assume
$\operatorname{Im}\zeta>0$, $\operatorname{Im}\eta>0$. The condition
$\boldsymbol{a}\in\boldsymbol{A}_{L,+}^{inv}\left(  C\right)  $ implies
$\tau_{\boldsymbol{a}}\left(  q_{\zeta}q_{\overline{\zeta}}p_{\eta
}p_{\overline{\eta}}\right)  \geq0$ for any $\eta$, $\zeta\in D_{-}%
\cap\mathbb{C}_{+}$. Assume $\tau_{\boldsymbol{a}}\left(  q_{\zeta_{0}%
}q_{\overline{\zeta_{0}}}p_{\eta_{0}}p_{\overline{\eta_{0}}}\right)  =0$ for
some $\eta_{0}$, $\zeta_{0}\in D_{-}$ and consider the analytic function%
\[
f(z)=\dfrac{m_{\boldsymbol{a}}\left(  z\right)  -m_{\boldsymbol{a}}\left(
\zeta_{0}\right)  }{m_{\boldsymbol{a}}\left(  z\right)  -\overline
{m_{\boldsymbol{a}}\left(  \zeta_{0}\right)  }}\dfrac{z^{2}-\overline
{\zeta_{0}}^{2}}{z^{2}-\zeta_{0}^{2}}\text{.}%
\]
The property $\tau_{\boldsymbol{a}}\left(  q_{\zeta_{0}}q_{\overline{\zeta
_{0}}}p_{z}p_{\overline{z}}\right)  \geq0$ shows%
\begin{equation}
\left\vert f(z)\right\vert \leq1\text{ \ \ for any }z\in D_{-}\cap
\mathbb{C}_{+}\text{,} \label{72}%
\end{equation}
and the assumption implies the equality at (\ref{72}) for $z=\eta_{0}$, which
concludes $f(z)=e^{i\alpha}$ with $\alpha\in\mathbb{R}$ identically on
$D_{-}\cap\mathbb{C}_{+}$. Then%
\[
m_{\boldsymbol{a}}\left(  z\right)  =\dfrac{\left(  m_{\boldsymbol{a}}\left(
\zeta_{0}\right)  -e^{i\alpha}\overline{m_{\boldsymbol{a}}\left(  \zeta
_{0}\right)  }\right)  z^{2}+e^{i\alpha}\overline{m_{\boldsymbol{a}}\left(
\zeta_{0}\right)  }\zeta_{0}^{2}-m_{\boldsymbol{a}}\left(  \zeta_{0}\right)
\overline{\zeta_{0}}^{2}}{\left(  1-e^{i\alpha}\right)  z^{2}+e^{i\alpha}%
\zeta_{0}^{2}-\overline{\zeta_{0}}^{2}}%
\]
holds, which contradicts $m_{\boldsymbol{a}}\left(  z\right)  =z+o\left(
1\right)  $ as $z\rightarrow\infty$. Therefore we have $\left\vert
f(z)\right\vert <1$ always, which is nothing but $\tau_{\boldsymbol{a}}\left(
r\right)  >0$.

The case (2) is already proved in Lemma \ref{l6} if $s\neq t$. Suppose $s=t$
and let $s_{n}\in D_{-}\cap\mathbb{R}$ be a sequence converging to $s$. Then,
$\tau_{\boldsymbol{a}}\left(  q_{s_{n}}p_{s}\right)  >0$ is valid due to
$s_{n}\neq s$ for any $n\geq1$. Moreover one can show easily $\tau
_{\boldsymbol{a}}\left(  q_{s}q_{s_{n}}^{-1}\right)  \rightarrow1$ as
$n\rightarrow\infty$. Taking sufficiently large $n$ such that $\tau
_{\boldsymbol{a}}\left(  q_{s}q_{s_{n}}^{-1}\right)  >0$ and fixing it we see
from the cocycle property%
\[
\tau_{\boldsymbol{a}}\left(  q_{s}p_{s}\right)  =\tau_{\boldsymbol{a}}\left(
r_{n}\left(  q_{s_{n}}p_{s}\right)  \right)  =\tau_{\boldsymbol{a}}\left(
r_{n}\right)  \tau_{r_{n}\boldsymbol{a}}\left(  q_{s_{n}}p_{s}\right)
>0\text{,}%
\]
where $r_{n}=q_{s}q_{s_{n}}^{-1}$, since $r_{n}\boldsymbol{a}\in
\boldsymbol{A}_{L,+}^{inv}\left(  C\right)  $ due to $\tau_{\boldsymbol{a}%
}\left(  r_{n}\right)  >0$.

Similarly one can show $\tau_{\boldsymbol{a}}\left(  q_{\zeta}q_{\overline
{\zeta}}p_{s}p_{s}\right)  >0$ as a limiting case of (1). Setting
$r_{1}=q_{\zeta}q_{\overline{\zeta}}p_{s}p_{s}$ we have%
\[
\tau_{\boldsymbol{a}}\left(  q_{\zeta}q_{\overline{\zeta}}p_{s}p_{t}\right)
=\tau_{\boldsymbol{a}}\left(  r_{1}q_{-s}p_{t}\right)  =\tau_{\boldsymbol{a}%
}\left(  r_{1}\right)  \tau_{r_{1}\boldsymbol{a}}\left(  q_{-s}p_{t}\right)
>0\text{,}%
\]
which shows the case (3). The case (4) can be shown similarly.\medskip
\end{proof}

The next task is to approximate general $g\in\Gamma_{n}$ by rational functions.

\begin{lemma}
\label{l18}Let $g\in\Gamma_{n}$. Then there exists a sequence of rational
functions $\left\{  r_{k}\right\}  _{k\geq1}\subset\Gamma_{0}^{\left(
0\right)  }$ such that $r_{k}\rightarrow g$ in the sense of (\ref{61}).
\end{lemma}

\begin{proof}
Let $U$ be a neighborhood of $\overline{D}_{+}$ whose boundary is described by
an equation $\left\vert x\right\vert =c\left\vert y\right\vert ^{-\left(
n-1\right)  }$ for large $\left\vert y\right\vert $ with sufficiently large
$c>0$. For integer $k\geq1$ let%
\[
\phi_{k}\left(  z\right)  =\left(  \frac{1+\dfrac{z}{2k}}{1-\dfrac{z}{2k}%
}\right)  ^{k}\text{.}%
\]
Note $\lim_{k\rightarrow\infty}\phi_{k}\left(  z\right)  =e^{z}$. For a
positive constant $a\leq k$ an inequality%
\begin{equation}
e^{-2a}\leq\left\vert \phi_{k}\left(  z\right)  \right\vert \leq e^{2a}\text{
\ \ \ if \ }\left\vert \operatorname{Re}z\right\vert \leq a \label{73}%
\end{equation}
holds. If $h(z)=c_{1}z^{n}+$ lower degree terms, then%
\begin{equation}
c_{2}\equiv\sup_{z\in U}\left\vert \operatorname{Re}h\left(  z\right)
\right\vert <\infty\label{74}%
\end{equation}
is valid. Define real rational functions by%
\[
r_{k}(z)=\phi_{k}\left(  h\left(  z\right)  \right)  \text{.}%
\]
The zeros and poles of $r_{k}(z)$ are determined by the equation $h\left(
z\right)  =\pm2k$. If $a$ is chosen so that $a>c_{2}$, then clearly there
exist a constant $c_{3}>1$ such that%
\[
c_{3}^{-1}\leq\left\vert r_{k}(z)\right\vert \leq c_{3}\text{ \ \ for \ }z\in
U\text{ and }k\geq1
\]
holds. Moreover%
\[
\left\vert r_{k}^{\prime}(z)\right\vert =\left\vert h^{\prime}(z)\right\vert
\left\vert \phi_{k}^{\prime}\left(  h\left(  z\right)  \right)  \right\vert
=\left\vert h^{\prime}(z)\right\vert \left\vert \frac{1+\dfrac{h\left(
z\right)  }{2k}}{1-\dfrac{h\left(  z\right)  }{2k}}\right\vert ^{k-1}%
\left\vert 1-\dfrac{h\left(  z\right)  }{2k}\right\vert ^{-2}%
\]
shows%
\[
\left\vert r_{k}^{\prime}(z)\right\vert \leq c_{4}\left\vert z\right\vert
^{n-1}\text{\ \ for \ }z\in U\text{ and }k\geq1\text{.}%
\]
Since $\lim_{k\rightarrow\infty}r_{k}(z)=e^{h(z)}=g(z)$, all conditions of
Lemma \ref{l11} are satisfied.\medskip
\end{proof}

Now we have

\begin{proposition}
\label{p6}Let $g\in\Gamma_{n}^{\left(  0\right)  }$ be real and
$\boldsymbol{a}\in\boldsymbol{A}_{L,+}^{inv}\left(  C\right)  $ with
$L\geq\max\left\{  n,2\right\}  $. Then, $\tau_{\boldsymbol{a}}^{\left(
2\right)  }\left(  g\right)  >0$ holds, hence $g\boldsymbol{a}\in
\boldsymbol{A}_{L,+}^{inv}\left(  C\right)  $ is valid.
\end{proposition}

\begin{proof}
$L\geq\max\left\{  n,1\right\}  $ is necessary for the definition of
$\tau_{\boldsymbol{a}}^{\left(  2\right)  }\left(  g\right)  $ for $g\in
\Gamma_{n}^{\left(  0\right)  }$ and $L\geq2$ is required to apply Lemma
\ref{l17}. First note that if $\tau_{\boldsymbol{a}}^{\left(  2\right)
}\left(  g\right)  >0$ holds, then $g\boldsymbol{a}\in\boldsymbol{A}%
_{L,+}^{inv}\left(  C\right)  $. To show this let $\left\{  r_{k}\right\}
_{k\geq1}$ be the sequence of Lemma \ref{l18} approximating $g$ and $r$ be any
real rational function of $\Gamma_{0}^{\left(  0\right)  }$. Then
$r_{k}\rightarrow g$ implies%
\[
\tau_{\boldsymbol{a}}^{\left(  2\right)  }\left(  gr\right)  =\lim
_{k\rightarrow\infty}\tau_{\boldsymbol{a}}^{\left(  2\right)  }\left(
r_{k}r\right)  \geq0
\]
and%
\[
\tau_{g\boldsymbol{a}}^{\left(  2\right)  }\left(  r\right)  =\dfrac
{\tau_{\boldsymbol{a}}^{\left(  2\right)  }\left(  gr\right)  }{\tau
_{\boldsymbol{a}}^{\left(  2\right)  }\left(  g\right)  }\exp\left(
E_{\boldsymbol{a}}\left(  g,r\right)  \right)  \geq0\text{,}%
\]
which means $g\boldsymbol{a}\in\boldsymbol{A}_{L,+}^{inv}\left(  C\right)  $.

Now $\left\{  g_{k}=gr_{k}^{-1}\right\}  _{k\geq1}$ also satisfies the
conditions of (\ref{61}) with $g=1$, hence Lemma \ref{l11} shows%
\[
\lim_{k\rightarrow\infty}\tau_{\boldsymbol{a}}^{\left(  2\right)  }\left(
gr_{k}^{-1}\right)  =\tau_{\boldsymbol{a}}^{\left(  2\right)  }\left(
1\right)  =1\text{.}%
\]
Fix a sufficiently large $k\geq1$ such that $\tau_{\boldsymbol{a}}^{\left(
2\right)  }\left(  gr_{k}^{-1}\right)  >0$. Then the above argument shows
$gr_{k}^{-1}\boldsymbol{a}\in\boldsymbol{A}_{L,+}^{inv}\left(  C\right)  $.
Applying Lemma \ref{l17} to $gr_{k}^{-1}\boldsymbol{a}\in\boldsymbol{A}%
_{L,+}^{inv}\left(  C\right)  $ and the rational function $r_{k}$ we have
$\tau_{gr_{k}^{-1}\boldsymbol{a}}^{\left(  2\right)  }\left(  r_{k}\right)
>0$. The cocycle property of tau-functions implies%
\[
\tau_{\boldsymbol{a}}^{\left(  2\right)  }\left(  g\right)  =\tau
_{\boldsymbol{a}}^{\left(  2\right)  }\left(  gr_{k}^{-1}r_{k}\right)
=\tau_{\boldsymbol{a}}^{\left(  2\right)  }\left(  gr_{k}^{-1}\right)
\tau_{gr_{k}^{-1}\boldsymbol{a}}^{\left(  2\right)  }\left(  r_{k}\right)
\exp\left(  -E_{\boldsymbol{a}}\left(  gr_{k}^{-1},r_{k}\right)  \right)
>0\text{.}%
\]
If $g=re^{h}$ with real $r\in\Gamma_{0}^{\left(  0\right)  }$, then%
\[
\tau_{\boldsymbol{a}}^{\left(  2\right)  }\left(  g\right)  =\tau
_{\boldsymbol{a}}^{\left(  2\right)  }\left(  r\right)  \tau_{r\boldsymbol{a}%
}^{\left(  2\right)  }\left(  e^{h}\right)  \exp\left(  -E_{\boldsymbol{a}%
}\left(  r,e^{h}\right)  \right)  >0\text{,}%
\]
which completes the proof.
\end{proof}

\subsection{$m$-function and Weyl function}

Since for $\boldsymbol{a}\in\boldsymbol{A}_{L,+}^{inv}\left(  C\right)  $ we
know $e_{x}\boldsymbol{a}\in\boldsymbol{A}_{L,+}^{inv}\left(  C\right)  $
($e_{x}(z)=e^{xz}$) from Theorem \ref{t1}, Proposition \ref{p1} shows that we
have a Schr\"{o}dinger equation%
\[
-\partial_{x}^{2}f_{\boldsymbol{a}}(x,z)+q(x)f_{\boldsymbol{a}}(x,z)=-z^{2}%
f_{\boldsymbol{a}}(x,z)\text{.}%
\]
Since $q$ is real valued, the Schr\"{o}dinger operator%
\[
L_{q}=-\partial_{x}^{2}+q
\]
turns to be a symmetric operator on $L^{2}\left(  \mathbb{R}\right)  $, one
can consider the Weyl's classification of the boundaries $\pm\infty$ and the
Weyl functions $m_{\pm}$ if the boundaries are of limit point type. In this
section we establish the connection between the $m$-function and the spectral
theory of $L_{q}$ initiated by Weyl. Necessary knowledge for this section will
be obtained in \cite{m1}.

\begin{lemma}
\label{l20}(Boundary classification) For any $z\in\mathbb{C}\backslash
\mathbb{R}$%
\[
\dim\left\{  f\in L^{2}\left(  \mathbb{R}_{+}\right)  \text{; \ }%
L_{q}f=zf\right\}
\]
is independent of $z$. According to $1$ or $2$ of the dimension the boundary
$+\infty$ is called limit point type or limit circle type respectively. It is
also valid that if $+\infty$ is of limit circle type, then%
\[
\dim\left\{  f\in L^{2}\left(  \mathbb{R}_{+}\right)  \text{; \ }%
L_{q}f=zf\right\}  =2
\]
for any $z\in\mathbb{C}$.\medskip
\end{lemma}

If the boundary $+\infty$ is of limit point type, then there exists a
non-trivial solution $f_{+}\left(  x,z\right)  $ to $L_{q}f_{+}=zf_{+}$ which
is in $L^{2}\left(  \mathbb{R}_{+}\right)  $. $f_{+}$ is unique up to constant
multiple. The Weyl function is defined by%
\[
m_{+}(z)=\dfrac{f_{+}^{\prime}(0,z)}{f_{+}(0,z)}%
\]
The boundary $-\infty$ also has the same classification, and if it is of limit
point type, the Weyl function at $-\infty$ is defined by%
\[
m_{-}(z)=-\dfrac{f_{-}^{\prime}(0,z)}{f_{-}(0,z)}%
\]
with the $L^{2}\left(  \mathbb{R}_{-}\right)  $ non-trivial solution
$f_{-}\left(  x,z\right)  $.

A general sufficient condition for the limit point type which is suitable for
our purpose is known by \cite{h}. A proof is given for completeness sake.

\begin{lemma}
\label{l21}If there exists a positive solution $f$ on $\mathbb{R}_{+}$ to
$L_{q}f=\lambda_{0}f$ for some $\lambda_{0}\in\mathbb{R}$, then the boundary
$+\infty$ is of limit point type.
\end{lemma}

\begin{proof}
Define%
\[
u(x)=f(x)\int_{0}^{x}f(y)^{-2}dy\text{.}%
\]
Then $u$ satisfies $L_{q}u=\lambda_{0}u$. We show $u\notin L^{2}\left(
\mathbb{R}_{+}\right)  $, which implies that $+\infty$ is of limit point type
due to Lemma \ref{l20}. For this purpose set%
\[
s(x)=%
{\displaystyle\int_{0}^{x}}
f\left(  y\right)  ^{-2}dy\text{ \ and \ }t(x)=\left(  -s(x)^{-1}\right)
^{\prime}\text{.}%
\]
Then $t(x)>0$ and%
\[
-s(x)^{-1}=c+%
{\displaystyle\int_{0}^{x}}
t\left(  y\right)  dy
\]
with some constant $c$. Since $s(x)^{-1}>0$, we have%
\[%
{\displaystyle\int_{0}^{x}}
t\left(  y\right)  dy<-c\text{,}%
\]
which implies%
\begin{equation}%
{\displaystyle\int_{0}^{\infty}}
t\left(  y\right)  dy\leq-c<\infty\text{.} \label{75}%
\end{equation}
On the other hand%
\[%
{\displaystyle\int_{0}^{\infty}}
u(x)^{2}dx=%
{\displaystyle\int_{0}^{\infty}}
s^{\prime}(x)^{-1}s(x)^{2}dx=%
{\displaystyle\int_{0}^{\infty}}
\frac{dx}{\left(  -s(x)^{-1}\right)  ^{\prime}}=%
{\displaystyle\int_{0}^{\infty}}
\frac{dx}{t(x)}%
\]
holds. Since for any $t(x)>0$%
\[%
{\displaystyle\int_{0}^{\infty}}
t\left(  y\right)  dy+%
{\displaystyle\int_{0}^{\infty}}
t\left(  y\right)  ^{-1}dy=%
{\displaystyle\int_{0}^{\infty}}
\left(  t\left(  y\right)  +t\left(  y\right)  ^{-1}\right)  dy\geq%
{\displaystyle\int_{0}^{\infty}}
2dy=\infty\text{,}%
\]
(\ref{75}) shows%
\[%
{\displaystyle\int_{0}^{\infty}}
u(x)^{2}dx=%
{\displaystyle\int_{0}^{\infty}}
\frac{dx}{t(x)}=\infty\text{,}%
\]
which completes the proof.\medskip
\end{proof}

Suppose the boundaries $\pm\infty$ are of limit point type. Then, it is known
that the symmetric operator $L_{q}$ has a unique self-adjoint extension in
$L^{2}\left(  \mathbb{R}\right)  $ and $m_{\pm}\left(  z\right)  $ are
analytic functions on $\mathbb{C}\backslash$\textrm{sp}$L_{q}$ satisfying
$\operatorname{Im}m_{\pm}\left(  z\right)  /\operatorname{Im}z>0$.

Now we consider relationship between $m_{\boldsymbol{a}}\left(  z\right)  $
and $m_{\pm}\left(  z\right)  $. The definition of $f_{\boldsymbol{a}}\left(
x,z\right)  $%
\[
f_{\boldsymbol{a}}\left(  x,z\right)  =\boldsymbol{a}\left(  z\right)  \left(
T\left(  e_{x}\boldsymbol{a}\right)  ^{-1}1\right)  \left(  z\right)
\]
implies%
\[
f_{e_{y}\boldsymbol{a}}\left(  x,z\right)  =e_{y}\left(  z\right)
\boldsymbol{a}\left(  z\right)  \left(  T\left(  e_{x}e_{y}\boldsymbol{a}%
\right)  ^{-1}1\right)  \left(  z\right)  =e_{y}\left(  z\right)
f_{\boldsymbol{a}}\left(  x+y,z\right)  \text{.}%
\]
Therefore, Corollary \ref{c1} shows%
\begin{equation}
-m_{e_{y}\boldsymbol{a}}\left(  z\right)  =\left.  \dfrac{\partial_{x}%
f_{e_{y}\boldsymbol{a}}\left(  x,z\right)  }{f_{e_{y}\boldsymbol{a}}\left(
x,z\right)  }\right\vert _{x=0}=\left.  \dfrac{\partial_{x}f_{\boldsymbol{a}%
}\left(  x+y,z\right)  }{f_{\boldsymbol{a}}\left(  x+y,z\right)  }\right\vert
_{x=0}=\dfrac{\partial_{y}f_{\boldsymbol{a}}\left(  y,z\right)  }%
{f_{\boldsymbol{a}}\left(  y,z\right)  }\text{,} \label{81}%
\end{equation}
hence it holds that%
\begin{equation}
f_{\boldsymbol{a}}\left(  x,z\right)  =f_{\boldsymbol{a}}\left(  0,z\right)
\exp\left(  -\int_{0}^{x}m_{e_{y}\boldsymbol{a}}\left(  z\right)  dy\right)
\text{.} \label{79}%
\end{equation}

\begin{proposition}
\label{p4}Let $\boldsymbol{a}\in\boldsymbol{A}_{L,+}^{inv}$ with $L\geq3$ and
$q$ be the associated potential. Then the boundaries $\pm\infty$ are of limit
point type for the Schr\"{o}dinger operator $L_{q}$. The $m$-function
$m_{\boldsymbol{a}}$ and the Weyl functions $m_{\pm}$ are connected with
$m_{\boldsymbol{a}}$ by%
\begin{equation}
m_{\boldsymbol{a}}\left(  z\right)  =\left\{
\begin{array}
[c]{cc}%
-m_{+}\left(  -z^{2}\right)  & \text{if \ }\operatorname{Re}z>0\\
m_{-}\left(  -z^{2}\right)  & \text{if \ }\operatorname{Re}z<0
\end{array}
\right.  \text{.} \label{80}%
\end{equation}

\end{proposition}

\begin{proof}
The key ingredient for the proof is (\ref{81}), which shows%
\begin{equation}
\partial_{x}m_{e_{x}\boldsymbol{a}}\left(  z\right)  =-z^{2}-q(x)+m_{e_{x}%
\boldsymbol{a}}\left(  z\right)  ^{2}\text{,} \label{82}%
\end{equation}
since $f_{\boldsymbol{a}}\left(  x,z\right)  $ satisfies $L_{q}%
f_{\boldsymbol{a}}\left(  x,z\right)  =-z^{2}f_{\boldsymbol{a}}\left(
x,z\right)  $. (\ref{82}) implies%
\[
\partial_{x}\operatorname{Im}m_{e_{x}\boldsymbol{a}}\left(  z\right)
=-\operatorname{Im}z^{2}+2\operatorname{Re}m_{e_{x}\boldsymbol{a}}\left(
z\right)  \operatorname{Im}m_{e_{x}\boldsymbol{a}}\left(  z\right)  \text{,}%
\]
which together with (\ref{79}) yields%
\begin{align*}
\left\vert f_{\boldsymbol{a}}\left(  x,z\right)  \right\vert ^{2}  &
=\left\vert f_{\boldsymbol{a}}\left(  0,z\right)  \right\vert ^{2}\exp\left(
-2\int_{0}^{x}\operatorname{Re}m_{e_{y}\boldsymbol{a}}\left(  z\right)
\right)  dy\\
&  =\left\vert f_{\boldsymbol{a}}\left(  0,z\right)  \right\vert ^{2}%
\dfrac{\operatorname{Im}m_{\boldsymbol{a}}\left(  z\right)  }%
{\operatorname{Im}m_{e_{x}\boldsymbol{a}}\left(  z\right)  }\exp\left(
-\int_{0}^{x}\dfrac{\operatorname{Im}z^{2}}{\operatorname{Im}m_{e_{y}%
\boldsymbol{a}}\left(  z\right)  }dy\right)  \text{.}%
\end{align*}
Then an identity%
\[
\int_{0}^{b}\left\vert f_{\boldsymbol{a}}\left(  x,z\right)  \right\vert
^{2}dx=\left\vert f_{\boldsymbol{a}}\left(  0,z\right)  \right\vert ^{2}%
\dfrac{\operatorname{Im}m_{\boldsymbol{a}}\left(  z\right)  }%
{\operatorname{Im}z^{2}}\left(  1-\exp\left(  -\int_{0}^{b}\dfrac
{\operatorname{Im}z^{2}}{\operatorname{Im}m_{e_{y}\boldsymbol{a}}\left(
z\right)  }dy\right)  \right)
\]
follows. Since $\operatorname{Im}z^{2}=2\operatorname{Re}z\operatorname{Im}z$
and $\operatorname{Im}m_{\boldsymbol{a}}\left(  z\right)  /\operatorname{Im}%
z>0$ hold due to (i) of Lemma \ref{l19}, if $\operatorname{Re}z>0$, we have%
\begin{equation}
\int_{0}^{\infty}\left\vert f_{\boldsymbol{a}}\left(  x,z\right)  \right\vert
^{2}dx\leq\left\vert f_{\boldsymbol{a}}\left(  0,z\right)  \right\vert
^{2}\dfrac{\operatorname{Im}m_{\boldsymbol{a}}\left(  z\right)  }%
{\operatorname{Im}z^{2}}<\infty\text{.} \label{83}%
\end{equation}
On the other hand, if $z=\lambda\in D_{-}\cap\mathbb{R}$, then (\ref{79})
implies $f_{\boldsymbol{a}}\left(  x,\lambda\right)  /f_{\boldsymbol{a}%
}\left(  0,\lambda\right)  $ is a positive solution to $L_{q}f=-\lambda^{2}f$,
hence the boundary $+\infty$ is of limit point type owing to Lemma \ref{l21}.
Since $f_{\boldsymbol{a}}\left(  x,z\right)  \in L^{2}\left(  \mathbb{R}%
_{+}\right)  $ if $\operatorname{Re}z>0$ and $\operatorname{Im}z\neq0$, the
uniqueness of such a solution justifies $f_{\boldsymbol{a}}\left(  x,z\right)
=f_{+}(x,-z^{2})$, which shows the identity $m_{\boldsymbol{a}}\left(
z\right)  =-m_{+}\left(  -z^{2}\right)  $ if $\operatorname{Re}z>0$. The
boundary $-\infty$ can be treated similarly, and we obtain $m_{\boldsymbol{a}%
}\left(  z\right)  =m_{-}\left(  -z^{2}\right)  $ if $\operatorname{Re}z<0$,
which completes the proof.\medskip
\end{proof}

This Proposition says that for $\boldsymbol{a}\in\boldsymbol{A}_{L,+}^{inv}$
its $m$-function $m_{\boldsymbol{a}}\left(  z\right)  $ is analytically
continuable up to $\mathbb{C}\backslash\left(  \left[  -\mu_{0},\mu
_{0}\right]  \cup i\mathbb{R}\right)  $ ($\mu_{0}=\sqrt{-\lambda_{0}}$)
although originally we knew its analyticity only on $D_{-}$.

The next issue is to show the converse statement of Proposition \ref{p4}. This
proposition and Lemma \ref{l19} implies that $m=m_{\boldsymbol{a}}$ for
$\boldsymbol{a}\in\boldsymbol{A}_{L,+}^{inv}\left(  C\right)  $ satisfies%
\begin{equation}
\left\{
\begin{array}
[c]{lcl}%
\dfrac{\operatorname{Im}m\left(  z\right)  }{\operatorname{Im}z}>0 & \text{on}
& \mathbb{C}\backslash\left(  \mathbb{R}\cup i\mathbb{R}\right)  \smallskip\\
\dfrac{m(x)-m(-x)}{x}>0 & \text{if} & x\in\mathbb{R}\text{ and\ }\left\vert
x\right\vert >\mu_{0}%
\end{array}
\right.  \text{.} \label{58}%
\end{equation}
It should be remarked that the analyticity of $m$ on $D_{-}$ implies
$1+\varphi_{\boldsymbol{a}}\left(  z\right)  \neq0$ on $D_{-}$ since
$1+\varphi_{\boldsymbol{a}}\left(  z\right)  $ and $z+\psi_{\boldsymbol{a}%
}\left(  z\right)  $ do not vanish simultaneously due to $\Delta
_{\boldsymbol{a}}\left(  z\right)  \neq0$. In the process of the proof we need
an operation%
\[
\left(  d_{\zeta}f\right)  \left(  z\right)  =\dfrac{z^{2}-\zeta^{2}}{f\left(
z\right)  -f\left(  \zeta\right)  }-f\left(  \zeta\right)
\]
for a function $f$ on $\mathbb{C}$ as long as they have meaning. Then
$\left\{  d_{\zeta}\right\}  _{\zeta\in D_{-}}$ is commutative and $d_{\zeta
}d_{-\zeta}=id$.

\begin{proposition}
\label{p5}Let $L\geq2$. For $\boldsymbol{a}\in\boldsymbol{A}_{L}^{inv}\left(
C\right)  $ suppose that $m_{\boldsymbol{a}}$ is analytic on $\mathbb{C}%
\backslash\left(  \left[  -\mu_{0},\mu_{0}\right]  \cup i\mathbb{R}\right)  $
and satisfies (\ref{58}). Then, $\boldsymbol{a}\in\boldsymbol{A}_{L,+}%
^{inv}\left(  C\right)  $ holds.
\end{proposition}

\begin{proof}
We have to show $\tau_{\boldsymbol{a}}\left(  r\right)  \geq0$ for any real
rational function $r\in\Gamma_{0}^{\left(  0\right)  }$. Since such $r$ is a
product of the 4 types of rational functions of (\ref{86}), first we prove
$\tau_{\boldsymbol{a}}\left(  r\right)  \geq0$ for $r$ of (\ref{86}). If $r$
is of the type (1) with $\eta$, $\zeta\in D_{-}\cap\mathbb{C}_{+}$, then
$\tau_{\boldsymbol{a}}\left(  r\right)  $ is given by (\ref{84}) and
$\tau_{\boldsymbol{a}}\left(  r\right)  >0$ is equivalent to $\left\vert
f(z)\right\vert <1$ on $D_{-}\cap\mathbb{C}_{+}$ with%
\[
f(z)=\dfrac{m_{\boldsymbol{a}}\left(  z\right)  -m_{\boldsymbol{a}}\left(
\zeta\right)  }{m_{\boldsymbol{a}}\left(  z\right)  -\overline
{m_{\boldsymbol{a}}\left(  \zeta\right)  }}\dfrac{z^{2}-\overline{\zeta}^{2}%
}{z^{2}-\zeta^{2}}\text{.}%
\]
Set $w=\zeta^{2}$. Then%
\[
f(\sqrt{z})=\dfrac{m_{\boldsymbol{a}}\left(  \sqrt{z}\right)
-m_{\boldsymbol{a}}\left(  \sqrt{w}\right)  }{m_{\boldsymbol{a}}\left(
\sqrt{z}\right)  -\overline{m_{\boldsymbol{a}}\left(  \sqrt{w}\right)  }%
}\dfrac{z-\overline{w}}{z-w}%
\]
holds. Since $m_{\boldsymbol{a}}\left(  \sqrt{z}\right)  $ is analytic on
$\mathbb{C}_{+}$ and $\operatorname{Im}m_{\boldsymbol{a}}\left(  \sqrt
{z}\right)  >0$, Schwarz lemma implies $\left\vert f(\sqrt{z})\right\vert <1$
for $z$, $\zeta\in\mathbb{C}_{+}$, unless $m_{\boldsymbol{a}}\left(  \sqrt
{z}\right)  $ is%
\[
m_{\boldsymbol{a}}\left(  \sqrt{z}\right)  =\dfrac{az+b}{cz+d}%
\]
with some constants $a$, $b$, $c$, $d$ satisfying $ad-bc\neq0$, which is
impossible since $m_{\boldsymbol{a}}\left(  \sqrt{z}\right)  =\sqrt
{z}+o\left(  1\right)  $ as $z\rightarrow\infty$. Therefore we have
$\left\vert f(z)\right\vert <1$ if $\operatorname{Re}z$, $\operatorname{Im}%
z>0$, $\operatorname{Re}\zeta$, $\operatorname{Im}\zeta>0$. On the other hand,
when $z\in\mathbb{C}_{-}$ we use the identity%
\[
f(-\sqrt{z})=\dfrac{m_{\boldsymbol{a}}\left(  -\sqrt{z}\right)
-m_{\boldsymbol{a}}\left(  \sqrt{w}\right)  }{m_{\boldsymbol{a}}\left(
-\sqrt{z}\right)  -\overline{m_{\boldsymbol{a}}\left(  \sqrt{w}\right)  }%
}\dfrac{z-\overline{w}}{z-w}\text{.}%
\]
If $\operatorname{Re}\zeta>0$, $\operatorname{Im}\zeta>0$, then
$\operatorname{Im}w>0$ implies%
\[
\left\vert \dfrac{z-\overline{w}}{z-w}\right\vert <1
\]
and%
\[
\operatorname{Im}m_{\boldsymbol{a}}\left(  -\sqrt{z}\right)  >0\text{,
\ }\operatorname{Im}m_{\boldsymbol{a}}\left(  \sqrt{w}\right)  >0
\]
due to $\operatorname{Im}\sqrt{z}<0$, $\operatorname{Im}\sqrt{w}>0$, hence%
\[
\left\vert \dfrac{m_{\boldsymbol{a}}\left(  -\sqrt{z}\right)
-m_{\boldsymbol{a}}\left(  \sqrt{w}\right)  }{m_{\boldsymbol{a}}\left(
-\sqrt{z}\right)  -\overline{m_{\boldsymbol{a}}\left(  \sqrt{w}\right)  }%
}\dfrac{z-\overline{w}}{z-w}\right\vert <\left\vert \dfrac{m_{\boldsymbol{a}%
}\left(  -\sqrt{z}\right)  -m_{\boldsymbol{a}}\left(  \sqrt{w}\right)
}{m_{\boldsymbol{a}}\left(  -\sqrt{z}\right)  -\overline{m_{\boldsymbol{a}%
}\left(  \sqrt{w}\right)  }}\right\vert <1\text{,}%
\]
which implies $\left\vert f(z)\right\vert <1$ if $\operatorname{Re}z<0$,
$\operatorname{Im}z>0$, $\operatorname{Re}\zeta$, $\operatorname{Im}\zeta>0$.
The rest of the cases can be proved by the symmetry and we have $\tau
_{\boldsymbol{a}}\left(  r\right)  >0$.

For the type (2) $r=q_{s}p_{t}$%
\[
\tau_{\boldsymbol{a}}\left(  r\right)  =\dfrac{\left(  1+\varphi
_{\boldsymbol{a}}\left(  s\right)  \right)  \left(  1+\varphi_{\boldsymbol{a}%
}\left(  -t\right)  \right)  }{\Delta_{\boldsymbol{a}}\left(  t\right)
}\dfrac{m_{\boldsymbol{a}}\left(  s\right)  -m_{\boldsymbol{a}}\left(
t\right)  }{s-t}\text{.}%
\]
Recall%
\[
\dfrac{\left(  1+\varphi_{\boldsymbol{a}}\left(  s\right)  \right)  \left(
1+\varphi_{\boldsymbol{a}}\left(  -t\right)  \right)  }{\Delta_{\boldsymbol{a}%
}\left(  t\right)  }>0\text{ \ if \ }\left\vert s\right\vert \text{,
}\left\vert t\right\vert >\mu_{0}\text{.}%
\]
On the other hand the property $\operatorname{Im}m_{\boldsymbol{a}}\left(
z\right)  /\operatorname{Im}z>0$ implies%
\[
m_{\boldsymbol{a}}^{\prime}\left(  t\right)  =\lim_{\epsilon\rightarrow0}%
\frac{m_{\boldsymbol{a}}\left(  t+i\epsilon\right)  -m_{\boldsymbol{a}}\left(
t-i\epsilon\right)  }{2i\epsilon}=\lim_{\epsilon\rightarrow0}\frac
{\operatorname{Im}m_{\boldsymbol{a}}\left(  t+i\epsilon\right)  }%
{\operatorname{Im}\left(  t+i\epsilon\right)  }\geq0\text{,}%
\]
which shows%
\[
\dfrac{m_{\boldsymbol{a}}\left(  s\right)  -m_{\boldsymbol{a}}\left(
t\right)  }{s-t}\geq0\text{ \ if \ }s,t\in\left(  -\infty,-\mu_{0}\right)
\text{ or }\left(  \mu_{0},\infty\right)  \text{.}%
\]
If $s\in\left(  -\infty,-\mu_{0}\right)  $ and $,t\in\left(  \mu_{0}%
,\infty\right)  $, then from $\left(  m\left(  x\right)  -m\left(  -x\right)
\right)  /x>0$ an inequality%
\[
m_{\boldsymbol{a}}\left(  s\right)  -m_{\boldsymbol{a}}\left(  t\right)
\leq0
\]
follows, hence we have $\tau_{\boldsymbol{a}}\left(  r\right)  \geq0$ for the
case (2). If $r=q_{\zeta}q_{\overline{\zeta}}p_{s}p_{t}$, the cocycle property
implies%
\[
\tau_{\boldsymbol{a}}\left(  r\right)  =\tau_{\boldsymbol{a}}\left(  q_{\zeta
}p_{\overline{\zeta}}p_{\overline{\zeta}}^{-1}q_{\overline{\zeta}}p_{s}%
p_{t}\right)  =\tau_{\boldsymbol{a}}\left(  q_{\zeta}p_{\overline{\zeta}%
}\right)  \tau_{q_{\zeta}p_{\overline{\zeta}}\boldsymbol{a}}\left(
p_{\overline{\zeta}}^{-1}q_{\overline{\zeta}}p_{s}p_{t}\right)  \text{.}%
\]
Since the $m$-function for $q_{\zeta}p_{\overline{\zeta}}\boldsymbol{a}$ is
$d_{\zeta}d_{\overline{\zeta}}m_{\boldsymbol{a}}$, which satisfies (\ref{58})
due to Lemma \ref{l25}, we have $\tau_{q_{\zeta}p_{\overline{\zeta}%
}\boldsymbol{a}}\left(  p_{s}p_{t}\right)  \neq0$ in view of the last
argument, and%
\[
\tau_{q_{\zeta}p_{\overline{\zeta}}\boldsymbol{a}}\left(  p_{\overline{\zeta}%
}^{-1}q_{\overline{\zeta}}p_{s}p_{t}\right)  =\tau_{q_{\zeta}p_{\overline
{\zeta}}\boldsymbol{a}}\left(  p_{\overline{\zeta}}^{-1}q_{\overline{\zeta}%
}\right)  \tau_{q_{\zeta}p_{\overline{\zeta}}\boldsymbol{a}}\left(  p_{s}%
p_{t}\right)  =\Delta_{q_{\zeta}p_{\overline{\zeta}}\boldsymbol{a}}\left(
\overline{\zeta}\right)  \tau_{q_{\zeta}p_{\overline{\zeta}}\boldsymbol{a}%
}\left(  p_{s}p_{t}\right)  \neq0
\]
is valid. Therefore we have $\tau_{\boldsymbol{a}}\left(  r\right)  >0$. The
case (4) $r=q_{s}q_{t}p_{\eta}p_{\overline{\eta}}$ can be treated similarly,
hence $\tau_{\boldsymbol{a}}\left(  r\right)  >0$ for $r$ of any type of
(\ref{86}).

The property $\tau_{\boldsymbol{a}}\left(  r\right)  \geq0$ for general real
$r\in\Gamma_{0}^{\left(  0\right)  }$ can be obtained by observing%
\[
\tau_{\boldsymbol{a}}\left(  r_{1}r_{2}\right)  =\tau_{\boldsymbol{a}}\left(
r_{1}\right)  \tau_{r_{1}\boldsymbol{a}}\left(  r_{2}\right)  \geq0
\]
if $\tau_{\boldsymbol{a}}\left(  r_{1}\right)  >0$ and the $m$-function
$m_{r_{1}\boldsymbol{a}}$ satisfies (\ref{58}) since $m_{r_{1}\boldsymbol{a}}$
is obtained by repeating the operation $d_{\zeta}d_{\overline{\zeta}}$,
$d_{t}$ to $m_{\boldsymbol{a}}$. If $\tau_{\boldsymbol{a}}\left(
r_{1}\right)  =0$, approximating $r_{1}$ by real rational functions $r_{n}$
such that $\tau_{\boldsymbol{a}}\left(  r_{n}\right)  >0$ one sees
$\tau_{\boldsymbol{a}}\left(  r_{1}r_{2}\right)  \geq0$, which completes the
proof.\medskip
\end{proof}

It is certainly better to prove $\tau_{\boldsymbol{a}}\left(  r\right)  >0$
directly by showing%
\[
\det\left(  \dfrac{m_{\boldsymbol{a}}\left(  \zeta_{i}\right)
-m_{\boldsymbol{a}}\left(  -\eta_{j}\right)  }{\zeta_{i}^{2}-\eta_{j}^{2}%
}\right)  \neq0\text{ \ (see Lemma \ref{l13})}%
\]
for $m_{\boldsymbol{a}}$ satisfying (\ref{58}), however the author has no such proof.

\subsection{Proof of Theorem \ref{t1}}

A more concrete criterion for an $m$ to be in $\boldsymbol{A}_{L,+}^{inv}$ is
given by the two conditions (M.1), (M.2) in Introduction. Recall the
definitions. Suppose an analytic function $m$ on $\mathbb{C}\backslash\left(
\left[  -\mu_{0},\mu_{0}\right]  \cup i\mathbb{R}\right)  $ ($\mu_{0}%
=\sqrt{-\lambda_{0}}$) satisfies\medskip\newline(M.1) \ $m\left(  z\right)  $
has the property (\ref{58}), namely%
\begin{equation}
\left\{
\begin{array}
[c]{lcl}%
\dfrac{\operatorname{Im}m\left(  z\right)  }{\operatorname{Im}z}>0 & \text{on}
& \mathbb{C}\backslash\left(  \mathbb{R}\cup i\mathbb{R}\right)  \smallskip\\
\dfrac{m(x)-m(-x)}{x}>0 & \text{if} & x\in\mathbb{R}\text{ and\ }\left\vert
x\right\vert >\mu_{0}%
\end{array}
\right.  \text{.} \label{42}%
\end{equation}
\newline(M.2) \ $m$ has an asymptotic behavior:%
\begin{equation}
m\left(  z\right)  =z+\sum_{1\leq k\leq L-2}m_{k}z^{-k}+O\left(
z^{-L+1}\right)  \text{ on }D_{-}\text{.} \label{77}%
\end{equation}
\newline

\textbf{Proof of Theorem \ref{t1}}

Suppose $m$ satisfies (M.1), (M.2) and set $\boldsymbol{m}\left(  z\right)
\equiv\left(  1,m\left(  z\right)  /z\right)  $. To see $\boldsymbol{m}%
\in\boldsymbol{M}_{L}\left(  C\right)  $ only the condition (ii) of (\ref{9}),
namely%
\begin{equation}
m_{1}(z)m_{2}(-z)+m_{1}(-z)m_{2}(z)\neq0\text{\ on\ }\mathbb{C}\backslash
\left(  \left[  -\mu_{0},\mu_{0}\right]  \cup i\mathbb{R}\right)  \label{85}%
\end{equation}
must be verified. For this $\boldsymbol{m}\left(  z\right)  $ the left hand
side of (\ref{85}) is $\left(  m\left(  z\right)  -m\left(  -z\right)
\right)  /z$, which is not $0$ since $\operatorname{Im}\left(  m\left(
z\right)  -m\left(  -z\right)  \right)  >0$ if $\operatorname{Im}z>0$ and
$m(x)-m(-x)\neq0$ if $\left\vert x\right\vert >\mu_{0}$. Since the identity
$m=m_{\boldsymbol{m}}$ is clear, Proposition \ref{p5} implies $\boldsymbol{m}%
\in\boldsymbol{A}_{L,+}^{inv}\left(  C\right)  $, which proves the
theorem.\medskip

It may be interesting to see to what extent $m_{\boldsymbol{a}}$ for
$\boldsymbol{a}\in\boldsymbol{A}_{L,+}^{inv}\left(  C\right)  $ has the
property (M.2).

\begin{proposition}
\label{p9}Let $C^{\prime}=\sigma C$ with $\sigma>1$. Then $m_{\boldsymbol{a}}$
for $\boldsymbol{a}\in\boldsymbol{A}_{L,+}^{inv}\left(  C\right)  $ satisfies
(M.2) on $D_{-}^{\prime}=\sigma D_{-}$ replacing $L$ by $L-1-\left(
n-1\right)  /2$.
\end{proposition}

\begin{proof}
To verify the property (M.2) recall the definition of the $m$-function
$m_{\boldsymbol{a}}$ with $\boldsymbol{a}\in\boldsymbol{A}_{L,+}^{inv}\left(
C\right)  $:%
\[
m_{\boldsymbol{a}}\left(  z\right)  =\dfrac{z+\psi_{\boldsymbol{a}}\left(
z\right)  }{1+\varphi_{\boldsymbol{a}}\left(  z\right)  }+\kappa_{1}\left(
\boldsymbol{a}\right)
\]
with%
\[
\left\{
\begin{array}
[c]{l}%
\varphi_{\boldsymbol{a}}\left(  z\right)  =\boldsymbol{a}\left(  z\right)
T\left(  \boldsymbol{a}\right)  ^{-1}1-1\\
\psi_{\boldsymbol{a}}\left(  z\right)  =\boldsymbol{a}\left(  z\right)
T\left(  \boldsymbol{a}\right)  ^{-1}z-z
\end{array}
\in H\left(  D_{-}\right)  \right.  \text{,}%
\]
and $\kappa_{1}\left(  \boldsymbol{a}\right)  =\lim_{z\rightarrow\infty
}z\varphi_{\boldsymbol{a}}\left(  z\right)  $. Set $u=T\left(  \boldsymbol{a}%
\right)  ^{-1}1\in H_{1}\left(  D_{+}\right)  $, $v=T\left(  \boldsymbol{a}%
\right)  ^{-1}z\in H_{2}\left(  D_{+}\right)  $. Since (\ref{e6}) implies%
\[
\mathfrak{p}_{-}\left(  \boldsymbol{a}u\right)  \left(  z\right)  =\dfrac
{1}{2\pi i}%
{\displaystyle\int_{C}}
\dfrac{\left(  \boldsymbol{a}\left(  \lambda\right)  -\boldsymbol{f}\left(
\lambda\right)  \right)  u\left(  \lambda\right)  }{z-\lambda}d\lambda\in
H\left(  D_{-}\right)
\]
for $u\in H_{L}\left(  D_{+}\right)  $, we have%
\[
\left\{
\begin{array}
[c]{c}%
\varphi_{\boldsymbol{a}}\left(  z\right)  =\mathfrak{p}_{-}\left(
\boldsymbol{a}u\right)  \left(  z\right)  =\dfrac{1}{2\pi i}%
{\displaystyle\int_{C}}
\dfrac{\left(  \boldsymbol{a}\left(  \lambda\right)  -\boldsymbol{f}\left(
\lambda\right)  \right)  u\left(  \lambda\right)  }{z-\lambda}d\lambda
\smallskip\\
\psi_{\boldsymbol{a}}\left(  z\right)  =\mathfrak{p}_{-}\left(  \boldsymbol{a}%
v\right)  \left(  z\right)  =\dfrac{1}{2\pi i}%
{\displaystyle\int_{C}}
\dfrac{\left(  \boldsymbol{a}\left(  \lambda\right)  -\boldsymbol{f}\left(
\lambda\right)  \right)  v\left(  \lambda\right)  }{z-\lambda}d\lambda
\end{array}
\right.  \text{.}%
\]
The expansion%
\[
\dfrac{1}{z-\lambda}=\sum_{k=1}^{M}z^{-k}\lambda^{k-1}+z^{-M}\dfrac
{\lambda^{M}}{z-\lambda}%
\]
implies%
\[
\left\{
\begin{array}
[c]{l}%
\varphi_{\boldsymbol{a}}\left(  z\right)  =\sum_{k=1}^{L-1}z^{-k}\ell
_{k-1}\left(  u\right)  +z^{-L+1}\delta_{1}\left(  z\right) \\
\psi_{\boldsymbol{a}}\left(  z\right)  =\sum_{k=1}^{L-2}z^{-k}\ell
_{k-1}\left(  v\right)  +z^{-L+2}\delta_{2}\left(  z\right)
\end{array}
\right.
\]
with%
\[
\left\{
\begin{array}
[c]{l}%
\ell_{k}\left(  u\right)  =\dfrac{1}{2\pi i}%
{\displaystyle\int_{C}}
\lambda^{k}\left(  \boldsymbol{a}\left(  \lambda\right)  -\boldsymbol{f}%
\left(  \lambda\right)  \right)  u\left(  \lambda\right)  d\lambda\\
\delta_{1}\left(  z\right)  =\dfrac{1}{2\pi i}%
{\displaystyle\int_{C}}
\dfrac{\lambda^{L-1}\left(  \boldsymbol{a}\left(  \lambda\right)
-\boldsymbol{f}\left(  \lambda\right)  \right)  u\left(  \lambda\right)
}{z-\lambda}d\lambda\in H\left(  D_{-}\right) \\
\delta_{2}\left(  z\right)  =\dfrac{1}{2\pi i}%
{\displaystyle\int_{C}}
\dfrac{\lambda^{L-2}\left(  \boldsymbol{a}\left(  \lambda\right)
-\boldsymbol{f}\left(  \lambda\right)  \right)  v\left(  \lambda\right)
}{z-\lambda}d\lambda\in H\left(  D_{-}\right)
\end{array}
\right.  \text{.}%
\]
Since $\lambda^{L-1}\left(  \boldsymbol{a}\left(  \lambda\right)
-\boldsymbol{f}\left(  \lambda\right)  \right)  u\left(  \lambda\right)  \in
L^{2}\left(  C\right)  $, Schwarz inequality implies%
\[
\left\vert \delta_{1}\left(  z\right)  \right\vert ^{2}\leq\dfrac{1}{4\pi^{2}}%
{\displaystyle\int_{C}}
\dfrac{\left\vert d\lambda\right\vert }{\left\vert z-\lambda\right\vert ^{2}}%
{\displaystyle\int_{C}}
\left\vert \lambda^{L-1}\left(  \boldsymbol{a}\left(  \lambda\right)
-\boldsymbol{f}\left(  \lambda\right)  \right)  u\left(  \lambda\right)
\right\vert ^{2}\left\vert d\lambda\right\vert \leq c\left\vert z\right\vert
^{n-1}%
\]
\ for $z\in C^{\prime}$ due to Lemma \ref{l15}. Similarly we have $\left\vert
\delta_{2}\left(  z\right)  \right\vert ^{2}\leq c\left\vert z\right\vert
^{n-1}$. Therefore%
\[
\left\{
\begin{array}
[c]{l}%
\varphi_{\boldsymbol{a}}\left(  z\right)  =\sum_{k=1}^{L^{\prime}-1}z^{-k}%
\ell_{k-1}\left(  u\right)  +O\left(  z^{-L^{\prime}}\right)  \smallskip\\
\psi_{\boldsymbol{a}}\left(  z\right)  =\sum_{k=1}^{L^{\prime}-2}z^{-k}%
\ell_{k-1}\left(  v\right)  +O\left(  z^{-L^{\prime}+1}\right)
\end{array}
\right.
\]
holds on $C^{\prime}$ with $L^{\prime}=L-1-\left(  n-1\right)  /2$, which
implies%
\[
m_{\boldsymbol{a}}\left(  z\right)  =z+\sum_{k=1}^{L^{\prime}-2}c_{k}%
z^{-k}+O\left(  z^{-L^{\prime}+1}\right)  \text{.}%
\]

\end{proof}

\section{KdV flow}

We are ready to construct the KdV flow.

\subsection{Definition of the flow and Theorem \ref{t2}}

We define the KdV flow by making use of the $m$-functions and the continuity
of $m_{\boldsymbol{a}}\left(  \zeta\right)  $ with respect to $\boldsymbol{a}$
is necessary, which is shown by representing $m_{\boldsymbol{a}}^{\prime}$ by
the tau-function. The identities in (\ref{66}) implies%
\begin{align*}
m_{\boldsymbol{a}}^{\prime}\left(  \zeta\right)   &  =\frac{\tau
_{\boldsymbol{a}}\left(  q_{\zeta}^{2}\right)  }{\tau_{\boldsymbol{a}}\left(
q_{\zeta}\right)  ^{2}}\\
&  =\frac{\tau_{\boldsymbol{a}}^{\left(  2\right)  }\left(  q_{\zeta}%
^{2}\right)  }{\tau_{\boldsymbol{a}}^{\left(  2\right)  }\left(  q_{\zeta
}\right)  ^{2}}\exp\mathrm{tr}\left(  \left(  q_{\zeta}^{-2}T\left(  q_{\zeta
}^{2}\boldsymbol{a}\right)  -2q_{\zeta}^{-1}T\left(  q_{\zeta}\boldsymbol{a}%
\right)  +T\left(  \boldsymbol{a}\right)  \right)  T\left(  \boldsymbol{a}%
\right)  ^{-1}\right)  \text{.}%
\end{align*}
The cocycle property shows%
\begin{align*}
m_{g\boldsymbol{a}}^{\prime}\left(  \zeta\right)   &  =\frac{\tau
_{g\boldsymbol{a}}^{\left(  2\right)  }\left(  q_{\zeta}^{2}\right)  }%
{\tau_{g\boldsymbol{a}}^{\left(  2\right)  }\left(  q_{\zeta}\right)  ^{2}%
}\exp\mathrm{tr}\left(  \left(  q_{\zeta}^{-2}T\left(  q_{\zeta}%
^{2}g\boldsymbol{a}\right)  -2q_{\zeta}^{-1}T\left(  q_{\zeta}g\boldsymbol{a}%
\right)  +T\left(  g\boldsymbol{a}\right)  \right)  T\left(  g\boldsymbol{a}%
\right)  ^{-1}\right) \\
&  =\frac{\tau_{\boldsymbol{a}}^{\left(  2\right)  }\left(  gq_{\zeta}%
^{2}\right)  }{\tau_{\boldsymbol{a}}^{\left(  2\right)  }\left(  gq_{\zeta
}\right)  ^{2}}\tau_{\boldsymbol{a}}^{\left(  2\right)  }\left(  g\right)
\exp\mathrm{tr}B_{1}%
\end{align*}
for $g\in\Gamma_{n}$ with%
\begin{align*}
B_{1}  &  =q_{\zeta}^{-2}T\left(  gq_{\zeta}^{2}\boldsymbol{a}\right)
T\left(  g\boldsymbol{a}\right)  ^{-1}-2q_{\zeta}^{-1}T\left(  gq_{\zeta
}\boldsymbol{a}\right)  T\left(  g\boldsymbol{a}\right)  ^{-1}+I\\
&  +\left(  \left(  gq_{\zeta}^{2}\right)  ^{-1}T\left(  gq_{\zeta}%
^{2}\boldsymbol{a}\right)  T\left(  g\boldsymbol{a}\right)  ^{-1}g-I\right)
\left(  g^{-1}T\left(  g\boldsymbol{a}\right)  T\left(  \boldsymbol{a}\right)
^{-1}-I\right) \\
&  -2\left(  \left(  gq_{\zeta}\right)  ^{-1}T\left(  gq_{\zeta}%
\boldsymbol{a}\right)  T\left(  g\boldsymbol{a}\right)  ^{-1}g-I\right)
\left(  g^{-1}T\left(  g\boldsymbol{a}\right)  T\left(  \boldsymbol{a}\right)
^{-1}-I\right) \\
&  =B_{2}-g^{-1}B_{2}g+g^{-1}B_{2}T\left(  g\boldsymbol{a}\right)  T\left(
\boldsymbol{a}\right)  ^{-1}\text{,}%
\end{align*}
where%
\[
B_{2}=q_{\zeta}^{-2}T\left(  gq_{\zeta}^{2}\boldsymbol{a}\right)  T\left(
g\boldsymbol{a}\right)  ^{-1}-2q_{\zeta}^{-1}T\left(  gq_{\zeta}%
\boldsymbol{a}\right)  T\left(  g\boldsymbol{a}\right)  ^{-1}+I\text{.}%
\]
Then we have%
\[
\mathrm{tr}B_{1}=\mathrm{tr}g^{-1}\left(  q_{\zeta}^{-2}T\left(  gq_{\zeta
}^{2}\boldsymbol{a}\right)  -2q_{\zeta}^{-1}T\left(  gq_{\zeta}\boldsymbol{a}%
\right)  +T\left(  g\boldsymbol{a}\right)  \right)  T\left(  \boldsymbol{a}%
\right)  ^{-1}\text{.}%
\]
Since for $v\in H_{N}\left(  D_{+}\right)  $%
\begin{align*}
&  q_{\zeta}^{-2}T\left(  gq_{\zeta}^{2}\boldsymbol{a}\right)  -2q_{\zeta
}^{-1}T\left(  gq_{\zeta}\boldsymbol{a}\right)  +T\left(  g\boldsymbol{a}%
\right)  v\left(  z\right) \\
&  =\frac{1}{2\pi i}\int_{C}\frac{g\left(  \lambda\right)  \left(  q_{\zeta
}\left(  z\right)  ^{-1}q_{\zeta}\left(  \lambda\right)  -1\right)
^{2}\widetilde{\boldsymbol{a}}\left(  \lambda\right)  v\left(  \lambda\right)
}{\lambda-z}d\lambda\\
&  =\frac{1}{2\pi i}\int_{C}\frac{\left(  \lambda-z\right)  g\left(
\lambda\right)  \widetilde{\boldsymbol{a}}\left(  \lambda\right)  v\left(
\lambda\right)  }{\left(  \zeta-\lambda\right)  ^{2}}d\lambda\text{,}%
\end{align*}
the above operator is of rank $2$, and the trace turns to%
\[
\mathrm{tr}B_{1}=\frac{1}{2\pi i}\int_{C}\frac{\Theta\left(  g,\lambda\right)
}{\left(  \zeta-\lambda\right)  ^{2}}d\lambda
\]
with%
\[
\Theta\left(  g,\lambda\right)  =g\left(  \lambda\right)  \left(
\lambda\widetilde{\boldsymbol{a}}\left(  \lambda\right)  \left(  T\left(
\boldsymbol{a}\right)  ^{-1}g^{-1}\right)  \left(  \lambda\right)
-\widetilde{\boldsymbol{a}}\left(  \lambda\right)  \left(  T\left(
\boldsymbol{a}\right)  ^{-1}g^{-1}z\right)  \left(  \lambda\right)  \right)
\text{,}%
\]
hence%
\begin{equation}
m_{g\boldsymbol{a}}^{\prime}\left(  \zeta\right)  =\frac{\tau_{\boldsymbol{a}%
}^{\left(  2\right)  }\left(  gq_{\zeta}^{2}\right)  }{\tau_{\boldsymbol{a}%
}^{\left(  2\right)  }\left(  gq_{\zeta}\right)  ^{2}}\tau_{\boldsymbol{a}%
}^{\left(  2\right)  }\left(  g\right)  \exp\left(  \frac{1}{2\pi i}\int
_{C}\frac{\Theta\left(  g,\lambda\right)  }{\left(  \zeta-\lambda\right)
^{2}}d\lambda\right)  \label{106}%
\end{equation}
holds.

\begin{lemma}
\label{l22}For $L\geq\max\left\{  n+1,3\right\}  $ let $\boldsymbol{a}_{1}$,
$\boldsymbol{a}_{2}\in\boldsymbol{A}_{L,+}^{inv}\left(  C\right)  $.
Then\newline(i) \ Suppose $m_{\boldsymbol{a}_{1}}=m_{\boldsymbol{a}_{2}}$.
Then $m_{g\boldsymbol{a}_{1}}=m_{g\boldsymbol{a}_{2}}$ for any $g\in\Gamma
_{n}$.\newline(ii) Suppose $\partial_{x}\kappa_{1}\left(  e_{x}\boldsymbol{a}%
_{1}\right)  =\partial_{x}\kappa_{1}\left(  e_{x}\boldsymbol{a}_{2}\right)  $
for any $x\in\mathbb{R}$. Then $\partial_{x}\kappa_{1}\left(  e_{x}%
g\boldsymbol{a}_{1}\right)  =\partial_{x}\kappa_{1}\left(  e_{x}%
g\boldsymbol{a}_{2}\right)  $ for any $g\in\Gamma_{n}$ and $x\in\mathbb{R}$.
\end{lemma}

\begin{proof}
Propositions \ref{p6}, \ref{p4}, \ref{p5} provide necessary ingredients.
Suppose $m_{\boldsymbol{a}_{1}}=m_{\boldsymbol{a}_{2}}$ for $\boldsymbol{a}%
_{1}$, $\boldsymbol{a}_{2}\in\boldsymbol{A}_{L,+}^{inv}\left(  C\right)  $.
Then Lemma \ref{l16} implies%
\[
m_{q_{\zeta}p_{\eta}\boldsymbol{a}_{1}}\left(  z\right)  =\left(  d_{\zeta
}d_{\eta}m_{\boldsymbol{a}_{1}}\right)  \left(  z\right)  =\left(  d_{\zeta
}d_{\eta}m_{\boldsymbol{a}_{2}}\right)  \left(  z\right)  =m_{q_{\zeta}%
p_{\eta}\boldsymbol{a}_{2}}\left(  z\right)  \text{.}%
\]
Repeating this operation finite times one can show $m_{r\boldsymbol{a}_{1}%
}=m_{r\boldsymbol{a}_{2}}$ for any real rational function $r\in\Gamma
_{0}^{\left(  0\right)  }$. For $g\in\Gamma_{n}$ we approximate it by real
rational functions $r_{k}\in\Gamma_{0}^{\left(  0\right)  }$, which is
possible by Lemma \ref{l18}. We show the convergence $m_{r_{k}\boldsymbol{a}%
}^{\prime}\left(  \zeta\right)  \rightarrow m_{g\boldsymbol{a}}^{\prime
}\left(  \zeta\right)  $ for each fixed $\zeta\in D_{-}$ by making use of
(\ref{106}). Lemma \ref{l11} shows%
\[
\lim_{k\rightarrow\infty}\frac{\tau_{\boldsymbol{a}}^{\left(  2\right)
}\left(  r_{k}q_{\zeta}^{2}\right)  }{\tau_{\boldsymbol{a}}^{\left(  2\right)
}\left(  r_{k}q_{\zeta}\right)  ^{2}}=\frac{\tau_{\boldsymbol{a}}^{\left(
2\right)  }\left(  gq_{\zeta}^{2}\right)  }{\tau_{\boldsymbol{a}}^{\left(
2\right)  }\left(  gq_{\zeta}\right)  ^{2}}\text{.}%
\]
On the other hand $\Theta\left(  r_{k},\lambda\right)  \rightarrow
\Theta\left(  g,\lambda\right)  $ in $L^{2}\left(  C\right)  $ is valid if
$L\geq\max\left\{  n+1,3\right\}  $, hence%
\begin{equation}
\lim_{k\rightarrow\infty}m_{r_{k}\boldsymbol{a}}^{\prime}\left(  \zeta\right)
=m_{g\boldsymbol{a}}^{\prime}\left(  \zeta\right)  \label{110}%
\end{equation}
follows. Consequently one sees%
\[
m_{g\boldsymbol{a}_{1}}^{\prime}\left(  \zeta\right)  =\lim_{k\rightarrow
\infty}m_{r_{k}\boldsymbol{a}_{1}}^{\prime}\left(  \zeta\right)
=\lim_{k\rightarrow\infty}m_{r_{k}\boldsymbol{a}_{2}}^{\prime}\left(
\zeta\right)  =m_{g\boldsymbol{a}_{2}}^{\prime}\left(  \zeta\right)
\]
for $\zeta\in D_{-}$. Since generally $m_{g\boldsymbol{a}}\left(
\zeta\right)  =\zeta+o\left(  1\right)  $, we have $m_{g\boldsymbol{a}_{1}%
}=m_{g\boldsymbol{a}_{2}}$.

(ii) is proved by the uniqueness of the correspondence between the Weyl
functions $m_{\pm}$ and the potential $q$. That is, $m_{\boldsymbol{a}_{1}%
}=m_{\boldsymbol{a}_{2}}$ follows from%
\[
-2\partial_{x}\kappa_{1}\left(  e_{x}\boldsymbol{a}_{1}\right)  =-2\partial
_{x}\kappa_{1}\left(  e_{x}\boldsymbol{a}_{2}\right)  =q(x)\text{.}%
\]
Then (i) yields $m_{g\boldsymbol{a}_{1}}=m_{g\boldsymbol{a}_{2}}$, which
implies again by the uniqueness%
\[
-2\partial_{x}\kappa_{1}\left(  e_{x}g\boldsymbol{a}_{1}\right)
=-2\partial_{x}\kappa_{1}\left(  e_{x}g\boldsymbol{a}_{2}\right)  \text{.}%
\]
We have used the condition $L\geq\max\left\{  n+1,3\right\}  $ to have the
differentiability.\medskip
\end{proof}

Set%
\[
\mathcal{Q}_{L}\left(  C\right)  =\left\{  q\text{; }q(x)=-2\partial_{x}%
\kappa_{1}\left(  e_{x}\boldsymbol{a}\right)  \text{ with real }%
\boldsymbol{a}\in A_{L,+}^{inv}\left(  C\right)  \right\}  \text{.}%
\]
Then this Lemma allows us to define%
\[
\left(  K\left(  g\right)  q\right)  (x)=-2\partial_{x}\kappa_{1}\left(
e_{x}g\boldsymbol{a}\right)  \text{ \ if }q=-2\partial_{x}\kappa_{1}\left(
e_{x}\boldsymbol{a}\right)  \in\mathcal{Q}_{L}\left(  C\right)
\]
for\ $g\in\Gamma_{n}$, and we have Theorem \ref{t2}. We call the flow
$\left\{  K\left(  g\right)  \right\}  _{g\in\Gamma_{n}}$ as \textbf{KdV flow}.

It might be helpful to remark that one can define an equivalent flow on the
space of $m$-functions. Let%
\[
\mathcal{M}_{L}\left(  C\right)  =\left\{  m\text{; }m(z)=m_{\boldsymbol{a}%
}(z)\text{ for }\boldsymbol{a}\in A_{L,+}^{inv}\left(  C\right)  \right\}  ,
\]
and define%
\[
g\cdot m_{\boldsymbol{a}}=m_{g\boldsymbol{a}}\text{ \ \ for }m_{\boldsymbol{a}%
}\in\mathcal{M}_{L}\left(  C\right)  \text{.}%
\]
(i) of Lemma \ref{l12} justifies this definition. The set%
\[
\mathcal{M}_{\infty}\left(  C\right)  =%
{\displaystyle\bigcap\limits_{L\geq1}}
\left(
{\displaystyle\bigcup\limits_{\sigma>1}}
\mathcal{M}_{L}\left(  \sigma C\right)  \right)
\]
is equal to the set of all functions $m$ satisfying (M.1), (M.2) on $\sigma
D_{-}$. Theorem \ref{t1} implies $g\cdot m\in\mathcal{M}_{\infty}\left(
C\right)  $ for $m\in\mathcal{M}_{\infty}\left(  C\right)  $.

\subsection{Tau-function representation of the flow}

Hirota introduced his tau-function as the function $u(t,x)$ such that
$-2\partial_{x}^{2}\log u(t,x)$ is a solution to the KdV equation and he tried
to find an equation satisfied by $u(t,x)$. Sato discovered the intrinsic
meaning of $u(t,x)$ and found that solutions to the KdV equation can be
described by the tau-functions. Although the theorems in this paper can be
proved without this representation, in view of its historical significance we
show it in the present framework.

To define the tau-function $\tau_{\boldsymbol{a}}\left(  g\right)  $ we have
to assume that the operator%
\[
g^{-1}T\left(  g\boldsymbol{a}\right)  T\left(  \boldsymbol{a}\right)
^{-1}-I:H_{N}\left(  D_{+}\right)  \rightarrow H_{N}\left(  D_{+}\right)
\]
is of trace class. We assume in this section the condition (\ref{39}) of Lemma
\ref{l3} for $\lambda\boldsymbol{a}\left(  \lambda\right)  $, $N=1$, that is%
\begin{equation}%
{\displaystyle\int_{C^{2}}}
\left\vert \dfrac{z^{2}\widetilde{a}_{j}\left(  z\right)  -\lambda
^{2}\widetilde{a}_{j}\left(  \lambda\right)  }{z-\lambda}\right\vert
^{2}\left\vert dz\right\vert \left\vert d\lambda\right\vert <\infty\text{ with
}\widetilde{a}_{j}\left(  z\right)  =a_{j}(z)-f_{j}(z)\text{,} \label{14}%
\end{equation}
which shows the operator $S_{\lambda\boldsymbol{a}}$ is of HS from
$H_{1}\left(  D_{+}\right)  $ to $H\left(  D_{-}\right)  $. This condition can
be verified in the same manner as Lemma \ref{l10} for $\boldsymbol{a}%
(z)=\left(  1,m\left(  z\right)  /z\right)  $ if $m$ satisfies (M.1), (M.2)
for sufficiently large $L$.\medskip

Let%
\[
e_{x}(z)=e^{xz}\in\Gamma_{1}\text{, \ \ }e_{t,x}(z)=e^{xz+tz^{3}}\in\Gamma
_{3}\text{.}%
\]
In Proposition \ref{p1} for $\boldsymbol{a}\in\boldsymbol{A}_{3}\left(
C\right)  $ satisfying $e_{x}\boldsymbol{a}\in\boldsymbol{A}_{3}^{inv}\left(
C\right)  $ for any $x\in\mathbb{R}$ we have introduced the potential $q$
associated with $\boldsymbol{a}$ of Schr\"{o}dinger operator by
$q(x)=-2\partial_{x}\kappa_{1}\left(  e_{x}\boldsymbol{a}\right)  $. The
$\kappa_{1}\left(  \boldsymbol{a}\right)  $ is described by the characteristic
functions as%
\[
\kappa_{1}\left(  \boldsymbol{a}\right)  =\lim_{\zeta\rightarrow\infty}%
\zeta\varphi_{\boldsymbol{a}}\left(  \zeta\right)  \text{.}%
\]

\begin{proposition}
\label{p3}Assume $\boldsymbol{a}\in\boldsymbol{A}_{L,+}^{inv}\left(  C\right)
$ and (\ref{14}). Then an identity%
\[
\kappa_{1}\left(  e_{x}\boldsymbol{a}\right)  =\partial_{x}\log\tau
_{\boldsymbol{a}}\left(  e_{x}\right)
\]
holds, which yields%
\[
q(x)=-2\partial_{x}^{2}\log\tau_{\boldsymbol{a}}\left(  e_{x}\right)  \text{,}%
\]
if $L\geq2$. Generally for $g\in\Gamma_{n}$ and $\boldsymbol{a}\in
\boldsymbol{A}_{L,+}^{inv}\left(  C\right)  $ with $L\geq\max\left\{
n,2\right\}  $
\[
\left(  K\left(  g\right)  q\right)  \left(  x\right)  =-2\partial_{x}^{2}%
\log\tau_{\boldsymbol{a}}\left(  ge_{x}\right)
\]
holds. Especially the solution $q(t,x)$ to the KdV equation starting from
$q(x)$ is given by%
\[
q(t,x)=-2\partial_{x}^{2}\log\tau_{\boldsymbol{a}}\left(  e_{t,x}\right)
\text{,}%
\]
if $\boldsymbol{a}\in\boldsymbol{A}_{L,+}^{inv}\left(  C\right)  $ for
$L\geq4$. The condition $L\geq4$ is necessary for the differentiability of
$q(t,x)$ in $t$ (see Proposition \ref{p2}).
\end{proposition}

\begin{proof}
The definition of $\varphi_{\boldsymbol{a}}$ implies%
\begin{align}
\kappa_{1}\left(  \boldsymbol{a}\right)   &  =\lim_{\zeta\rightarrow\infty
}\zeta\varphi_{\boldsymbol{a}}\left(  \zeta\right)  =\lim_{\zeta
\rightarrow\infty}\zeta\frac{1}{2\pi i}\int_{C}\frac{\varphi_{\boldsymbol{a}%
}\left(  \lambda\right)  }{\zeta-\lambda}d\lambda\nonumber\\
&  =\lim_{\zeta\rightarrow\infty}\zeta\frac{1}{2\pi i}\int_{C}\frac
{\widetilde{\boldsymbol{a}}\left(  \lambda\right)  \left(  T\left(
\boldsymbol{a}\right)  ^{-1}1\right)  \left(  \lambda\right)  }{\zeta-\lambda
}d\lambda\nonumber\\
&  =\frac{1}{2\pi i}\int_{C}\widetilde{\boldsymbol{a}}\left(  \lambda\right)
\left(  T\left(  \boldsymbol{a}\right)  ^{-1}1\right)  \left(  \lambda\right)
d\lambda\text{.} \label{113}%
\end{align}
On the other hand the formal identity%
\[
\tau_{\boldsymbol{a}}\left(  e_{\epsilon}\right)  =\det\left(  e_{\epsilon
}^{-1}T\left(  e_{\epsilon}\boldsymbol{a}\right)  T\left(  \boldsymbol{a}%
\right)  ^{-1}\right)  =\det\left(  I+e_{\epsilon}^{-1}H_{e_{\epsilon}%
}S_{\boldsymbol{a}}T\left(  \boldsymbol{a}\right)  ^{-1}\right)
\]
is justified by reproving below in a way different from that of Lemma \ref{l7}
that $e_{\epsilon}^{-1}H_{e_{\epsilon}}S_{\boldsymbol{a}}T\left(
\boldsymbol{a}\right)  ^{-1}$ defines a trace class operator on $H_{1}\left(
D_{+}\right)  $. For $v\in H_{1}\left(  D_{+}\right)  $ it holds that%
\[
S_{\boldsymbol{a}}v\left(  z\right)  =\frac{1}{2\pi i}\int_{C}\frac
{\widetilde{\boldsymbol{a}}\left(  \lambda\right)  v\left(  \lambda\right)
}{z-\lambda}d\lambda=\ell_{1}\left(  \boldsymbol{a}\text{,}v\right)
z^{-1}+z^{-1}\left(  S_{\boldsymbol{a}}^{\left(  1\right)  }v\right)  \left(
z\right)
\]
with%
\[
\left\{
\begin{array}
[c]{l}%
\ell_{1}\left(  \boldsymbol{a}\text{,}v\right)  =\dfrac{1}{2\pi i}%
{\displaystyle\int_{C}}
\widetilde{\boldsymbol{a}}\left(  \lambda\right)  v\left(  \lambda\right)
d\lambda\\
S_{\boldsymbol{a}}^{\left(  1\right)  }v\left(  z\right)  =\dfrac{1}{2\pi i}%
{\displaystyle\int_{C}}
\dfrac{\lambda\widetilde{\boldsymbol{a}}\left(  \lambda\right)  v\left(
\lambda\right)  }{z-\lambda}d\lambda\\
\text{ \ \ \ \ \ \ \ \ \ }=\dfrac{1}{2\pi i}%
{\displaystyle\int_{C}}
\dfrac{\left(  \lambda\widetilde{\boldsymbol{a}}\left(  \lambda\right)
-z\widetilde{\boldsymbol{a}}\left(  z\right)  \right)  v\left(  \lambda
\right)  }{z-\lambda}d\lambda\in H\left(  D_{-}\right)
\end{array}
\right.  \text{,}%
\]
since $\widetilde{\boldsymbol{a}}\left(  \lambda\right)  =O\left(
\lambda^{-L}\right)  $ for $L\geq2$. Hence for $z\in D_{+}$%
\begin{align}
e_{\epsilon}^{-1}H_{e_{\epsilon}}S_{\boldsymbol{a}}v\left(  z\right)   &
=\frac{\ell_{1}\left(  \boldsymbol{a}\text{,}v\right)  }{2\pi i}\int_{C}%
\frac{e^{\epsilon\left(  \lambda-z\right)  }}{\left(  \lambda-z\right)
\lambda}d\lambda+e_{\epsilon}^{-1}H_{e_{\epsilon}}z^{-1}S_{\boldsymbol{a}%
}^{\left(  1\right)  }v\left(  z\right) \nonumber\\
&  =\ell_{1}\left(  \boldsymbol{a}\text{,}v\right)  \frac{1-e^{-\epsilon z}%
}{z}+e_{\epsilon}^{-1}H_{e_{\epsilon}}z^{-1}S_{\boldsymbol{a}}^{\left(
1\right)  }v\left(  z\right)  \label{116}%
\end{align}
holds. Since $S_{\boldsymbol{a}}^{\left(  1\right)  }$ defines an HS operator
from $H_{1}\left(  D_{+}\right)  $ to $H\left(  D_{-}\right)  $ under
(\ref{14}), $e_{\epsilon}^{-1}H_{e_{\epsilon}}z^{-1}S_{\boldsymbol{a}%
}^{\left(  1\right)  }$ turns to a trace class operator on $H_{1}\left(
D_{+}\right)  $, which makes it possible to define $\tau_{\boldsymbol{a}%
}\left(  e_{\epsilon}\right)  $ rigorously. Moreover in this case for $w\in
H\left(  D_{-}\right)  $%
\[
e_{\epsilon}^{-1}H_{e_{\epsilon}}z^{-1}w\left(  z\right)  =\frac{1}{2\pi
i}\int_{C}\left(  \frac{e^{\epsilon\left(  \lambda-z\right)  }-1}{\lambda
-z}-\epsilon\right)  \lambda^{-1}w\left(  \lambda\right)  d\lambda
\]
holds due to%
\[
\frac{1}{2\pi i}\int_{C}\frac{1}{\lambda-z}\lambda^{-1}w\left(  \lambda
\right)  d\lambda=\frac{1}{2\pi i}\int_{C}\lambda^{-1}w\left(  \lambda\right)
d\lambda=0\text{,}%
\]
Hence the square of the HS-norm of $\epsilon^{-1}e_{\epsilon}^{-1}%
H_{e_{\epsilon}}z^{-1}$ is%
\[
\delta_{\epsilon}\equiv\left(  2\pi\right)  ^{-2}\int_{C^{2}}\left\vert
\frac{e^{\epsilon\left(  \lambda-z\right)  }-1}{\epsilon\left(  \lambda
-z\right)  }-1\right\vert ^{2}\left\vert \lambda\right\vert ^{-2}\left\vert
z\right\vert ^{-2}\left\vert d\lambda\right\vert \left\vert dz\right\vert
\text{.}%
\]
Since there exists a constant $c$ such that%
\[
\left\vert \frac{e^{\epsilon\left(  \lambda-z\right)  }-1}{\epsilon\left(
\lambda-z\right)  }-1\right\vert =\left\vert \int_{0}^{1}\left(
e^{t\epsilon\left(  \lambda-z\right)  }-1\right)  dt\right\vert \leq c
\]
holds for $\epsilon\in\mathbb{R}$, $\lambda$, $z\in C$, the dominated
convergence theorem shows $\lim_{\epsilon\rightarrow0}\delta_{\epsilon}=0$.
Consequently we have%
\begin{equation}
\left\Vert e_{\epsilon}^{-1}H_{e_{\epsilon}}z^{-1}S_{\boldsymbol{a}}^{\left(
1\right)  }T\left(  \boldsymbol{a}\right)  ^{-1}\right\Vert _{trace}\leq
c_{1}\epsilon\sqrt{\delta_{\epsilon}}=o\left(  \epsilon\right)  \text{.}
\label{117}%
\end{equation}
The first term of (\ref{116}) generates a rank $1$ operator and%
\[
\lim_{\epsilon\rightarrow0}\frac{1}{\epsilon}\ell_{1}\left(  \boldsymbol{a}%
\text{,}v\right)  \frac{1-e^{-\epsilon z}}{z}=\ell_{1}\left(  \boldsymbol{a}%
\text{,}v\right)  =\dfrac{1}{2\pi i}%
{\displaystyle\int_{C}}
\widetilde{\boldsymbol{a}}\left(  \lambda\right)  v\left(  \lambda\right)
d\lambda\text{.}%
\]
Noting an identity (see (\cite{si}))%
\[
\det\left(  I+A\right)  =\exp\left(  \mathrm{tr}\log\left(  I+A\right)
\right)  =\exp\left(  \mathrm{tr}A+O\left(  \left\Vert A\right\Vert
_{trace}^{2}\right)  \right)
\]
if $\left\Vert A\right\Vert _{trace}<1$, we have%
\begin{align*}
\tau_{\boldsymbol{a}}\left(  e_{\epsilon}\right)   &  =\exp\left(
\mathrm{tr}\left(  e_{\epsilon}^{-1}H_{e_{\epsilon}}S_{\boldsymbol{a}}T\left(
\boldsymbol{a}\right)  ^{-1}\right)  +o\left(  \epsilon\right)  \right) \\
&  =\exp\left(  \epsilon\dfrac{1}{2\pi i}%
{\displaystyle\int_{C}}
\widetilde{\boldsymbol{a}}\left(  \lambda\right)  \left(  T\left(
\boldsymbol{a}\right)  ^{-1}1\right)  \left(  \lambda\right)  d\lambda
+o\left(  \epsilon\right)  \right) \\
&  =1+\epsilon\dfrac{1}{2\pi i}%
{\displaystyle\int_{C}}
\widetilde{\boldsymbol{a}}\left(  \lambda\right)  \left(  T\left(
\boldsymbol{a}\right)  ^{-1}1\right)  \left(  \lambda\right)  d\lambda
+o\left(  \epsilon\right)
\end{align*}
for sufficiently small $\epsilon$ due to (\ref{117}). Then (\ref{113}) implies%
\[
\lim_{\epsilon\rightarrow0}\frac{\tau_{\boldsymbol{a}}\left(  e_{\epsilon
}\right)  -1}{\epsilon}=\kappa_{1}\left(  \boldsymbol{a}\right)  \text{,}%
\]
and the cocycle property shows%
\[
\lim_{\epsilon\rightarrow0}\frac{\tau_{\boldsymbol{a}}\left(  e_{x+\epsilon
}\right)  -\tau_{\boldsymbol{a}}\left(  e_{x}\right)  }{\epsilon}%
=\lim_{\epsilon\rightarrow0}\frac{\tau_{e_{x}\boldsymbol{a}}\left(
e_{\epsilon}\right)  -1}{\epsilon}\tau_{\boldsymbol{a}}\left(  e_{x}\right)
=\kappa_{1}\left(  e_{x}\boldsymbol{a}\right)  \tau_{\boldsymbol{a}}\left(
e_{x}\right)  \text{,}%
\]
hence $\partial_{x}\log\tau_{\boldsymbol{a}}\left(  e_{x}\right)  =\kappa
_{1}\left(  e_{x}\boldsymbol{a}\right)  $ holds. The rest of the proof is automatic.
\end{proof}

\section{Sufficient conditions for $q\in\mathcal{Q}_{L}\left(  C\right)  $}

A sufficient condition for $q\in\mathcal{Q}_{L}\left(  C\right)  $ was given
in (\ref{77}) in terms of its Weyl functions $m_{\pm}$ (see Theorem \ref{t1}).
In this section we provide concrete examples of this class including the
well-known cases. Throughout this section we treat $g=e^{h}$ with real odd
polynomial $h$ of degree $n$, hence the curve $C$ is taken so that $g$ is
bounded on $D_{+}$, namely%
\[
C=\left\{
\begin{array}
[c]{c}%
\pm\omega\left(  y\right)  +iy\text{; \ }y\in\mathbb{R}\text{, }\omega\left(
y\right)  >0\text{, }\omega\left(  y\right)  =\omega\left(  -y\right)
\text{,}\\
\omega\text{ is smooth and satisfies }\omega\left(  y\right)  =O\left(
y^{-\left(  n-1\right)  }\right)  \text{ as }y\rightarrow\infty
\end{array}
\right\}  \text{.}%
\]

\subsection{Decaying potentials}

If a potential $q$ satisfies $q\in L^{1}\left(  \mathbb{R}_{+}\right)  $, it
is known that for $0\neq k\in\overline{\mathbb{C}}_{+}\equiv\left\{
z\in\mathbb{C}\text{; }\operatorname{Im}z\geq0\right\}  $ there exists the
Jost solution $f_{+}\left(  x,k\right)  $ of%
\[
-\partial_{x}^{2}f_{+}\left(  x,k\right)  +q(x)f_{+}\left(  x,k\right)
=k^{2}f_{+}\left(  x,k\right)
\]
such that%
\[
\left\{
\begin{array}
[c]{l}%
f_{+}\left(  x,k\right)  =e^{ixk}+o\left(  1\right) \\
f_{+}^{\prime}\left(  x,k\right)  =ike^{ixk}+o\left(  1\right)
\end{array}
\right.  \text{ \ as }x\rightarrow\infty\text{,}%
\]
where $^{\prime}$ denotes the derivative with respect to $x$. Therefore%
\[
m_{+}\left(  z\right)  =\dfrac{f_{+}^{\prime}\left(  0,\sqrt{z}\right)
}{f_{+}\left(  0,\sqrt{z}\right)  }%
\]
and one can see that $m_{+}\left(  z\right)  $ is extendable to $\overline
{\mathbb{C}}_{+}\backslash\left\{  0\right\}  $ as a continuous function.
$f_{+}\left(  x,k\right)  $ is obtained as a unique solution to an integral
equation%
\[
e^{-ixk}f_{+}\left(  x,k\right)  =1+\int_{x}^{\infty}\dfrac{e^{2ik\left(
s-x\right)  }-1}{2ik}q\left(  s\right)  e^{-iks}f_{+}\left(  s,k\right)
ds\text{.}%
\]
Rybkin (\cite{r1}) showed%
\begin{equation}
e^{-ixk}f_{+}\left(  x,k\right)  =1+\sum_{j=1}^{N+1}f_{j}\left(  x\right)
\left(  2ik\right)  ^{-j}+o\left(  k^{-N-1}\right)  \label{41}%
\end{equation}
for $q$ such that $q^{\left(  j\right)  }\in L^{1}\left(  \mathbb{R}\right)  $
for $j=0$, $1$,$\cdots$, $L$. The small $o$ is uniform with respect to
$x\geq0$. The coefficients $\left\{  f_{j}\left(  x\right)  \right\}  $ are
determined inductively by%
\[
\left\{
\begin{array}
[c]{l}%
f_{1}(x)=-Q(x)\equiv-\int_{x}^{\infty}q\left(  s\right)  ds\\
f_{j+1}(x)=-f_{j}^{\prime}\left(  x\right)  -\int_{x}^{\infty}q\left(
s\right)  f_{j}\left(  s\right)  ds\text{, \ }\left(  j\geq1\right)
\end{array}
\right.  \text{.}%
\]
Therefore one can show $f_{j+1}$ is $L-j$ times differentiable and
$f_{j+1}^{\left(  L-j+1\right)  }\in L^{1}\left(  \mathbb{R}_{+}\right)  $.
Since%
\[
\left(  e^{-ixk}f_{+}\left(  x,k\right)  \right)  ^{\prime}=\int_{x}^{\infty
}e^{2ik\left(  s-x\right)  }q\left(  s\right)  e^{-iks}f_{+}\left(
s,k\right)  ds\text{,}%
\]
substituting (\ref{41}) we have asymptotic behavior%
\[
\left(  e^{-ixk}f_{+}\left(  x,k\right)  \right)  ^{\prime}=\sum_{j=1}%
^{L+1}g_{j}\left(  x\right)  \left(  2ik\right)  ^{-j}+o\left(  k^{-N-1}%
\right)  \text{,}%
\]
which leads us to%
\[
m_{+}\left(  z\right)  =-\sqrt{-z}+\sum_{j=1}^{L+1}c_{j}\left(  -z\right)
^{-j/2}+o\left(  z^{-\left(  L+1\right)  /2}\right)  \text{ \ if }%
z\rightarrow\infty\text{ on }\overline{\mathbb{C}}_{+}\text{.}%
\]
An analogous asymptotic behavior for $m_{-}(z)$ is possible if $q^{\left(
j\right)  }\in L^{1}\left(  \mathbb{R}\right)  $ for $j=0$, $1$,$\cdots$, $L$
by replacing $c_{j}$ by $\left(  -1\right)  ^{j+1}c_{j}$, which shows

\begin{proposition}
\label{p10}If $q^{\left(  j\right)  }\in L^{1}\left(  \mathbb{R}\right)  $ for
$j=0$, $1$,$\cdots$, $L$, then (M.2) is satisfied with $L+2$ for any curve
$C$, and we have $q\in\mathcal{Q}_{L+1}\left(  C\right)  $. Therefore one can
define the KdV flow $K(g)q$ for $g\in\Gamma_{n}$ if $L\geq1$. For the KdV
equation $L\geq3$ is required to guarantee the differentiability.
\end{proposition}

One cannot apply this proposition to the interesting case $q(x)=\varphi\left(
x\right)  /x$ with smooth periodic function $\varphi$ satisfying
$\varphi\left(  0\right)  =0$, however there is a possibility of estimating
directly $m_{\pm}$ in this case by a sort of perturbation.

\subsection{Reflection coefficients}

To obtain another class of $q$ satisfying (M.2) we prepare the necessary
terminologies from the spectral theory of Schr\"{o}dinger operators. Since
$m_{\pm}\left(  z\right)  $ take values in $\mathbb{C}_{+}$ for $z\in
\mathbb{C}_{+}$, we start from

\begin{lemma}
\label{l23}For any complex numbers $m_{\pm}\in\mathbb{C}_{+}$ set%
\[
m_{1}=-\dfrac{1}{m_{+}+m_{-}}\text{, \ }m_{2}=\dfrac{m_{+}m_{-}}{m_{+}+m_{-}%
}\text{, \ }R=\dfrac{\overline{m_{+}}+m_{-}}{m_{+}+m_{-}}\text{.}%
\]
Then, $m_{1}$, $m_{2}\in\mathbb{C}_{+}$, $\left\vert R\right\vert \leq1$ hold,
and $\xi_{j}=\left(  \arg m_{j}\right)  /\pi\in\left[  0,1\right]  $ ($j=1$,
$2$) satisfy%
\begin{equation}
\left\vert \xi_{1}-\dfrac{1}{2}\right\vert \text{, \ }\ \left\vert \xi
_{2}-\dfrac{1}{2}\right\vert \leq\dfrac{1}{2}\left\vert R\right\vert \text{.}
\label{88}%
\end{equation}

\end{lemma}

For\ $z\in\mathbb{C}_{+}$ set%
\begin{equation}
\left\{
\begin{array}
[c]{l}%
m_{1}(z)=-\dfrac{1}{m_{+}(z)+m_{-}(z)}\text{, \ }m_{2}(z)=\dfrac{m_{+}%
(z)m_{-}(z)}{m_{+}(z)+m_{-}(z)}\\
\xi_{j}\left(  z\right)  =\dfrac{1}{\pi}\arg m_{j}\left(  z\right)  \left(
=\dfrac{1}{\pi}\operatorname{Im}\log m_{j}\left(  z\right)  \right)  \text{,
(}j=1\text{, }2\text{)}%
\end{array}
\right.  \text{.} \label{92}%
\end{equation}
Then, $\left\{  m_{j},j=1,2\right\}  $ are Herglotz functions. $m\left(
z\right)  $ defined by%
\[
m\left(  z\right)  =\left\{
\begin{array}
[c]{c}%
-m_{+}\left(  -z^{2}\right)  \text{ \ \ if \ }\operatorname{Re}z>0\\
m_{-}\left(  -z^{2}\right)  \text{ \ \ \ if \ }\operatorname{Re}z<0
\end{array}
\right.
\]
satisfies the asymptotic behavior (\ref{77}) if and only if $m_{\pm}$ satisfy%
\begin{equation}
\left\{
\begin{array}
[c]{l}%
m_{+}\left(  -z^{2}\right)  =-z\left(  1-\sum_{k=2}^{L-1}c_{k}z^{-k}+O\left(
z^{-L}\right)  \right)  \smallskip\\
m_{-}\left(  -z^{2}\right)  =-z\left(  1-\sum_{k=2}^{L-1}c_{k}\left(
-z\right)  ^{-k}+O\left(  z^{-L}\right)  \right)
\end{array}
\right.  \text{.} \label{15}%
\end{equation}
as $z\rightarrow\infty$ on $D_{-}\cap\left\{  \operatorname{Re}z>0\right\}  $.
It is known that if $q\in C^{L-2}\left(  [0,\delta)\right)  $ for some
$\delta>0$, then defining inductively the functions $c_{j}\left(  x\right)  $
by%
\[
\left\{
\begin{array}
[c]{l}%
c_{1}\left(  x\right)  =0,\ \ \ \ \ c_{2}\left(  x\right)  =q(x)/2,\smallskip
\\
c_{j}\left(  x\right)  =(c_{j-1}^{^{\prime}}\left(  x\right)  -\sum_{\ell
=1}^{j-1}c_{\ell}\left(  x\right)  c_{j-\ell}\left(  x\right)  )/2,\ \ \ j\geq
3
\end{array}
\right.
\]
one has $c_{k}=c_{k}(0)$. The coefficients for $m_{-}\left(  -z^{2}\right)  $
can be obtained by considering $q\left(  -x\right)  $ in place of $q(x)$.
Then, if $q\in C^{L-2}(-\delta,\delta)$, (\ref{15}) implies%
\begin{equation}
\left\{
\begin{array}
[c]{l}%
m_{1}(-z^{2})=\dfrac{1}{2}z^{-1}\left(  1+\sum_{k=1}^{M}a_{k}z^{-2k}+O\left(
z^{-L}\right)  \right)  \smallskip\\
m_{2}(-z^{2})=-\dfrac{1}{2}z\left(  1+\sum_{k=1}^{M}b_{k}z^{-2k}+O\left(
z^{-L}\right)  \right)
\end{array}
\right.  \label{89}%
\end{equation}
on $D_{-}\cap\left\{  \operatorname{Re}z>0\right\}  $ with some $a_{k}$,
$b_{k}\in\mathbb{R}$, where $M=\left[  \left(  L-1\right)  /2\right]  $
($\left[  x\right]  $ denotes the integer part of $x$). $m_{\pm}$ can be
recovered from $m_{1}$, $m_{2}$ by%
\[
m_{\pm}=-\dfrac{1}{2m_{1}}\left(  1\pm\sqrt{1+4m_{1}m_{2}}\right)  \text{.}%
\]
Observe%
\begin{align*}
1+4m_{1}(-z^{2})m_{2}(-z^{2})  &  =\left(  \dfrac{m_{+}(-z^{2})-m_{-}(-z^{2}%
)}{m_{+}(-z^{2})+m_{-}(-z^{2})}\right)  ^{2}\\
&  =\left(  \dfrac{f_{o}(z)+O\left(  z^{-L}\right)  }{1-f_{e}(z)+O\left(
z^{-L}\right)  }\right)  ^{2}%
\end{align*}
with $f(z)=\sum_{k=2}^{L-1}c_{k}z^{-k}$. Let $N$ be the least number such that
$\left(  1-\left(  -1\right)  ^{k}\right)  c_{k}\neq0$. $N$ should be odd and
$N\geq3$, since $f\left(  z\right)  =c_{2}z^{-2}+O\left(  z^{-3}\right)  $.
Then one sees%
\begin{equation}
1+4m_{1}(-z^{2})m_{2}(-z^{2})=-\sum_{k=N}^{M}d_{k}z^{-2k}+O\left(
z^{-L-N}\right)  \label{118}%
\end{equation}
with $d_{k}\in\mathbb{R}$, $d_{N}\neq0$. However one cannot expect to have
(\ref{118}) from (\ref{89}).

We would like to have a partial converse:

\begin{lemma}
\label{l24}For Herglotz functions $m_{\pm}$ suppose $m_{1}$, $m_{2}$ are given
by (\ref{92}) and satisfy (\ref{89}) with $M^{\prime}$, $L^{\prime}$ such that
$2M^{\prime}<L^{\prime}$. Then $m$ satisfies the conditions (M.1), (M.2) in
(\ref{77}) with $L=1+\left[  L^{\prime}/2\right]  $.
\end{lemma}

\begin{proof}
We have only to verify (M.2). Since%
\[
m_{+}+m_{-}=-m_{1}^{-1}\text{, \ \ \ }m_{+}m_{-}=-m_{1}^{-1}m_{2}\text{,}%
\]
we have%
\[
m_{\pm}=\dfrac{1}{2}\left(  -m_{1}^{-1}\mp\sqrt{m_{1}^{-2}+4m_{1}^{-1}m_{2}%
}\right)  =-\dfrac{1}{2m_{1}}\left(  1\pm\sqrt{1+4m_{1}m_{2}}\right)  \text{.}%
\]
Then%
\begin{align*}
\sqrt{\frac{1+4m_{1}\left(  -z^{2}\right)  m_{2}\left(  -z^{2}\right)  }%
{d_{N}}}  &  =z^{-N}\sqrt{1+\sum_{k=N+1}^{M^{\prime}}\dfrac{d_{k}}{d_{N}%
}z^{-2\left(  k-N\right)  }+O\left(  z^{-L^{\prime}+2N}\right)  }\\
&  =z^{-N}\left(  1+\sum_{k=N+1}^{M^{\prime}}d_{k}^{\prime}z^{-2\left(
k-N\right)  }+O\left(  z^{-L^{\prime}+2N}\right)  \right) \\
&  =z^{-N}+\sum_{k=N+1}^{M^{\prime}}d_{k}^{\prime}z^{-2k+N}+O\left(
z^{-L^{\prime}+N}\right)
\end{align*}
holds with other constants $d_{k}^{\prime}$. Since $L^{\prime}-N\geq
L^{\prime}-M^{\prime}>L^{\prime}/2$, we have the lemma.\medskip
\end{proof}

The asymptotic behavior of $m_{j}$ is translated to that of $\xi_{j}$ as
follows. If $\operatorname{Im}z>0$, then $\operatorname{Im}m_{j}(z)>0$, hence
$0\leq\xi_{j}\left(  z\right)  \leq1$ holds. $\log m_{j}$ are of Herglotz as
well, since $\operatorname{Im}\log m_{j}(z)=\pi\xi_{j}\left(  z\right)  \geq
0$. On the other hand (ii) of (\ref{84}) implies that $m_{1}(z)$, $m_{2}(z)$
take real values on $\left(  -\infty,\lambda_{0}\right)  $ and $m_{1}%
(\lambda)>0$ there. Let $\lambda_{1}\leq\lambda_{0}$ be a unique zero of
$m_{2}(z)$ if it has, and set $\lambda_{1}=\lambda_{0}$ if it has not. Assume
in the sequel%
\begin{equation}
\int_{0}^{\infty}\left\vert \xi_{j}\left(  \lambda\right)  -\dfrac{1}%
{2}\right\vert d\lambda<\infty\text{ \ for }j=1\text{, }2\text{.} \label{96}%
\end{equation}
Then $m_{j}$ are represented as%
\begin{equation}
\left\{
\begin{array}
[c]{l}%
m_{1}(z)=\dfrac{1}{2\sqrt{-z}}\exp\left(
{\displaystyle\int_{\lambda_{0}}^{\infty}}
\dfrac{\xi_{1}\left(  \lambda\right)  -I_{\lambda>0}/2}{\lambda-z}%
d\lambda\right)  \smallskip\\
m_{2}(z)=-\dfrac{\sqrt{-z}}{2}\dfrac{\lambda_{1}-z}{-z}\exp\left(
{\displaystyle\int_{\lambda_{0}}^{\infty}}
\dfrac{\xi_{2}\left(  \lambda\right)  -I_{\lambda>0}/2}{\lambda-z}%
d\lambda\right)
\end{array}
\right.  \text{.} \label{94}%
\end{equation}
The function $\arg m_{+}\left(  \lambda+i0\right)  $ is intensively
investigated by Gesztesy-Simon \cite{g-s} in connection with inverse spectral
problems, and they call $\arg m_{+}\left(  \lambda+i0\right)  /\pi$ as
\textbf{xi-function}. To have (\ref{89}) for some $L\geq1$ it is sufficient
that%
\begin{equation}
\int_{0}^{\infty}\lambda^{M}\left\vert \xi_{j}\left(  \lambda\right)
-\dfrac{1}{2}\right\vert d\lambda<\infty\text{, \ \ (}j=1\text{, }2\text{)}
\label{91}%
\end{equation}
hold for an $M\geq1$. This is due to the expansion%
\begin{equation}
-\int_{\lambda_{0}}^{\infty}\dfrac{f\left(  \lambda\right)  }{\lambda
-z}d\lambda=\sum_{k=0}^{M-1}z^{-k-1}\int_{\lambda_{0}}^{\infty}\lambda
^{k}f\left(  \lambda\right)  d\lambda+z^{-M}\int_{\lambda_{0}}^{\infty}%
\dfrac{\lambda^{M}f\left(  \lambda\right)  }{\lambda-z}d\lambda\text{,}
\label{100}%
\end{equation}
and the estimate%
\begin{equation}
\left\vert \int_{\lambda_{0}}^{\infty}\dfrac{\lambda^{M}f\left(
\lambda\right)  }{\lambda-z}d\lambda\right\vert \leq\int_{\lambda_{0}}%
^{\infty}\dfrac{\left\vert \lambda^{M}f\left(  \lambda\right)  \right\vert
}{\left\vert \lambda-z\right\vert }d\lambda\leq\dfrac{1}{\left\vert
\operatorname{Im}z\right\vert }\int_{\lambda_{0}}^{\infty}\left\vert
\lambda^{M}f\left(  \lambda\right)  \right\vert d\lambda\text{.} \label{107}%
\end{equation}
We control $\left\vert \xi_{j}\left(  \lambda\right)  -1/2\right\vert $ by
another quantity. The \textbf{reflection coefficient} $R(z)$ is defined by%

\[
R(z)=\dfrac{\overline{m_{+}(z)}+m_{-}(z)}{m_{+}(z)+m_{-}(z)}\text{.}%
\]
This quantity was considered by Gesztesy-Simon, Rybkin and others as a
generalization of the conventional reflection coefficient defined for decaying potentials.

\subsection{Proof of Theorems \ref{t3}, \ref{t4}}

Assume%
\begin{equation}
\int_{0}^{\infty}\lambda^{M}\left\vert R\left(  \lambda\right)  \right\vert
d\lambda<\infty. \label{105}%
\end{equation}
Then Lemma \ref{l23} implies that (\ref{91}) holds, and from (\ref{100}),
(\ref{107}) one has%
\[
-\int_{\lambda_{0}}^{\infty}\dfrac{\xi_{j}\left(  \lambda\right)  -\dfrac
{1}{2}}{\lambda-z}d\lambda=\sum_{k=0}^{M-1}z^{-k-1}\int_{\lambda_{0}}^{\infty
}\lambda^{k}\left(  \xi_{j}\left(  \lambda\right)  -\dfrac{1}{2}\right)
d\lambda+O\left(  z^{-M+n/2-1}\right)
\]
if $z\in D_{-}$. Applying Lemma \ref{l24} with $L^{\prime}=2M-n+2$ we see that
this $m$ satisfies (M.2) with%
\[
L=1+\left[  L^{\prime}/2\right]  =1+\left[  M-n/2+1\right]  =M+1-\left(
n-1\right)  /2\text{,}%
\]
which yields Theorem \ref{t3}.

The condition (\ref{105}) implies that the ac spectrum is large, which
restricts the possible potential class strongly. We try to relax the condition
(\ref{105}) by replacing it with a condition on a curve surrounding
$[\lambda_{0},\infty)$.

We prepare two curves%
\[
\left\{
\begin{array}
[c]{l}%
C=\left\{  \pm\omega\left(  y\right)  +iy\text{; \ }y\in\mathbb{R}\right\}
\text{ with }\omega\left(  y\right)  =cy^{-\left(  n-1\right)  }\text{ for
}\left\vert y\right\vert \geq1\\
C_{1}=\left\{  \pm\omega_{1}\left(  y\right)  +iy\text{; \ }y\in
\mathbb{R}\right\}  \text{ with }\omega_{1}\left(  y\right)  =cy^{-\left(
n_{1}-1\right)  }\text{ for }\left\vert y\right\vert \geq1
\end{array}
\right.
\]
for $n_{1}>n$, and%
\[
\left\{
\begin{array}
[c]{l}%
\widehat{C}=\left\{  -z^{2}\text{; }z\in C\text{, }\operatorname{Re}%
z>0\right\}  =\left\{  x\pm i\widehat{\omega}\left(  x\right)  \text{; }%
x\in\mathbb{R}\text{, \ }x\geq\lambda_{0}\right\} \\
\widehat{C}_{1}=\left\{  -z^{2}\text{; }z\in C_{1}\text{, }\operatorname{Re}%
z>0\right\}  =\left\{  x\pm i\widehat{\omega}_{1}\left(  x\right)  \text{;
}x\in\mathbb{R}\text{, \ }x\geq\lambda_{0}\right\}
\end{array}
\right.
\]
with $\widehat{\omega}\left(  x\right)  =\widehat{c}x^{1-n/2}$, $\widehat
{\omega}_{1}\left(  x\right)  =\widehat{c}_{1}x^{1-n_{1}/2}$ for $x\geq1$.
Then one can assume $D_{-}\subset D_{1,-}$.

\begin{lemma}
\label{l26}Let $M\in\mathbb{Z}$ be $M\geq n/2$. Assume $q\in C^{2M-n}\left(
-\delta,\delta\right)  $ and%
\begin{equation}
\int_{\widehat{C}_{1}}\left\vert z^{M}R(z)\right\vert \left\vert dz\right\vert
<\infty\label{171}%
\end{equation}
holds. Then, the $m$ satisfies (M.2) with $L\leq M+1-\left(  n-1\right)  /2$
on the curve $C$.
\end{lemma}

\begin{proof}
Let $\phi$ be $\phi\left(  z\right)  =-\phi_{k}\left(  -z\right)  $ in Lemma
\ref{l51} with $k=\left(  n_{1}-1\right)  /2$. Then $\phi$ maps $\mathbb{C}%
\backslash\lbrack0,\infty)$ onto $\widehat{D}_{1,-}=\left\{  z\text{;
}\left\vert \operatorname{Im}z\right\vert >\widehat{\omega}_{1}\left(
\operatorname{Re}z\right)  \text{, }\operatorname{Re}z>\lambda_{0}\right\}  $
conformally. Without loss of generality we can assume $-a_{k}^{2}<\lambda_{0}%
$. (\ref{171}) implies%
\[
\int_{0}^{\infty}\left\vert \phi\left(  \lambda\right)  ^{M}R\left(
\phi\left(  \lambda\right)  \right)  )\right\vert \left\vert d\phi\left(
\lambda\right)  \right\vert <\infty\text{.}%
\]
which is equivalent to%
\[
\int_{0}^{\infty}\left\vert \lambda^{M}R(\phi\left(  \lambda\right)
)\right\vert d\lambda<\infty
\]
due to $\phi^{\prime}\left(  \lambda\right)  =1+o\left(  1\right)  $. Hence
Lemma \ref{l23} implies%
\begin{equation}
\int_{0}^{\infty}\lambda^{M}\left\vert \xi_{j}\left(  \phi\left(
\lambda\right)  \right)  -\dfrac{1}{2}\right\vert d\lambda<\infty\text{.}
\label{172}%
\end{equation}
Since $m_{j}\left(  \phi\left(  z\right)  \right)  $ ($j=1$, $2$) are Herglotz
functions and its argument on $\mathbb{R}$ is $\pi\xi_{j}\left(  \phi\left(
\lambda\right)  \right)  $, applying the formula (\ref{81}) to $m_{j}\left(
\phi\left(  z\right)  \right)  $ yields%
\[
\left\{
\begin{array}
[c]{c}%
m_{1}(\phi\left(  z\right)  )=\dfrac{1}{2\sqrt{-z}}\exp\left(
{\displaystyle\int_{0}^{\infty}}
\dfrac{\xi_{1}\left(  \phi\left(  \lambda\right)  \right)  -1/2}{\lambda
-z}d\lambda\right)  \smallskip\\
m_{2}(\phi\left(  z\right)  )=-\dfrac{\sqrt{-z}}{2}\exp\left(
{\displaystyle\int_{0}^{\infty}}
\dfrac{\xi_{2}\left(  \phi\left(  \lambda\right)  \right)  -1/2}{\lambda
-z}d\lambda\right)
\end{array}
\right.  \text{.}%
\]
We have assumed here $-a_{k}^{2}\leq\lambda_{1}$ for simplicity. Then for
$z\in\mathbb{C}_{+}$%
\[
m_{1}(z)=\dfrac{1}{2\sqrt{-\phi^{-1}\left(  z\right)  }}\exp\left(
{\displaystyle\int_{0}^{\infty}}
\dfrac{f\left(  \lambda\right)  }{\lambda-\phi^{-1}\left(  z\right)  }%
d\lambda\right)
\]
with $f\left(  \lambda\right)  =\xi_{1}\left(  \phi\left(  \lambda\right)
\right)  -1/2$. Lemma \ref{l51} implies%
\begin{equation}
\phi^{-1}\left(  z\right)  =z-g_{1}\left(  -z\right)  -\left(  -z\right)
^{-\left(  n_{1}-1\right)  /2+1/2}g_{2}(-z) \label{97}%
\end{equation}
with functions $g_{j}$ analytic near $z=\infty$ taking real values on
$\mathbb{R}$, hence%
\[
u\left(  z\right)  \equiv\sum_{k=0}^{M-1}\phi^{-1}\left(  z^{2}\right)
^{-k-1}\int_{0}^{\infty}\lambda^{k}f\left(  \lambda\right)  d\lambda
\]
is analytic at $z=\infty$, and the identity%
\[
-%
{\displaystyle\int_{0}^{\infty}}
\dfrac{f\left(  \lambda\right)  }{\lambda-\phi^{-1}\left(  z\right)  }%
d\lambda=u\left(  \sqrt{z}\right)  +\phi^{-1}\left(  z\right)  ^{-M}\int
_{0}^{\infty}\dfrac{\lambda^{M}f\left(  \lambda\right)  }{\lambda-\phi
^{-1}\left(  z\right)  }d\lambda
\]
holds. The estimate%
\[
\left\vert \int_{\lambda_{0}}^{\infty}\dfrac{\lambda^{M}f\left(
\lambda\right)  }{\lambda-\phi^{-1}\left(  z\right)  }d\lambda\right\vert
\leq\dfrac{1}{\left\vert \operatorname{Im}\phi^{-1}\left(  z\right)
\right\vert }\int_{\lambda_{0}}^{\infty}\left\vert \lambda^{M}f\left(
\lambda\right)  \right\vert d\lambda
\]
for $z\in\mathbb{C}_{+}$. Since $n_{1}>n$, (\ref{97}) shows for some $c>0$%
\[
\left\vert \operatorname{Im}\phi^{-1}\left(  z\right)  \right\vert \geq
c\left\vert z\right\vert ^{-n/2+1}\text{ \ \ if }z\in\widehat{D}_{-}\text{.}%
\]
Therefore%
\[
m_{1}(z)=\dfrac{1}{2\sqrt{-z}}\left(  \sum_{j=1}^{2M-1}\widetilde{a}_{j}%
\sqrt{-z}^{-j}+O\left(  z^{-M+n/2-1}\right)  \right)
\]
on $\widehat{D}_{-}$. However the assumption $q\in C^{2M-n}\left(
-\delta,\delta\right)  $ and (\ref{89}) imply%
\[
m_{1}(z)=\dfrac{1}{2\sqrt{-z}}\left(  1+\sum_{k=1}^{M-\left(  n-1\right)
/2}a_{k}\left(  -z\right)  ^{-k}+O\left(  z^{-M+n/2-1}\right)  \right)
\]
on a sector $\left\{  \epsilon<\arg z<\pi-\epsilon\right\}  $, hence
$\widetilde{a}_{j}=0$ for even $j$, which implies%
\[
m_{1}(-z^{2})=\dfrac{1}{2z}\left(  1+\sum_{k=1}^{M-\left(  n-1\right)
/2}a_{k}z^{-2k}+O\left(  z^{-2M+n-2}\right)  \right)
\]
on $D_{-}$. A similar calculation for $m_{2}(-z^{2})$ is possible and one can
obtain%
\[
m_{2}(-z^{2})=-\dfrac{1}{2}z\left(  1+\sum_{k=1}^{M-\left(  n-1\right)
/2}b_{k}z^{-2k}+O\left(  z^{-2M+n-2}\right)  \right)  \text{,}%
\]
which together with Lemma \ref{l24} for $L=2M+2-n$ completes the
proof.\medskip
\end{proof}

To apply Lemma \ref{l26} to ergodic potentials we need a lemma. The necessary
terminologies can be found in Appendix.

\begin{lemma}
\label{l31} Suppose the Lyapunov exponent $\gamma\left(  \lambda\right)  $
satisfies%
\begin{equation}
\int_{0}^{\infty}\lambda^{m}\gamma\left(  \lambda\right)  d\lambda
<\infty\label{173}%
\end{equation}
for some $m>4$. Then, for a.e. $\omega\in\Omega$ the condition (\ref{171}) is
fulfilled on the curve $\widehat{C}_{1}$ by any integer $M$ such that%
\begin{equation}
M<\min\left\{  \frac{n_{1}}{4}-1,\frac{m-n_{1}}{2}\right\}  \text{.}
\label{98}%
\end{equation}

\end{lemma}

\begin{proof}
Set%
\[
\left\{
\begin{array}
[c]{l}%
\rho\left(  \lambda\right)  =\sqrt{-\lambda}N\left(  \lambda\right)
I_{\left[  \lambda_{0},0\right]  }\left(  \lambda\right)  +\dfrac{1}{\pi}%
\sqrt{-\lambda}\gamma\left(  \lambda\right)  I_{\left(  0,\infty\right)
}\left(  \lambda\right)  \smallskip\\
c=\mathbb{E}\left(  q_{\omega}(0)/2\right)  \smallskip\\
w\left(  z\right)  =\mathbb{E}\left(  m_{\pm}\left(  z,\omega\right)  \right)
\end{array}
\right.  \text{.}%
\]
Since $\sqrt{-z}w(z)$ is of Herglotz, one has%
\[
w\left(  z\right)  =-\sqrt{-z}-\dfrac{c}{\sqrt{-z}}+\dfrac{1}{\sqrt{-z}}%
\int_{\lambda_{0}}^{\infty}\dfrac{\rho\left(  \lambda\right)  }{\lambda
-z}d\lambda\text{,}%
\]
hence%
\begin{equation}
w\left(  z\right)  =\sum_{k=-1}^{m_{1}}w_{k}(z)+\widetilde{w}(z) \label{174}%
\end{equation}
with%
\[
\left\{
\begin{array}
[c]{l}%
w_{k}(z)=c_{k}\left(  -z\right)  ^{-k-1/2}\smallskip\\
c_{-1}=-1\text{, }c_{0}=-c\text{, }c_{k}=\left(  -1\right)  ^{k-1}%
\int_{\lambda_{0}}^{\infty}\lambda^{k-1}\rho\left(  \lambda\right)
d\lambda\smallskip\\
\widetilde{w}(z)=\left(  -1\right)  ^{m_{1}}\left(  -z\right)  ^{-m_{1}%
-1/2}\int_{\lambda_{0}}^{\infty}\dfrac{\lambda^{m_{1}}\rho\left(
\lambda\right)  }{\lambda-z}d\lambda
\end{array}
\right.
\]
holds due to the assumption (\ref{173}), where $m_{1}=\left[  m\right]  $. Set%
\[
\left\{
\begin{array}
[c]{l}%
\chi(z)=\dfrac{-\operatorname{Re}w(z)}{\operatorname{Im}z}-\operatorname{Im}%
w^{\prime}(z)\smallskip\\
\chi_{k}(z)=\dfrac{-\operatorname{Re}w_{k}(z)}{\operatorname{Im}%
z}-\operatorname{Im}w_{k}^{\prime}(z)\smallskip\\
\widetilde{\chi}(z)=\dfrac{-\operatorname{Re}\widetilde{w}(z)}%
{\operatorname{Im}z}-\operatorname{Im}\widetilde{w}^{\prime}(z)
\end{array}
\right.  \text{.}%
\]
Noting%
\begin{align*}
\left(  -x-iy\right)  ^{-k-1/2}  &  =ix^{-k-1/2}\left(  -1\right)  ^{k}\left(
1+iyx^{-1}\right)  ^{-k-1/2}\\
&  =ix^{-k-1/2}\left(  -1\right)  ^{k}\left(  1-\left(  k+1/2\right)
iyx^{-1}+O\left(  yx^{-1}\right)  ^{2}\right)
\end{align*}
for $x\geq1$, $0<y<1$, we have%
\[
\chi_{k}\left(  x+iy\right)  =O\left(  x^{-k-5/2}y\right)  \text{.}%
\]
On the other hand, the estimates%
\[
\left\vert \int_{\lambda_{0}}^{\infty}\dfrac{\lambda^{m_{1}}\rho\left(
\lambda\right)  }{\left(  \lambda-z\right)  ^{j}}d\lambda\right\vert \leq
y^{-j}\int_{\lambda_{0}}^{\infty}\left\vert \lambda\right\vert ^{m_{1}}%
\rho\left(  \lambda\right)  d\lambda\text{ \ (}j=1\text{, }2\text{)}%
\]
yield a bound for the last term $\widetilde{w}$ of (\ref{174}). Then, we have%
\[
\widetilde{\chi}\left(  x+iy\right)  =O\left(  y^{-2}x^{-m_{1}-1/2}\right)
\text{,}%
\]
imply%
\[
\chi(x+iy)=\sum_{k=-1}^{m_{1}}\chi_{k}(x+iy)+\widetilde{\chi}\left(
x+iy\right)  =O\left(  yx^{-3/2}+y^{-2}x^{-m_{1}-1/2}\right)  \text{.}%
\]
This together with $\operatorname{Im}w\left(  x+iy\right)  =O\left(
x^{1/2}\right)  $ (due to $N\left(  \lambda\right)  \sim\sqrt{\lambda}$ as
$\lambda\rightarrow\infty$) yields%
\[
\sqrt{2\chi\left(  z\right)  \operatorname{Im}w\left(  z\right)  }=O\left(
yx^{-1}+y^{-2}x^{-m_{1}}\right)  ^{1/2}\text{.}%
\]
Therefore, if the curve is parametrized as $x+ix^{-\left(  n_{1}/2-1\right)
}$ near $x=\infty$, applying (\ref{90}) we have%
\begin{align*}
\mathbb{E}\left(  \int_{\widehat{C}}\left\vert z\right\vert ^{M}\left\vert
R\left(  z,\omega\right)  \right\vert \left\vert dz\right\vert \right)   &
\leq\int_{\widehat{C}}\left\vert z\right\vert ^{M}\sqrt{2\chi\left(  z\right)
\operatorname{Im}w\left(  z\right)  }\left\vert dz\right\vert \\
&  \leq c\int_{1}^{\infty}x^{M}\left(  x^{-\left(  n_{1}/2-1\right)  }%
x^{-1}+x^{\left(  n_{1}-2\right)  }x^{-m_{1}}\right)  ^{1/2}dx\text{,}%
\end{align*}
which is finite for $M$ such that%
\[
M-\left(  n_{1}/2-1\right)  /2-1/2<-1\text{\ and\ }M+\left(  n_{1}/2-1\right)
-m_{1}/2<-1\text{.}%
\]
Then, Fubini's theorem implies the condition (\ref{171}).\medskip
\end{proof}

Now one can prove Theorem \ref{t4}. Suppose the Lyapunov exponent
$\gamma\left(  \lambda\right)  $ satisfies%
\begin{equation}
\int_{0}^{\infty}\lambda^{m}\gamma\left(  \lambda\right)  d\lambda
<\infty\text{.} \label{178}%
\end{equation}
Since and Lemma \ref{l26} and (\ref{98}) require%
\[
L+\left(  n-1\right)  /2-1\leq M<\min\left(  \dfrac{n_{1}}{4}-1,\dfrac
{m-n_{1}}{2}\right)  \text{,}%
\]
$n_{1}$ should satisfy%
\begin{equation}
4L+2n-2<n_{1}<m-2L-n+3\text{.} \label{99}%
\end{equation}
If%
\[
m-2L-n+3-\left(  4L+2n-2\right)  =m-6L-3n+5>2\text{,}%
\]
one can choose an odd integer $n_{1}$ satisfying (\ref{99}). Then applying
Lemma \ref{l26} and Lemma \ref{l31} we have $q_{\omega}\in\mathcal{Q}%
_{L}\left(  C\right)  $ for a.e. $\omega$ if $L<\left(  m-3\left(  n-1\right)
\right)  /6$. On the other hand from
\[
q_{\theta_{x}\omega}\left(  y\right)  =q_{\omega}\left(  x+y\right)  =\left(
K\left(  e_{x}\right)  q_{\omega}\right)  \left(  y\right)
\]
the identity%
\[
f_{g}\left(  \theta_{x}\omega\right)  =\left(  K(g)q_{\theta_{x}\omega
}\right)  \left(  0\right)  =\left(  K(g)K\left(  e_{x}\right)  q_{\omega
}\right)  \left(  0\right)  =\left(  K\left(  e_{x}\right)  K(g)q_{\omega
}\right)  \left(  0\right)  =K(g)q_{\omega}\left(  x\right)
\]
follows. Moreover Kotani-Krishna \cite{k2} showed that $q_{\omega}\in
C_{b}^{m}\left(  \mathbb{R}\right)  $ implies
\[
\int_{0}^{\infty}\lambda^{m+1/2}\gamma\left(  \lambda\right)  d\lambda
<\infty\text{,}%
\]
which is sufficient for (\ref{178}) and completes the proof of Theorem
\ref{t4}.

\section{Appendix}

\subsection{Calculation of $\tau_{\boldsymbol{a}}\left(  r\right)  $,
$m_{r\boldsymbol{a}}\left(  z\right)  $ for rational functions $r$}

This section is devoted to the calculation of $\tau_{\boldsymbol{a}}\left(
r\right)  $, $m_{r\boldsymbol{a}}\left(  z\right)  $ for general rational
function $r$ in terms of the characteristic functions $\left\{  \varphi
_{\boldsymbol{a}},\psi_{\boldsymbol{a}}\right\}  $ and $m$-function
$m_{\boldsymbol{a}}$.

The simplest rational functions are%
\[
\left\{
\begin{array}
[c]{l}%
p_{\zeta}\left(  z\right)  =1+z\zeta^{-1}\in\Gamma_{0}^{\left(  1\right)  }\\
q_{\zeta}\left(  z\right)  =\left(  1-z\zeta^{-1}\right)  ^{-1}\in\Gamma
_{0}^{\left(  -1\right)  }%
\end{array}
\text{ for }\zeta\in D_{-}\right.  \text{.}%
\]
Any rational function can be represented as a product of these simple functions.

We treat $r\in\Gamma_{0}^{\left(  m\right)  }$ with $m\leq0$. Recall%
\[
\tau_{\boldsymbol{a}}\left(  g\right)  =\det\left(  g^{-1}T\left(
g\boldsymbol{a}\right)  T(\boldsymbol{a})^{-1}\right)
\]
for $\boldsymbol{a}\in\boldsymbol{A}_{L}^{inv}\left(  C\right)  $, which is
well-defined if $g^{-1}T\left(  g\boldsymbol{a}\right)  T(\boldsymbol{a}%
)^{-1}-I$ is of trace class in some space $H_{N}\left(  D_{+}\right)  $. Let
$M=L-N\geq0$ and observe for $u\in H_{N}\left(  D_{+}\right)  $%
\begin{align*}
T\left(  g\boldsymbol{a}\right)  T(\boldsymbol{a})^{-1}u  &  =\mathfrak{p}%
_{+}\left(  g\mathfrak{p}_{+}\boldsymbol{a}T(\boldsymbol{a})^{-1}u\right)
+\mathfrak{p}_{+}\left(  g\mathfrak{p}_{-}\boldsymbol{a}T(\boldsymbol{a}%
)^{-1}u\right) \\
&  =gu+H_{g}\mathfrak{p}_{-}\left(  \boldsymbol{a}T(\boldsymbol{a}%
)^{-1}u\right)
\end{align*}
with $H_{g}:H_{M}\left(  D_{-}\right)  \left(  =z^{-M}H\left(  D_{-}\right)
\right)  \rightarrow H\left(  D_{+}\right)  $ defined by%
\begin{equation}
H_{g}v\left(  z\right)  =\mathfrak{p}_{+}\left(  gv\right)  \left(  z\right)
=\dfrac{1}{2\pi i}\int_{C}\dfrac{g\left(  \lambda\right)  }{\lambda-z}v\left(
\lambda\right)  d\lambda\text{.} \label{36}%
\end{equation}
Here $v\in H\left(  D_{-}\right)  $ is%
\[
v=\mathfrak{p}_{-}\left(  \boldsymbol{a}T(\boldsymbol{a})^{-1}u\right)
=\boldsymbol{a}T(\boldsymbol{a})^{-1}u-u\text{.}%
\]
The key identity is (\ref{62}):%
\[
T\left(  \boldsymbol{a}\right)  ^{-1}\dfrac{1}{z+b}=\dfrac{\left(
\psi_{\boldsymbol{a}}\left(  b\right)  +b\right)  u-\left(  \varphi
_{\boldsymbol{a}}\left(  b\right)  +1\right)  v}{\Delta_{\boldsymbol{a}%
}\left(  b\right)  \left(  z^{2}-b^{2}\right)  }\text{ for }b\in D_{-}\text{.}%
\]

\begin{lemma}
\label{l13}Let $\boldsymbol{a}\in\boldsymbol{A}_{L}^{inv}\left(  C\right)  $
and $r\in\Gamma_{0}^{\left(  m\right)  }$.\newline(i) \ For any $N$ such that
$-m\leq N\leq L$ the operator $r^{-1}T\left(  r\boldsymbol{a}\right)
T(\boldsymbol{a})^{-1}-I$ has a finite rank on $H_{N}\left(  D_{+}\right)  $
not greater than the numbers of the poles of $r$. Hence $\tau_{\boldsymbol{a}%
}\left(  r\right)  $ is well-defined for any $r\in\Gamma_{0}^{\left(
m\right)  }$.\newline(ii) Suppose $L\geq2$. For $\zeta$, $\zeta_{1}$,
$\zeta_{2}\in C\backslash\left(  \left[  -\mu_{0},\mu_{0}\right]  \cup
i\mathbb{R}\right)  $%
\begin{equation}
\left\{
\begin{array}
[c]{l}%
\tau_{\boldsymbol{a}}\left(  q_{\zeta}\right)  =1+\varphi_{\boldsymbol{a}%
}\left(  \zeta\right) \\
\tau_{\boldsymbol{a}}\left(  q_{\zeta_{1}}q_{\zeta_{2}}\right)  =\left(
1+\varphi_{\boldsymbol{a}}\left(  \zeta_{1}\right)  \right)  \left(
1+\varphi_{\boldsymbol{a}}\left(  \zeta_{2}\right)  \right)  \dfrac
{m_{\boldsymbol{a}}\left(  \zeta_{1}\right)  -m_{\boldsymbol{a}}\left(
\zeta_{2}\right)  }{\zeta_{1}-\zeta_{2}}%
\end{array}
\right.  \text{.} \label{66}%
\end{equation}
(iii) Suppose $L\geq2$. Suppose $r\in\Gamma_{0}^{\left(  0\right)  }$ has
simple zeroes $\left\{  \eta_{j}\right\}  _{1\leq j\leq n}$ and simple poles
$\left\{  \zeta_{j}\right\}  _{1\leq j\leq n}$ in $D_{-}$. Then $H_{r}$ has a
finite rank not greater than $n$, and it holds that%
\begin{align}
&  \tau_{\boldsymbol{a}}\left(  r\right) \nonumber\\
&  =\left(  \prod_{j=1}^{n}\dfrac{\left(  \varphi_{\boldsymbol{a}}\left(
\zeta_{j}\right)  +1\right)  \left(  \varphi_{\boldsymbol{a}}\left(  -\eta
_{j}\right)  +1\right)  }{\Delta_{\boldsymbol{a}}\left(  \eta_{j}\right)
r^{\prime}\left(  \eta_{j}\right)  \widehat{r}^{\prime}\left(  \zeta
_{j}\right)  }\right)  \det\left(  \dfrac{1}{\eta_{i}-\zeta_{j}}\right)
\det\left(  \dfrac{m_{\boldsymbol{a}}\left(  \zeta_{i}\right)
-m_{\boldsymbol{a}}\left(  -\eta_{j}\right)  }{\zeta_{i}^{2}-\eta_{j}^{2}%
}\right)  , \label{60}%
\end{align}
with $\widehat{r}\left(  z\right)  =r(z)^{-1}$.
\end{lemma}

\begin{proof}
Let $\left\{  \zeta_{j}\right\}  _{j=1}^{n}$ be the poles of $r$. For
simplicity assume they are simple. Then $r$ can be expressed as%
\[
r\left(  z\right)  =r\left(  \infty\right)  +\sum_{j=1}^{n}\dfrac{r_{j}%
}{z-\zeta_{j}}\text{ with }r_{j}=\lim_{z\rightarrow\zeta_{j}}\left(
z-\zeta_{j}\right)  r\left(  z\right)  =\dfrac{1}{\widehat{r}^{\prime}\left(
\zeta_{j}\right)  }\text{,}%
\]
hence for $v\in H\left(  D_{-}\right)  $%
\[
\left(  H_{r}v\right)  \left(  z\right)  =\dfrac{1}{2\pi i}\int_{C}%
\dfrac{r\left(  \lambda\right)  }{\lambda-z}v\left(  \lambda\right)
d\lambda=\sum_{j=1}^{n}\dfrac{r_{j}v\left(  \zeta_{j}\right)  }{z-\zeta_{j}}%
\]
and%
\begin{equation}
r(z)^{-1}\left(  H_{r}v\right)  \left(  z\right)  =\sum_{j=1}^{n}%
f_{j}(z)v\left(  \zeta_{j}\right)  \text{ \ with\ \ }f_{j}(z)=\dfrac{r_{j}%
}{\left(  z-\zeta_{j}\right)  r(z)}\text{,} \label{20}%
\end{equation}
which means that the map $r^{-1}T\left(  r\boldsymbol{a}\right)
T(\boldsymbol{a})^{-1}-I$ is of finite rank.

To compute $\tau_{\boldsymbol{a}}\left(  r\right)  $, we take the independent
vectors $f_{j}$ and obtain the coefficients of the image of $r^{-1}H_{r}$ for
$u=T\left(  \boldsymbol{a}\right)  ^{-1}f_{j}$.

If $L\geq2$, one can use the characteristic functions and $m$-function. Let
$r=q_{\zeta}$. Then $n=1$ and $f_{1}=1$, hence $r(z)^{-1}\left(
H_{r}v\right)  \left(  z\right)  =v\left(  \zeta\right)  $. For $u=1$%
\[
v=\boldsymbol{a}T(\boldsymbol{a})^{-1}1-1=\varphi_{\boldsymbol{a}}\text{,}%
\]
which yields $\tau_{\boldsymbol{a}}\left(  q_{\zeta}\right)  =1+\varphi
_{\boldsymbol{a}}\left(  \zeta\right)  $. If $r=q_{\zeta_{1}}q_{\zeta_{2}}$,
then $n=2$ and%
\[
f_{1}\left(  z\right)  =\dfrac{r_{1}}{\left(  z-\zeta_{1}\right)  r(z)}%
=\dfrac{\zeta_{2}-z}{\zeta_{2}-\zeta_{1}}\text{, \ }f_{2}\left(  z\right)
=\dfrac{r_{2}}{\left(  z-\zeta_{2}\right)  r(z)}=\dfrac{\zeta_{1}-z}{\zeta
_{1}-\zeta_{2}}\text{.}%
\]
For $u=f_{1}$, $u=f_{2}$%
\[
\left\{
\begin{array}
[c]{c}%
v_{1}=\boldsymbol{a}T(\boldsymbol{a})^{-1}f_{1}-f_{1}=\dfrac{\zeta_{2}}%
{\zeta_{2}-\zeta_{1}}\varphi_{\boldsymbol{a}}-\dfrac{1}{\zeta_{2}-\zeta_{1}%
}\psi_{\boldsymbol{a}}\\
v_{2}=\boldsymbol{a}T(\boldsymbol{a})^{-1}f_{2}-f_{2}=\dfrac{\zeta_{1}}%
{\zeta_{1}-\zeta_{2}}\varphi_{\boldsymbol{a}}-\dfrac{1}{\zeta_{1}-\zeta_{2}%
}\psi_{\boldsymbol{a}}%
\end{array}
\right.  \text{,}%
\]
hence%
\begin{align*}
\tau_{\boldsymbol{a}}\left(  q_{\zeta_{1}}q_{\zeta_{2}}\right)   &
=\det\left(
\begin{array}
[c]{cc}%
1+v_{1}\left(  \zeta_{1}\right)  & v_{1}\left(  \zeta_{2}\right) \\
v_{2}\left(  \zeta_{1}\right)  & 1+v_{2}\left(  \zeta_{2}\right)
\end{array}
\right) \\
&  =\dfrac{\left(  \zeta_{1}+\psi_{\boldsymbol{a}}\left(  \zeta_{1}\right)
\right)  \left(  1+\varphi_{\boldsymbol{a}}\left(  \zeta_{2}\right)  \right)
-\left(  1+\varphi_{\boldsymbol{a}}\left(  \zeta_{1}\right)  \right)  \left(
\zeta_{2}+\psi_{\boldsymbol{a}}\left(  \zeta_{2}\right)  \right)  }{\zeta
_{1}-\zeta_{2}}\text{,}%
\end{align*}
which is (\ref{66}).

Now go back to (\ref{60}). Since $f_{j}(z)$ is a rational function with poles
at $\eta_{i}$ and $f_{j}(\infty)=0$, an identity%
\[
f_{j}(z)=\dfrac{r_{j}}{\left(  z-\zeta_{j}\right)  r(z)}=\sum_{i=1}^{n}%
\dfrac{r_{ij}}{\left(  z-\zeta_{j}\right)  \left(  z-\eta_{i}\right)  }%
\]
with%
\[
r_{ij}=\lim_{z\rightarrow\eta_{i}}\dfrac{r_{j}\left(  z-\eta_{i}\right)
}{\left(  z-\zeta_{j}\right)  r(z)}=\dfrac{r_{j}}{\left(  \eta_{i}-\zeta
_{j}\right)  r^{\prime}\left(  \eta_{i}\right)  }%
\]
is valid. Then (\ref{62}) yields%
\[
T(\boldsymbol{a})^{-1}f_{j}=\sum_{i}r_{ij}T(\boldsymbol{a})^{-1}\dfrac
{1}{z-\eta_{i}}=\sum_{i}r_{ij}\dfrac{\left(  \varphi_{\boldsymbol{a}}\left(
-\eta_{i}\right)  +1\right)  v-\left(  \psi_{\boldsymbol{a}}\left(  -\eta
_{i}\right)  -\eta_{i}\right)  u}{\Delta_{\boldsymbol{a}}\left(  \eta
_{i}\right)  \left(  z^{2}-\eta_{i}^{2}\right)  }%
\]
and%
\begin{align*}
\boldsymbol{a}T(\boldsymbol{a})^{-1}f_{j}  &  =\sum_{i}r_{ij}\dfrac{\left(
\varphi_{\boldsymbol{a}}\left(  -\eta_{i}\right)  +1\right)  \left(
\psi_{\boldsymbol{a}}+z\right)  -\left(  \psi_{\boldsymbol{a}}\left(
-\eta_{i}\right)  -\eta_{i}\right)  \left(  \varphi_{\boldsymbol{a}}+1\right)
}{\Delta_{\boldsymbol{a}}\left(  \eta_{i}\right)  \left(  z^{2}-\eta_{i}%
^{2}\right)  }\\
&  =\left(  \varphi_{\boldsymbol{a}}+1\right)  \sum_{i}r_{ij}\left(
\varphi_{\boldsymbol{a}}\left(  -\eta_{i}\right)  +1\right)  \dfrac
{m_{\boldsymbol{a}}-m_{\boldsymbol{a}}\left(  -\eta_{i}\right)  }%
{\Delta_{\boldsymbol{a}}\left(  \eta_{i}\right)  \left(  z^{2}-\eta_{i}%
^{2}\right)  }\text{.}%
\end{align*}
Therefore, noting $f_{j}(\zeta_{i})=\delta_{ij}$, we have%
\begin{align*}
\tau_{\boldsymbol{a}}\left(  r\right)   &  =\det\left(  \left(  \delta
_{ij}+\left(  \boldsymbol{a}T(\boldsymbol{a})^{-1}f_{j}\right)  \left(
\zeta_{i}\right)  -f_{j}\left(  \zeta_{i}\right)  \right)  _{1\leq i,j\leq
n}\right) \\
&  =\det\left(  \left(  \boldsymbol{a}T(\boldsymbol{a})^{-1}f_{j}\right)
\left(  \zeta_{i}\right)  \right) \\
&  =\det\left(  \left(  \varphi_{\boldsymbol{a}}\left(  \zeta_{i}\right)
+1\right)  \sum_{k=1}^{n}r_{kj}\dfrac{m_{\boldsymbol{a}}\left(  \zeta
_{i}\right)  -m_{\boldsymbol{a}}\left(  -\eta_{k}\right)  }{\zeta_{i}^{2}%
-\eta_{k}^{2}}\dfrac{\varphi_{\boldsymbol{a}}\left(  -\eta_{k}\right)
+1}{\Delta_{\boldsymbol{a}}\left(  \eta_{k}\right)  }\right)  \text{.}\\
&  =\left(  \prod_{j=1}^{n}\dfrac{\left(  \varphi_{\boldsymbol{a}}\left(
\zeta_{i}\right)  +1\right)  \left(  \varphi_{\boldsymbol{a}}\left(  -\eta
_{j}\right)  +1\right)  }{\Delta_{\boldsymbol{a}}\left(  \eta_{j}\right)
}\right)  \det\left(  r_{ij}\right)  \det\left(  \dfrac{m_{\boldsymbol{a}%
}\left(  \zeta_{i}\right)  -m_{\boldsymbol{a}}\left(  -\eta_{j}\right)
}{\zeta_{i}^{2}-\eta_{j}^{2}}\right)  \text{,}%
\end{align*}
where%
\[
\det\left(  r_{ij}\right)  =\det\left(  \dfrac{1}{\left(  \eta_{i}-\zeta
_{j}\right)  r^{\prime}\left(  \eta_{i}\right)  \widehat{r}^{\prime}\left(
\zeta_{j}\right)  }\right)  =\left(  \prod_{j=1}^{n}\dfrac{1}{r^{\prime
}\left(  \eta_{i}\right)  \widehat{r}^{\prime}\left(  \zeta_{j}\right)
}\right)  \det\left(  \dfrac{1}{\eta_{i}-\zeta_{j}}\right)  \text{,}%
\]
which is (\ref{60}).\medskip
\end{proof}

It should be remarked that if $\eta_{i}\neq\zeta_{j}$, then%
\[
\det\left(  \dfrac{1}{\eta_{i}-\zeta_{j}}\right)  \neq0\text{.}%
\]
This is because the identity%
\[
0=\sum_{j=1}^{n}\dfrac{u_{j}}{\eta_{i}-\zeta_{j}}\ \ \text{for any }i
\]
implies the rational function $f(z)=\sum_{j=1}^{n}u_{j}\left(  z-\zeta
_{j}\right)  ^{-1}$ satisfies $f\left(  \eta_{i}\right)  =0$ for $i=1$, $2$,
$\cdots$, $n$, which shows $f(z)=0$ identically.

The next task is to compute $m_{r\boldsymbol{a}}$ for rational functions $r$.
In principle the computation for general rational $r$ is possible similarly as
the previous lemma, however to grasp the picture it is enough to know the
change of $m_{r\boldsymbol{a}}$ for $r=q_{\zeta}p_{\eta}$, since the formula
of $m_{r\boldsymbol{a}}$ for general $r$ can be obtained by iteration of
$q_{\zeta}p_{\eta}$.

\begin{lemma}
\label{l16}Let $\boldsymbol{a}\in\boldsymbol{A}_{2}^{inv}$ and $\zeta$,
$\eta\in D_{-}$ and assume $\tau_{\boldsymbol{a}}\left(  q_{\zeta}p_{\eta
}\right)  \neq0$. Then%
\[
m_{q_{\zeta}p_{\eta}\boldsymbol{a}}\left(  z\right)  =\left(  d_{\zeta}%
d_{\eta}m_{\boldsymbol{a}}\right)  \left(  z\right)  \text{.}%
\]

\end{lemma}

\begin{proof}
Let $r=q_{\zeta}p_{\eta}$ with $\zeta$, $\eta\in D_{-}$. First we have to
compute $\varphi_{r\boldsymbol{a}}$, $\psi_{r\boldsymbol{a}}.$ Then%
\begin{equation}
r(z)^{-1}=r_{1}+\dfrac{r_{2}}{z+\eta}\text{ \ with \ }r_{1}=-\dfrac{\eta
}{\zeta}\text{, \ }r_{2}=\eta\left(  1+\dfrac{\eta}{\zeta}\right)  \label{67}%
\end{equation}
holds. To compute $\varphi_{r\boldsymbol{a}}$ set $w_{1}=T\left(
r\boldsymbol{a}\right)  ^{-1}1$. The definition implies $1+\varphi
_{r\boldsymbol{a}}=raw_{1}$, hence (\ref{67}) yields%
\begin{align}
\left(  \boldsymbol{a}w_{1}\right)  \left(  z\right)   &  =r\left(  z\right)
^{-1}\left(  1+\varphi_{r\boldsymbol{a}}\left(  z\right)  \right) \label{68}\\
&  =r_{1}+\dfrac{r_{2}\left(  1+\varphi_{r\boldsymbol{a}}\left(  -\eta\right)
\right)  }{z+\eta}+\dfrac{r_{2}\left(  \varphi_{r\boldsymbol{a}}\left(
z\right)  -\varphi_{r\boldsymbol{a}}\left(  -\eta\right)  \right)  }{z+\eta
}\text{,}\nonumber
\end{align}
which is a decomposition in $H_{1}\left(  D_{+}\right)  \oplus H\left(
D_{-}\right)  $. Applying $\mathfrak{p}_{+}$, we have%
\[
T\left(  \boldsymbol{a}\right)  w_{1}=r_{1}+\dfrac{r_{2}\left(  1+\varphi
_{r\boldsymbol{a}}\left(  -\eta\right)  \right)  }{z+\eta}\text{.}%
\]
Therefore (\ref{62}) implies%
\begin{align*}
w_{1}  &  =r_{1}T\left(  \boldsymbol{a}\right)  ^{-1}1+r_{2}\left(
1+\varphi_{r\boldsymbol{a}}\left(  -\eta\right)  \right)  T\left(
\boldsymbol{a}\right)  ^{-1}\dfrac{1}{z+\eta}\\
&  =r_{1}u+\mu_{1}\dfrac{\left(  m_{\boldsymbol{a}}\left(  \eta\right)
-\kappa_{1}\left(  \boldsymbol{a}\right)  \right)  u-v}{z^{2}-\eta^{2}}%
\end{align*}
with%
\[
\mu_{1}=\dfrac{r_{2}\left(  1+\varphi_{\boldsymbol{a}}\left(  \eta\right)
\right)  \left(  1+\varphi_{r\boldsymbol{a}}\left(  -\eta\right)  \right)
}{\Delta_{\boldsymbol{a}}\left(  \eta\right)  }\text{,}%
\]
where $u=T\left(  \boldsymbol{a}\right)  ^{-1}1$, $v=T\left(  \boldsymbol{a}%
\right)  ^{-1}z$, from which it follows that%
\begin{equation}
\boldsymbol{a}w_{1}=\left(  1+\varphi_{\boldsymbol{a}}\right)  \left(
r_{1}-\mu_{1}\dfrac{m_{\boldsymbol{a}}-m_{\boldsymbol{a}}\left(  \eta\right)
}{z^{2}-\eta^{2}}\right)  \text{.} \label{69}%
\end{equation}
The identity (\ref{68}) shows the left hand side is meromorphic on $D_{-}$
vanishing at $z=\zeta$, hence%
\[
\mu_{1}\dfrac{m_{\boldsymbol{a}}\left(  \zeta\right)  -m_{\boldsymbol{a}%
}\left(  \eta\right)  }{\zeta^{2}-\eta^{2}}=r_{1}%
\]
holds, which together with (\ref{69}) shows%
\begin{equation}
1+\varphi_{r\boldsymbol{a}}=r_{1}r\left(  1+\varphi_{\boldsymbol{a}}\right)
\left(  1-\dfrac{\zeta^{2}-\eta^{2}}{m_{\boldsymbol{a}}\left(  \zeta\right)
-m_{\boldsymbol{a}}\left(  \eta\right)  }\dfrac{m_{\boldsymbol{a}%
}-m_{\boldsymbol{a}}\left(  \eta\right)  }{z^{2}-\eta^{2}}\right)  \text{.}
\label{70}%
\end{equation}
This identity also shows%
\begin{equation}
\kappa_{1}\left(  r\boldsymbol{a}\right)  =\lim_{z\rightarrow\infty}%
z\varphi_{r\boldsymbol{a}}\left(  z\right)  =\zeta+\eta+\kappa_{1}\left(
\boldsymbol{a}\right)  -\dfrac{\zeta^{2}-\eta^{2}}{m_{\boldsymbol{a}}\left(
\zeta\right)  -m_{\boldsymbol{a}}\left(  \eta\right)  }\text{.} \label{71}%
\end{equation}
Similarly one has $\psi_{r\boldsymbol{a}}$, namely $w_{2}=T\left(
r\boldsymbol{a}\right)  ^{-1}z$ satisfies%
\begin{align*}
&  \left(  \boldsymbol{a}w_{2}\right)  \left(  z\right) \\
&  =r\left(  z\right)  ^{-1}\left(  z+\psi_{r\boldsymbol{a}}\left(  z\right)
\right) \\
&  =r_{1}z+r_{2}+r_{2}\dfrac{-\eta+\psi_{r\boldsymbol{a}}\left(  -\eta\right)
}{z+\eta}+r_{2}\dfrac{\psi_{r\boldsymbol{a}}\left(  z\right)  -\psi
_{r\boldsymbol{a}}\left(  -\eta\right)  }{z+\eta}\text{,}%
\end{align*}
hence%
\[
T\left(  \boldsymbol{a}\right)  w_{2}=r_{1}z+r_{2}+r_{2}\dfrac{-\eta
+\psi_{r\boldsymbol{a}}\left(  -\eta\right)  }{z+\eta}\text{,}%
\]
which yields%
\[
w_{2}=r_{1}v+r_{2}u+\mu_{2}\dfrac{\left(  m_{\boldsymbol{a}}\left(
\eta\right)  -\kappa_{1}\left(  \boldsymbol{a}\right)  \right)  u-v}%
{z^{2}-\eta^{2}}%
\]
with%
\[
\mu_{2}=\dfrac{r_{2}\left(  \varphi_{\boldsymbol{a}}\left(  \eta\right)
+1\right)  \left(  -\eta+\psi_{r\boldsymbol{a}}\left(  -\eta\right)  \right)
}{\Delta_{\boldsymbol{a}}\left(  \eta\right)  }\text{.}%
\]
Then%
\[
\boldsymbol{a}w_{2}=r_{1}\left(  z+\psi_{\boldsymbol{a}}\right)  +\left(
1+\varphi_{\boldsymbol{a}}\right)  \left(  r_{2}-\mu_{2}\dfrac
{m_{\boldsymbol{a}}-m_{\boldsymbol{a}}\left(  \eta\right)  }{z^{2}-\eta^{2}%
}\right)
\]
with%
\[
\mu_{2}=\dfrac{r_{2}\left(  \varphi_{\boldsymbol{a}}\left(  \eta\right)
+1\right)  \left(  -\eta+\psi_{r\boldsymbol{a}}\left(  -\eta\right)  \right)
}{\Delta_{\boldsymbol{a}}\left(  \eta\right)  }\text{,}%
\]
hence setting $z=\zeta$ we have%
\[
r_{1}\left(  \zeta+\psi_{\boldsymbol{a}}\left(  \zeta\right)  \right)
+\left(  1+\varphi_{\boldsymbol{a}}\left(  \zeta\right)  \right)  \left(
r_{2}-\mu_{2}\dfrac{m_{\boldsymbol{a}}\left(  \zeta\right)  -m_{\boldsymbol{a}%
}\left(  \eta\right)  }{\zeta^{2}-\eta^{2}}\right)  =0\text{,}%
\]
which yields%
\[
\mu_{2}=r_{2}\dfrac{\zeta^{2}-\eta^{2}}{m_{\boldsymbol{a}}\left(
\zeta\right)  -m_{\boldsymbol{a}}\left(  \eta\right)  }\left(  1-\dfrac
{m_{\boldsymbol{a}}\left(  \zeta\right)  -\kappa_{1}\left(  \boldsymbol{a}%
\right)  }{\zeta+\eta}\right)  \text{,}%
\]
and%
\begin{align}
&  z+\psi_{r\boldsymbol{a}}\left(  z\right) \nonumber\\
&  =\dfrac{r_{2}}{r_{1}}\left(  1-\dfrac{m_{\boldsymbol{a}}\left(
\zeta\right)  -\kappa_{1}\left(  \boldsymbol{a}\right)  }{\zeta+\eta}\right)
\left(  1+\varphi_{r\boldsymbol{a}}\left(  z\right)  \right)  -r_{2}r\left(
1+\varphi_{\boldsymbol{a}}\right)  \dfrac{m_{\boldsymbol{a}}-m_{\boldsymbol{a}%
}\left(  \zeta\right)  }{\zeta+\eta} \label{65}%
\end{align}
Consequently $m_{r\boldsymbol{a}}$ is computed from (\ref{69}), (\ref{71}) and
(\ref{65}) as%
\begin{align*}
m_{r\boldsymbol{a}}\left(  z\right)   &  =\dfrac{z+\psi_{r\boldsymbol{a}%
}\left(  z\right)  }{1+\varphi_{r\boldsymbol{a}}\left(  z\right)  }+\kappa
_{1}\left(  r\boldsymbol{a}\right) \\
&  =\dfrac{\left(  m_{\boldsymbol{a}}\left(  z\right)  -m_{\boldsymbol{a}%
}\left(  \zeta\right)  \right)  \dfrac{m_{\boldsymbol{a}}\left(  \zeta\right)
-m_{\boldsymbol{a}}\left(  \eta\right)  }{\zeta^{2}-\eta^{2}}}{\dfrac
{m_{\boldsymbol{a}}\left(  \zeta\right)  -m_{\boldsymbol{a}}\left(
\eta\right)  }{\zeta^{2}-\eta^{2}}-\dfrac{m_{\boldsymbol{a}}\left(  z\right)
-m_{\boldsymbol{a}}\left(  \eta\right)  }{z^{2}-\eta^{2}}}+m_{\boldsymbol{a}%
}\left(  \zeta\right)  -\dfrac{m_{\boldsymbol{a}}\left(  \zeta\right)
-m_{\boldsymbol{a}}\left(  \eta\right)  }{\zeta^{2}-\eta^{2}}\\
&  =\left(  d_{\zeta}d_{\eta}m_{\boldsymbol{a}}\right)  \left(  z\right)
\end{align*}
\medskip
\end{proof}

This formula can be easily understood if we first show the identity
$m_{q_{\zeta}\boldsymbol{a}}=d_{\zeta}m_{\boldsymbol{a}}$. The reason why we
do not take this procedure is that we have defined $m_{\boldsymbol{a}}$ for
$\boldsymbol{a}$ such that $T\left(  \boldsymbol{a}\right)  $ maps
$H_{N}\left(  D_{+}\right)  $ to $H_{N}\left(  D_{+}\right)  $ bijectively,
hence $m_{q_{\zeta}\boldsymbol{a}}$, $m_{p_{\zeta}\boldsymbol{a}}$ are out of
the present framework. However, a slight modification of the definition of
$m_{\boldsymbol{a}}$ might allow us to show%
\[
m_{q_{\zeta}\boldsymbol{a}}=m_{p_{\zeta}\boldsymbol{a}}=d_{\zeta
}m_{\boldsymbol{a}}\text{,}%
\]
and the identity of the Lemma would be more straightly understandable.

\subsection{Properties of $m$-functions and Herglotz functions}

The $m$-function $m_{\boldsymbol{a}}$ for $\boldsymbol{a}\in\boldsymbol{A}%
_{L,+}^{inv}\left(  C\right)  $ has the following properties:%
\begin{equation}%
\begin{tabular}
[c]{l}%
$m_{\boldsymbol{a}}$ is analytic on $\mathbb{C}\backslash\left(  \left[
-\mu_{0},\mu_{0}\right]  \cup i\mathbb{R}\right)  $\ and $m_{\boldsymbol{a}%
}\left(  z\right)  =\overline{m_{\boldsymbol{a}}\left(  \overline{z}\right)
}$\\
$\dfrac{\operatorname{Im}m_{\boldsymbol{a}}\left(  z\right)  }%
{\operatorname{Im}z}>0\text{ for }z\in\mathbb{C}\backslash\mathbb{R}\text{,
\ }\dfrac{m_{\boldsymbol{a}}\left(  x\right)  -m_{\boldsymbol{a}}\left(
-x\right)  }{x}>0\text{ \ for }\left\vert x\right\vert >\mu_{0}\text{.}$%
\end{tabular}
\ \label{47}%
\end{equation}
with $\mu_{0}=\sqrt{-\lambda_{0}}$. For such function $m$ define%
\[
m_{+}\left(  z\right)  =-m\left(  \sqrt{-z}\right)  \text{, \ \ }m_{-}\left(
z\right)  =m\left(  -\sqrt{-z}\right)  \text{ \ for }z\in\mathbb{C}%
_{+}\text{.}%
\]
Then $m_{\pm}$ become Herglotz functions, namely they analytic functions on
$\mathbb{C}\backslash\lbrack\lambda_{0},\infty)$ satisfying%
\[%
\begin{tabular}
[c]{l}%
$m_{\pm}\left(  z\right)  =\overline{m_{\pm}\left(  \overline{z}\right)
}\text{ \ and \ }\dfrac{\operatorname{Im}m_{\pm}(z)}{\operatorname{Im}z}%
\geq0\text{ \ for any }z\in\mathbb{C}\backslash\mathbb{R}$\\
$m_{+}\left(  x\right)  +m_{-}\left(  x\right)  <0\text{ \ \ for \ }%
x<\lambda_{0}$%
\end{tabular}
\ \text{.}%
\]
A necessary and sufficient condition for $m$ to be a Herglotz function is that
$m$ has a representation%
\[
m\left(  z\right)  =\alpha+\beta z+\int_{-\infty}^{\infty}\left(  \dfrac
{1}{\xi-z}-\dfrac{\xi}{\xi^{2}+1}\right)  \sigma\left(  d\xi\right)
\]
with a real $\alpha$, non-negative $\beta$ and measure $\sigma$ on
$\mathbb{R}$ satisfying%
\[
\int_{-\infty}^{\infty}\dfrac{1}{\xi^{2}+1}\sigma\left(  d\xi\right)
<\infty\text{.}%
\]
The present $m_{\pm}$ are represented as%
\[
m_{\pm}\left(  z\right)  =\alpha_{\pm}+\beta_{\pm}z+\int_{\lambda_{0}}%
^{\infty}\left(  \dfrac{1}{\xi-z}-\dfrac{\xi}{\xi^{2}+1}\right)  \sigma_{\pm
}\left(  d\xi\right)  \text{.}%
\]
The original $m$ is recovered by $m_{\pm}$ by%
\[
m\left(  z\right)  =\left\{
\begin{array}
[c]{cc}%
-m_{+}\left(  -z^{2}\right)  & \text{if \ }\operatorname{Re}z>0\\
m_{-}\left(  -z^{2}\right)  & \text{if \ }\operatorname{Re}z<0
\end{array}
\right.  \text{.}%
\]

\begin{lemma}
\label{l25}If $m$ satisfies the property (\ref{47}), then so do $d_{s}%
m$,\ $d_{\zeta}d_{\overline{\zeta}}m$ if $\left\vert s\right\vert >\mu_{0}$,
$\zeta\in\mathbb{C}\backslash\left(  \mathbb{R}\cup i\mathbb{R}\right)  $.
\end{lemma}

\begin{proof}
If $s>0$, $\operatorname{Re}z>0$, then setting $w\equiv-z^{2}\in\mathbb{C}%
_{-}$, $u=-s^{2}<\lambda_{0}$ we have%
\[
d_{s}m\left(  z\right)  =\dfrac{z^{2}-s^{2}}{-m_{+}\left(  -z^{2}\right)
+m_{+}\left(  -s^{2}\right)  }+m_{+}\left(  -s^{2}\right)  =\dfrac{w-u}%
{m_{+}\left(  w\right)  -m_{+}\left(  u\right)  }+m_{+}\left(  u\right)
\text{.}%
\]
Since $m_{+}$ is of Herglotz,%
\[
\dfrac{m_{+}\left(  w\right)  -m_{+}\left(  u\right)  }{w-u}=\beta_{+}%
+\int_{\lambda_{0}}^{\infty}\dfrac{1}{\left(  \lambda-w\right)  \left(
\lambda-u\right)  }\sigma_{+}\left(  d\lambda\right)
\]
and we see $\operatorname{Im}\left(  m_{+}\left(  w\right)  -m_{+}\left(
u\right)  \right)  /\left(  w-u\right)  >0$ due to $\operatorname{Im}w<0$,
$u<\lambda_{0}$, which implies $\operatorname{Im}d_{s}m\left(  z\right)  >0$.
If $s<0$, $\operatorname{Re}z>0$, then%
\[
d_{s}m_{\boldsymbol{a}}\left(  z\right)  =\dfrac{z^{2}-s^{2}}{-m_{+}\left(
-z^{2}\right)  -m_{-}\left(  -s^{2}\right)  }-m_{-}\left(  -s^{2}\right)
=\dfrac{w-u}{m_{+}\left(  w\right)  +m_{-}\left(  u\right)  }-m_{-}\left(
u\right)  \text{.}%
\]
Since%
\[
\dfrac{m_{+}\left(  w\right)  +m_{-}\left(  u\right)  }{w-u}=\dfrac
{m_{+}\left(  w\right)  -m_{+}\left(  u\right)  }{w-u}+\dfrac{m_{+}\left(
u\right)  +m_{-}\left(  u\right)  }{w-u}\in\mathbb{C}_{-}%
\]
due to $m_{+}\left(  u\right)  +m_{-}\left(  u\right)  <0$, we have
$d_{s}m_{\boldsymbol{a}}\left(  z\right)  \in\mathbb{C}_{+}$. The cases
($s<0$, $\operatorname{Re}z<0$), ($s>0$, $\operatorname{Re}z<0$) can be
treated similarly.

On the other hand note%
\[
d_{\zeta}d_{\overline{\zeta}}m\left(  z\right)  =\dfrac{z^{2}-\overline{\zeta
}^{2}}{\dfrac{z^{2}-\zeta^{2}}{m\left(  z\right)  -m\left(  \zeta\right)
}-\dfrac{\zeta^{2}-\overline{\zeta}^{2}}{m\left(  \zeta\right)  -\overline
{m\left(  \zeta\right)  }}}-\left(  \dfrac{\zeta^{2}-\overline{\zeta}^{2}%
}{m\left(  \zeta\right)  -\overline{m\left(  \zeta\right)  }}-m\left(
\zeta\right)  \right)  \text{.}%
\]
We can assume $\operatorname{Im}z$, $\operatorname{Im}\zeta>0$. To compute the
imaginary part the term $\left(  \zeta^{2}-\overline{\zeta}^{2}\right)
/\left(  m\left(  \zeta\right)  -\overline{m\left(  \zeta\right)  }\right)  $
can be neglected, and%
\[
\dfrac{z^{2}-\overline{\zeta}^{2}}{\dfrac{z^{2}-\zeta^{2}}{m\left(  z\right)
-m\left(  \zeta\right)  }-\dfrac{\zeta^{2}-\overline{\zeta}^{2}}{m\left(
\zeta\right)  -\overline{m\left(  \zeta\right)  }}}+m\left(  \zeta\right)
=\dfrac{m\left(  z\right)  w-\overline{m\left(  \zeta\right)  }a}{w-a}%
\]
with%
\[
w=\dfrac{z^{2}-\zeta^{2}}{m\left(  z\right)  -m\left(  \zeta\right)  }\text{,
\ }a=\dfrac{\zeta^{2}-\overline{\zeta}^{2}}{m\left(  \zeta\right)
-\overline{m\left(  \zeta\right)  }}\text{.}%
\]
Hence%
\begin{align*}
\operatorname{Im}d_{\zeta}d_{\overline{\zeta}}m\left(  z\right)   &
=\operatorname{Im}\dfrac{\left(  m\left(  z\right)  w-\overline{m\left(
\zeta\right)  }a\right)  \left(  \overline{w}-a\right)  }{\left\vert
w-a\right\vert ^{2}}\\
&  =\dfrac{\left\vert w\right\vert ^{2}\operatorname{Im}%
m(z)-a\operatorname{Im}w\left(  m(z)-m\left(  \zeta\right)  \right)
-a^{2}\operatorname{Im}m\left(  \zeta\right)  }{\left\vert w-a\right\vert
^{2}}\\
&  =\dfrac{\left\vert w\right\vert ^{2}\operatorname{Im}%
m(z)-a\operatorname{Im}\left(  z^{2}-\zeta^{2}\right)  -a^{2}\operatorname{Im}%
m\left(  \zeta\right)  }{\left\vert w-a\right\vert ^{2}}%
\end{align*}
Suppose $\operatorname{Re}z$, $\operatorname{Re}\zeta>0$, then%
\begin{align*}
&  \left\vert w\right\vert ^{2}\operatorname{Im}m(z)-a\operatorname{Im}\left(
z^{2}-\zeta^{2}\right)  -a^{2}\operatorname{Im}m\left(  \zeta\right) \\
&  =\dfrac{\operatorname{Im}v\operatorname{Im}u}{\operatorname{Im}m_{+}\left(
v\right)  }-\left\vert \dfrac{u-v}{m_{+}(u)-m_{+}\left(  v\right)
}\right\vert ^{2}\operatorname{Im}m_{+}(u)
\end{align*}
with $u=-z^{2}$, $v=-\zeta^{2}\in\mathbb{C}_{-}$. The Herglotz representation
for $m_{+}$ shows%
\[
\left\{
\begin{array}
[c]{l}%
\dfrac{m_{+}(u)-m_{+}\left(  v\right)  }{u-v}=\beta_{+}+%
{\displaystyle\int_{\lambda_{0}}^{\infty}}
\dfrac{1}{\left(  \lambda-u\right)  \left(  \lambda-v\right)  }\sigma
_{+}\left(  d\lambda\right) \\
\operatorname{Im}m_{+}(u)=\beta_{+}+%
{\displaystyle\int_{\lambda_{0}}^{\infty}}
\dfrac{1}{\left\vert \lambda-u\right\vert ^{2}}\sigma_{+}\left(
d\lambda\right)
\end{array}
\right.  \text{,}%
\]
which implies%
\[
\dfrac{\operatorname{Im}v\operatorname{Im}u}{\operatorname{Im}m_{+}\left(
v\right)  }-\left\vert \dfrac{u-v}{m_{+}(u)-m_{+}\left(  v\right)
}\right\vert ^{2}\operatorname{Im}m_{+}(u)>0
\]
if $u$, $v\in\mathbb{C}_{-}$. The case $\operatorname{Re}z>0$,
$\operatorname{Re}\zeta<0$ is computed similarly, that is%
\begin{align*}
&  \left\vert w\right\vert ^{2}\operatorname{Im}m(z)-a\operatorname{Im}\left(
z^{2}-\zeta^{2}\right)  -a^{2}\operatorname{Im}m\left(  \zeta\right) \\
&  =-\left\vert \dfrac{u-v}{m_{+}(u)+m_{-}\left(  v\right)  }\right\vert
^{2}\operatorname{Im}m_{+}(u)-\dfrac{\operatorname{Im}v\operatorname{Im}%
u}{\operatorname{Im}m_{-}\left(  v\right)  }>0\text{,}%
\end{align*}
since $\operatorname{Im}m_{+}(u)<0$, $\operatorname{Im}m_{-}\left(  v\right)
>0$ if $\operatorname{Im}u<0$, $\operatorname{Im}v>0$, which completes the proof.
\end{proof}

\subsection{Conformal maps}

Although Riemann mapping theorem says that every simply connected domain on
$\mathbb{C}$ can be an image of a conformal map on $\mathbb{C}_{+}$, sometimes
a quantitative estimate of it is necessary. In this section we provide a model
of conformal map from $\mathbb{C}\backslash(-\infty,0]$ to $D_{\omega}^{-}$ of
(\ref{33}).

A conformal map $\psi$ on $\mathbb{C}_{+}$ is easily obtained if
$\operatorname{Im}\psi^{\prime}\left(  z\right)  $ has a definite sign on
$\mathbb{C}_{+}$. A simple such example is $\psi_{\infty}\left(  z\right)
=\sqrt{z}$, and a more general conformal map in this framework can be
constructed for an integer $k\geq1$ by an integral%
\[
\psi_{k}\left(  z\right)  =\frac{k-1/2}{k+1/2}\sqrt{z}+\int_{0}^{\infty}%
\sqrt{z+t}\left(  1+t\right)  ^{-k-3/2}dt\text{.}%
\]
This $\psi_{k}$ satisfies%
\[
\operatorname{Re}\psi_{k}\left(  z\right)  >0\text{, }\operatorname{Im}%
\psi_{k}\left(  z\right)  >0\text{, }\operatorname{Re}\psi_{k}^{\prime}\left(
z\right)  >0\text{, }\operatorname{Im}\psi_{k}^{\prime}\left(  z\right)  <0
\]
for $z\in\mathbb{C}_{+}$, hence $\psi_{k}$ maps $\mathbb{C}_{+}$ to $\psi
_{k}\left(  \mathbb{C}_{+}\right)  \left(  \subset\mathbb{C}_{+}\right)  $,
and $\phi_{k}\left(  z\right)  =\psi_{k}\left(  z\right)  ^{2}$ maps
$\mathbb{C}_{+}$ to $\phi_{k}\left(  \mathbb{C}_{+}\right)  \left(
\subset\mathbb{C}_{+}\right)  $ conformally. Since $\psi_{k}\left(  z\right)
$ takes real values on $[0,\infty)$, $\psi_{k}\left(  z\right)  $ and
$\phi_{k}\left(  z\right)  $ can be extended as conformal maps from
$\mathbb{C}\backslash(-\infty,0]$ to a domain in $\left\{  \operatorname{Re}%
z>0\right\}  $ and a domain in $\mathbb{C}$ respectively. Set%
\[
\left\{
\begin{array}
[c]{l}%
a_{k}=2\int_{0}^{1}s^{2}\left(  1-s^{2}\right)  ^{k-1}ds=\dfrac{\sqrt{\pi
}\Gamma\left(  k\right)  }{2\Gamma\left(  k+3/2\right)  }\\
b_{k}=2a_{k}\left(  2a_{k}^{2}k\left(  k+1/2\right)  ^{2}/\left(
k-1/2\right)  ^{2}+1\right)  ^{-1/2}%
\end{array}
\right.  \text{.}%
\]

\begin{lemma}
\label{l51}The image $\phi_{k}\left(  \mathbb{C}\backslash(-\infty,0]\right)
$ is described as follows:
\[
\phi_{k}\left(  \mathbb{C}\backslash(-\infty,0]\right)  =\mathbb{C}%
\backslash\left\{  z\in\mathbb{C}\text{; }\left\vert \operatorname{Im}%
z\right\vert \leq\omega\left(  \operatorname{Re}z\right)  \text{,
}\operatorname{Re}z\leq a_{k}^{2}\right\}
\]
with positive smooth function $\omega\left(  x\right)  $ on $(-\infty
,a_{k}^{2})$ such that%
\begin{equation}
\omega\left(  x\right)  =\left\{
\begin{array}
[c]{l}%
2a_{k}\left(  -x\right)  ^{-k+1/2}\left(  1+O\left(  x^{-1}\right)  \right)
\text{ \ \ \ \ \ \ \ \ \ as \ }x\rightarrow-\infty\\
b_{k}\left(  a_{k}^{2}-x\right)  ^{1/2}\left(  1+O\left(  a_{k}^{2}-x\right)
\right)  \text{ \ \ \ \ \ as \ }x\rightarrow a_{k}^{2}-0
\end{array}
\right.  \text{.} \label{101}%
\end{equation}
Moreover, $\phi_{k}$ takes a form of%
\begin{equation}
\phi_{k}\left(  z\right)  =z+f_{1}(z)+z^{-k+1/2}f_{2}(z) \label{102}%
\end{equation}
with some real rational functions $f_{1}$, $f_{2}$ (that is, $f_{j}\left(
z\right)  =\overline{f_{j}\left(  \overline{z}\right)  }$ for\ $j=1$, $2$)
satisfying%
\[
\left\{
\begin{array}
[c]{l}%
f_{1}(\infty)=\left(  k^{2}-1/4\right)  ^{-1}\\
f_{2}(\infty)=2\left(  -1\right)  ^{k}a_{k}%
\end{array}
\right.  \text{.}%
\]
Conversely, $\phi_{k}^{-1}\left(  w\right)  $ has an expression%
\begin{equation}
\phi_{k}^{-1}\left(  w\right)  =w+g_{1}\left(  w\right)  +w^{-k+1/2}g_{2}(w)
\label{164}%
\end{equation}
with real $g_{1}$, $g_{2}$ analytic in a neighborhood of $\infty$. Moreover,
it holds that%
\[
g_{1}\left(  \infty\right)  =-\left(  k^{2}-1/4\right)  ^{-1}\text{,
\ \ \ }g_{2}\left(  \infty\right)  =-2\left(  -1\right)  ^{k}a_{k}.
\]

\end{lemma}

\begin{proof}
Setting $s=\sqrt{\left(  z+t\right)  /\left(  1+t\right)  }$, we have%
\[
\psi_{k}\left(  z\right)  =\frac{k-1/2}{k+1/2}\sqrt{z}+2\left(  z-1\right)
^{-k}\int_{1}^{\sqrt{z}}s^{2}\left(  s^{2}-1\right)  ^{k-1}ds\text{.}%
\]
Since the integral $\int_{0}^{z}s^{2}\left(  s^{2}-1\right)  ^{k-1}ds$ is an
odd polynomial of degree $2k+1$, the integral
\[
p\left(  z\right)  =\frac{k-1/2}{k+1/2}\left(  z-1\right)  ^{k}+2\sqrt{z}%
^{-1}\int_{0}^{\sqrt{z}}s^{2}\left(  s^{2}-1\right)  ^{k-1}ds
\]
defines a polynomial of degree $k$, and%
\begin{equation}
\psi_{k}\left(  z\right)  =\left(  z-1\right)  ^{-k}\left(  \sqrt{z}p\left(
z\right)  -p\left(  1\right)  \right)  \label{103}%
\end{equation}
holds. It should be noted that $\sqrt{z}p\left(  z\right)  -p\left(  1\right)
$ has zero of degree $k$ at $z=1$, so $\psi_{k}\left(  z\right)  $ has no
singularity at $z=1$. Set%
\[
s(x)=\operatorname{Re}\psi_{k}\left(  x+i0\right)  \text{, \ \ }t\left(
x\right)  =\operatorname{Im}\psi_{k}\left(  x+i0\right)
\]
for $x\in\mathbb{R}$. Then, (\ref{103}) implies%
\begin{align*}
s(x)  &  =\left\{
\begin{array}
[c]{l}%
-p(1)\left(  x-1\right)  ^{-k}\text{ \ \ \ \ \ \ \ \ \ \ \ \ \ \ \ \ for
}x<0\\
\left(  x-1\right)  ^{-k}\left(  \sqrt{x}p\left(  x\right)  -p\left(
1\right)  \right)  \text{ \ for }x\geq0
\end{array}
\right.  \text{,}\\
t(x)  &  =\left\{
\begin{array}
[c]{l}%
\left(  x-1\right)  ^{-k}\sqrt{-x}p\left(  x\right)  \text{ \ \ \ \ for }x<0\\
0\text{ \ \ \ \ \ \ \ \ \ \ \ \ \ \ \ \ \ \ \ \ \ \ \ \ \ \ for }x\geq0
\end{array}
\right.  \text{,}%
\end{align*}
and their asymptotics are%
\begin{align*}
s(x)  &  =\left\{
\begin{array}
[c]{l}%
a_{k}\left(  1+kx+O\left(  x^{2}\right)  \right)  \text{
\ \ \ \ \ \ \ \ \ \ \ \ \ \ \ as \ }x\rightarrow-0\\
a_{k}\left(  -x\right)  ^{-k}\left(  1+kx^{-1}+O\left(  x^{-2}\right)
\right)  \text{ \ \ as \ }x\rightarrow-\infty
\end{array}
\right. \\
t\left(  x\right)   &  =\left\{
\begin{array}
[c]{l}%
\dfrac{k-1/2}{k+1/2}\sqrt{-x}\left(  1+O\left(  x\right)  \right)  \text{
\ \ as \ }x\rightarrow-0\\
\sqrt{-x}\left(  1+O\left(  x^{-1}\right)  \right)  \text{
\ \ \ \ \ \ \ \ \ \ \ as \ }x\rightarrow-\infty
\end{array}
\right.  \text{,}%
\end{align*}
where we have used%
\[
\left\{
\begin{array}
[c]{l}%
p\left(  1\right)  =\left(  -1\right)  ^{k-1}a_{k}\text{, \ \ \ \ }%
p(0)=\dfrac{k-1/2}{k+1/2}\left(  -1\right)  ^{k}\\
p(z)=z^{k}-\left(  k-\dfrac{2}{4k^{2}-1}\right)  z^{k-1}+\cdots
\end{array}
\right.  \text{.}%
\]
From (\ref{103})
\[
\phi_{k}\left(  z\right)  =\psi_{k}\left(  z\right)  ^{2}=z+f_{1}%
(z)+z^{-k+1/2}f_{2}(z)
\]
follows, which yields (\ref{102}) with%
\[
\left\{
\begin{array}
[c]{l}%
f_{1}(z)=\left(  z-1\right)  ^{-2k}\left(  p\left(  1\right)  ^{2}%
+zp(z)^{2}\right)  -z\\
f_{2}(z)=-2\left(  z-1\right)  ^{-2k}p\left(  1\right)  p\left(  z\right)
z^{k}%
\end{array}
\right.  \text{,}%
\]
and
\begin{equation}
\left\{
\begin{array}
[c]{l}%
\operatorname{Re}\phi_{k}\left(  x+i0\right)  =x+f_{1}(x)+f_{2}(x)\times
\left\{
\begin{array}
[c]{l}%
x^{-k+1/2}\text{ \ \ \ \ if \ }x>0\\
0\text{ \ \ \ \ \ \ \ \ \ \ \ if \ }x<0
\end{array}
\right. \\
\operatorname{Im}\phi_{k}\left(  x+i0\right)  =\left\{
\begin{array}
[c]{l}%
0\text{ \ \ \ \ \ \ \ \ \ \ \ \ \ \ \ \ \ \ \ \ \ \ \ \ \ \ \ \ \ if \ }x>0\\
\sqrt{-x}x^{-k}f_{2}(x)\text{\ \ \ if \ }x<0
\end{array}
\right.
\end{array}
\right.  \label{104}%
\end{equation}
is valid, hence (\ref{57}) shows%
\[
\left\{
\begin{array}
[c]{l}%
\operatorname{Re}\phi_{k}\left(  x+i0\right)  =x+\left(  k^{2}-1/4\right)
^{-1}+O\left(  \left(  -x\right)  ^{-1}\right) \\
\operatorname{Im}\phi_{k}\left(  x+i0\right)  =2a_{k}\left(  -x\right)
^{-k+1/2}+O\left(  \left(  -x\right)  ^{-k-1/2}\right)
\end{array}
\right.
\]
as $x\rightarrow-\infty$, and%
\[
\left\{
\begin{array}
[c]{l}%
\operatorname{Re}\phi_{k}\left(  x+i0\right)  =a_{k}^{2}+\left(  2ka_{k}%
^{2}+\left(  \dfrac{k-1/2}{k+1/2}\right)  ^{2}\right)  x+O\left(  x^{2}\right)
\\
\operatorname{Im}\phi_{k}\left(  x+i0\right)  =2\dfrac{k-1/2}{k+1/2}a_{k}%
\sqrt{-x}+O\left(  \left(  -x\right)  ^{3/2}\right)
\end{array}
\right.
\]
as $x\rightarrow-0$.\ Since $\operatorname{Re}\phi_{k}\left(  x+i0\right)
=s(x)^{2}-t\left(  x\right)  ^{2}$, $\operatorname{Im}\phi_{k}\left(
x+i0\right)  =2s\left(  x\right)  t\left(  x\right)  $, (\ref{52}) implies
\[
\left\{
\begin{array}
[c]{l}%
\operatorname{Re}\phi_{k}\left(  x+i0\right)  \text{ is increasing and moving
from}-\infty\text{ to }\infty\\
\operatorname{Im}\phi_{k}\left(  x+i0\right)  \text{ }>0\text{ on }\left(
-\infty,0\right)  \text{ and }0\text{ on }[0,\infty)
\end{array}
\right.  \text{.}%
\]
Therefore, $\omega$ can be defined by an equation%
\[
\omega\left(  \operatorname{Re}\phi_{k}\left(  x+i0\right)  \right)
=\operatorname{Im}\phi_{k}\left(  x+i0\right)  \text{.}%
\]
due to (\ref{56}), and (\ref{104}), (\ref{55}), (\ref{56}) show $\omega\left(
x\right)  $ satisfies (\ref{101}).

We use (\ref{103}) to show (\ref{164}). Set $\vartheta\left(  z\right)
=z^{2}$. $\vartheta$ is a conformal map from $\left\{  \operatorname{Re}%
z>0\right\}  $ to $\mathbb{C}\backslash(-\infty,0]$ and define $\widetilde
{\psi}_{k}\left(  s\right)  =\psi_{k}\left(  \vartheta\left(  s\right)
\right)  $. Then the function%
\begin{align*}
F(s)  &  =\widetilde{\psi}_{k}\left(  s\right)  -s\\
&  =-p\left(  1\right)  \left(  s^{2}-1\right)  ^{-k}+s\left(  \left(
s^{2}-1\right)  ^{-k}p\left(  s^{2}\right)  -1\right)  \text{.}%
\end{align*}
is a rational function whose poles only at $s=\pm1$ and has expansion%
\[
F(s)=c_{1}s^{-1}+c_{2}s^{-3}+\cdots+c_{k}s^{-2k+1}+c_{k+1}s^{-2k}+
\]
at $s=\infty$ with $c_{1}=\left(  2k^{2}-1/2\right)  ^{-1}$ and $c_{k+1}%
=-p(1)$, namely the first coefficient of even order starts from $2k$. We
consider an equation for a given $t$:%
\begin{equation}
s+F(s)=t \label{170}%
\end{equation}
and find a solution of a form%
\[
s=t+G(t)\text{.}%
\]
Since the even coefficients of the power series of $F$ vanish up to $2\left(
k-1\right)  $, Lemma \ref{l52} shows that there exists uniquely such $G$ that
$G$ is real and analytic near $t=\infty$ and the odd coefficients of $G$
vanishes up to $2\left(  k-1\right)  $. $\widetilde{\psi}_{k}\left(  z\right)
$ is one-to-one on $\left\{  \left\vert z\right\vert >r_{1}\right\}  $ and its
inverse is given by $w+G(w)$ on $\left\{  \left\vert w\right\vert
>r_{2}\right\}  $. Since $\phi_{k}\left(  z\right)  =$ $\vartheta\left(
\psi_{k}\left(  z\right)  \right)  $ is a conformal map from $\mathbb{C}%
\backslash(-\infty,0]$ to $\phi_{k}\left(  \mathbb{C}\backslash(-\infty
,0]\right)  $, its inverse is given by $\phi_{k}^{-1}\left(  w\right)
=\left(  \vartheta\widetilde{\psi}^{-1}\vartheta^{-1}\right)  \left(
w\right)  $ for $w\in\phi_{m}\left(  \mathbb{C}\backslash(-\infty,0]\right)
$. Let
\[
\left\{
\begin{array}
[c]{l}%
G_{1}\left(  t\right)  =\dfrac{1}{2}\left(  G\left(  \sqrt{t}\right)
+G\left(  -\sqrt{t}\right)  \right)  \text{ \ }\left(  =G_{e}\left(  \sqrt
{t}\right)  \right) \\
G_{2}\left(  t\right)  =\dfrac{1}{2\sqrt{t}}\left(  G\left(  \sqrt{t}\right)
-G\left(  -\sqrt{t}\right)  \right)  \text{ \ }\left(  =\dfrac{1}{\sqrt{t}%
}G_{o}\left(  \sqrt{t}\right)  \right)
\end{array}
\right.  \text{.}%
\]
Then $G\left(  t\right)  =G_{1}\left(  t^{2}\right)  +tG_{2}(t^{2})$, and we
have%
\begin{align*}
&  \left(  \vartheta\widetilde{\psi}^{-1}\vartheta^{-1}\right)  \left(
w\right) \\
&  =\left(  \sqrt{w}+G_{1}\left(  w\right)  +\sqrt{w}G_{2}(w)\right)  ^{2}\\
&  =w+G_{1}\left(  w\right)  ^{2}+w\left(  \left(  G_{2}(w)+1\right)
^{2}-1\right)  +2\sqrt{w}G_{1}\left(  w\right)  \left(  G_{2}(w)+1\right)
\text{.}%
\end{align*}
Since $G$ has the even coefficients vanishing up to $2\left(  k-1\right)  $,
$w^{k}G_{1}\left(  w\right)  $ is analytic near $w=\infty$. Therefore,
setting
\[
\left\{
\begin{array}
[c]{l}%
g_{1}(w)=G_{1}\left(  w\right)  ^{2}+wG_{2}(w)\left(  G_{2}(w)+2\right) \\
g_{2}(w)=2w^{k}G_{1}\left(  w\right)  \left(  G_{2}(w)+1\right)
\end{array}
\right.  \text{,}%
\]
we have
\[
\phi_{k}^{-1}\left(  w\right)  =w+g_{1}(w)+w^{-k+1/2}g_{2}(w)
\]
with some $g_{1}$, $g_{2}$ analytic in a neighborhood of $\infty$ satisfying%
\[
g_{1}(\infty)=-\left(  k^{2}-1/4\right)  ^{-1}\text{, \ }g_{2}\left(
\infty\right)  =-2\left(  -1\right)  ^{k}a_{k}\text{,}%
\]
\ which completes the proof.
\end{proof}

\begin{lemma}
\label{l52}Let $F$ be a power series of $s^{-1}$ given by $F(s)=\sum
_{j=1}^{\infty}a_{j}s^{-j}$ and assume it has the positive radius of
convergence and consider an equation:%
\begin{equation}
t=s+F(s)\text{.} \label{166}%
\end{equation}
(i) This equation is uniquely solvable if $\left\vert t^{-1}\right\vert $ is
sufficiently small and it has a form:%
\[
s=t+G(t)
\]
with a convergent power series of $t^{-1}$ given by%
\begin{equation}
G(t)=\sum_{j=1}^{\infty}x_{j}t^{-j}\text{.} \label{169}%
\end{equation}
(ii) $x_{n}$ is determined from $\left\{  a_{j}\right\}  _{j=1}^{n}$ for each
$n\geq1$. The first three coefficients are%
\[
x_{1}=-a_{1}\text{, \ \ }x_{2}=-a_{2}\text{, \ \ }x_{3}=-a_{1}^{2}%
-a_{3}\text{.}%
\]
\newline(iii) Suppose $F(s)$ has a form
\begin{equation}
F(s)=\sum_{j=1}^{k}a_{2j-1}s^{-2j+1}+\sum_{j=2k}^{\infty}a_{j}s^{-j}
\label{167}%
\end{equation}
for an $k\geq1$. Then, the coefficients $x_{j}$ of $G(t)$ vanish for even $j $
up to $2\left(  k-1\right)  $. Moreover, if $a_{2j}\neq0$, then $x_{2j}%
=-a_{2j}$.
\end{lemma}

\begin{proof}
Replacing $s$ by $s^{-1}$ and $t$ by $t^{-1}$ we see the equation (\ref{166})
is equivalent to%
\begin{equation}
t=\dfrac{s}{1+sF(s^{-1})}\text{.} \label{168}%
\end{equation}
The condition on $F$ implies%
\[
t(0)=0\text{, }\dfrac{dt}{ds}\left(  0\right)  =1\text{, }\dfrac{d^{2}%
t}{ds^{2}}\left(  0\right)  =0\text{,}%
\]
hence the complex function theory shows the existence of the solution $s(t)$
of (\ref{168}) in a neighborhood of $0$ satisfying%
\[
s(0)=0\text{, }\dfrac{ds}{dt}\left(  0\right)  =1\text{, }\dfrac{d^{2}%
s}{dt^{2}}\left(  0\right)  =0\text{,}%
\]
which implies the existence of $G$ of the form of (\ref{169}). One can show
inductively that the coefficient $x_{n}$ is determined from $\left\{
a_{j}\right\}  _{j=1}^{n}$. To show (iii) one can assume $a_{j}=0$ for every
$j\geq2k$ owing to (ii). The relation between $F$, $G$ is rewritten as%
\[
F(s)+G\left(  F(s)+s\right)  =0\text{.}%
\]
If we define $\widehat{f}(s)=-f\left(  -s\right)  $, then the above equation
turns to%
\[
\widehat{F}\left(  s\right)  +\widehat{G}\left(  \widehat{F}(s)+s\right)
=0\text{.}%
\]
Therefore, if $\widehat{F}\left(  s\right)  =F(s)$, the uniqueness implies
$\widehat{G}\left(  s\right)  =G(s)$, which shows the first part of (iii). To
show the second part we note that if%
\[
\left\{
\begin{array}
[c]{l}%
F(s)=\sum_{j=1}^{k-1}a_{j}s^{-j}+a_{k}s^{-k}\equiv F_{1}(s)+a_{k}s^{-k}\\
G(s)=\sum_{j=1}^{k-1}x_{j}s^{-j}+x_{k}s^{-k}+\sum_{j=k+1}^{\infty}x_{j}%
s^{-j}\\
\text{ \ \ \ \ \ \ }\equiv G_{\left(  1\right)  }(s)+x_{k}s^{-k}+G_{\left(
2\right)  }(s)
\end{array}
\right.  \text{,}%
\]
and with some $b_{m}$
\[
F_{1}(s)+G_{\left(  1\right)  }\left(  s+F_{1}(s)\right)  =b_{k}%
s^{-k}+O\left(  s^{-k-1}\right)
\]
holds, which is verified by induction, then the identity%
\begin{gather*}
F_{1}(s)+a_{k}s^{-k}+G_{\left(  1\right)  }(s+F_{1}(s)+a_{k}s^{-k}%
)+x_{k}\left(  s+F_{1}(s)+a_{k}s^{-k}\right)  ^{-k}\\
+G_{\left(  2\right)  }(s+F_{1}(s)+a_{k}s^{-k})=0
\end{gather*}
together with%
\[
G_{\left(  1\right)  }(s+F_{1}(s)+a_{k}s^{-k})=G_{\left(  1\right)  }%
(s+F_{1}(s))+O\left(  s^{-k-2}\right)
\]
implies $x_{k}=-a_{k}-b_{k}$. Since, if $k$ is even and $a_{2j}=0$ for any
$j\leq k/2$, then (ii) implies $x_{k}=0$, and hence $b_{k}=0$. However,
clearly $b_{k}$ is determined from $\left\{  a_{j}\right\}  _{1\leq j\leq
k-1}$, hence $b_{k}=0$ is valid if $a_{2j}=0$ for any $j\leq k/2-1$ regardless
of the value $a_{k}$. Consequently we have $x_{k}=-a_{k}$ if $k$ is even and
$a_{2j}=0$ for any $j\leq k/2-1$ holds.
\end{proof}

\subsection{Estimates of relevant integral}

Suppose the curve $C$ is of the form:%
\[
C=\left\{  \pm\omega\left(  y\right)  +iy\text{; \ }y\in\mathbb{R}\text{,
\ }\omega\left(  y\right)  =O\left(  y^{-\left(  n-1\right)  }\right)
\right\}
\]
with $\omega\left(  y\right)  >0$, $\omega\left(  y\right)  =\omega\left(
-y\right)  $. Assume (\ref{32}), namely%
\[
\sup_{z\in C}%
{\displaystyle\int_{\left\vert z-\lambda\right\vert \leq1,\lambda\in C}}
\left\vert d\lambda\right\vert <\infty\text{.}%
\]

\begin{lemma}
\label{l15}Let $C^{\prime}=\sigma C$ with $\sigma>1$. Then%
\[%
{\displaystyle\int_{C}}
\frac{1}{\left\vert \lambda-z\right\vert ^{2}}\left\vert d\lambda\right\vert
=O\left(  \left\vert z\right\vert ^{n-1}\right)  \text{ \ for }z\in C^{\prime
}\text{.}%
\]

\end{lemma}

\begin{proof}
Let%
\begin{align*}%
{\displaystyle\int_{C}}
\frac{1}{\left\vert \lambda-z\right\vert ^{2}}\left\vert d\lambda\right\vert
&  =\int_{\left\vert \operatorname{Im}\left(  \lambda-z\right)  \right\vert
\leq\delta\left\vert \operatorname{Im}z\right\vert }\dfrac{\left\vert
d\lambda\right\vert }{\left\vert \lambda-z\right\vert ^{2}}+\int_{\left\vert
\operatorname{Im}\left(  z-\lambda\right)  \right\vert >\delta\left\vert
\operatorname{Im}z\right\vert }\dfrac{\left\vert d\lambda\right\vert
}{\left\vert \lambda-z\right\vert ^{2}}\\
&  \equiv I_{1}+I_{2}%
\end{align*}
with $\delta\in\left(  0,1\right)  $ specified later. Since $C$ is
parametrized as $\omega\left(  t\right)  +it$, we see%
\[
I_{1}=\int_{\left\vert t-\operatorname{Im}z\right\vert \leq\delta\left\vert
\operatorname{Im}z\right\vert }\dfrac{\left(  1+\omega^{\prime}(t)^{2}\right)
^{1/2}}{\left(  t-\operatorname{Im}z\right)  ^{2}+\left(  \omega
(t)-\operatorname{Re}z\right)  ^{2}}dt\leq c_{1}\pi\rho\left(  z\right)
^{-1}\text{,}%
\]
where%
\[
\left\{
\begin{array}
[c]{l}%
\rho\left(  z\right)  =\inf_{\lambda\in C;\left\vert \operatorname{Im}\left(
\lambda-z\right)  \right\vert \leq\delta\left\vert \operatorname{Im}%
z\right\vert }\left\vert \operatorname{Re}\left(  z-\lambda\right)
\right\vert \\
c_{1}=\sup_{t\in\mathbb{R}}\left(  1+\omega^{\prime}(t)^{2}\right)  ^{1/2}%
\end{array}
\right.  \text{.}%
\]
$I_{2}$ is estimated as%
\begin{align*}
I_{2}  &  =\int_{\left\vert t-\operatorname{Im}z\right\vert >\left\vert
\operatorname{Im}z\right\vert \delta}\dfrac{\left(  1+\omega^{\prime}%
(t)^{2}\right)  ^{1/2}dt}{\left\vert \left(  \omega\left(  t\right)
-\operatorname{Re}z\right)  +i\left(  t-\operatorname{Im}z\right)  \right\vert
^{2}}\\
&  \leq c_{1}\int_{\left\vert t-\operatorname{Im}z\right\vert >\left\vert
\operatorname{Im}z\right\vert \delta}\dfrac{dt}{\left\vert t-\operatorname{Im}%
z\right\vert ^{2}}=c_{1}\int_{\left\vert x\right\vert >\left\vert
\operatorname{Im}z\right\vert \delta}\dfrac{dx}{\left\vert x\right\vert ^{2}%
}=c_{2}\left\vert \operatorname{Im}z\right\vert ^{-1}.
\end{align*}
Therefore, we have%
\[
\int_{C}\dfrac{\left\vert d\lambda\right\vert }{\left\vert \lambda
-z\right\vert ^{2}}\leq c_{1}\pi\rho\left(  z\right)  ^{-1}+c_{2}\left\vert
\operatorname{Im}z\right\vert ^{-1}\text{.}%
\]
We have to show%
\begin{equation}
\rho\left(  z\right)  \geq c_{3}\left\vert \operatorname{Im}z\right\vert
^{-\left(  n-1\right)  }\text{,} \label{95}%
\end{equation}
if $\delta$ is chosen suitably. Assume $\operatorname{Re}z$,
$\operatorname{Im}z>0$. Note that in the region%
\[
\left\{  \lambda\in C;\left\vert \operatorname{Im}\left(  z-\lambda\right)
\right\vert \leq\delta\operatorname{Im}z\right\}  =\left\{  \lambda\in
C;\left(  1-\delta\right)  \operatorname{Im}z\leq\operatorname{Im}\lambda
\leq\left(  1+\delta\right)  \operatorname{Im}z\right\}  \text{.}%
\]
Since $\omega\left(  t\right)  =t^{-(n-1)}$ for sufficiently large $t$,
inequalities%
\begin{align*}
\left\vert \operatorname{Re}\left(  z-\lambda\right)  \right\vert  &
\geq\operatorname{Re}z-\operatorname{Re}\lambda\geq\sigma\omega\left(
\sigma^{-1}\left\vert \operatorname{Im}z\right\vert \right)  -\omega\left(
\left(  1-\delta\right)  \left\vert \operatorname{Im}z\right\vert \right) \\
&  =\left(  \sigma^{n}-\left(  1-\delta\right)  ^{-\left(  n-1\right)
}\right)  \left\vert \operatorname{Im}z\right\vert ^{-\left(  n-1\right)  }%
\end{align*}
are valid for $z\in C^{\prime}$. Therefore, if%
\[
\sigma^{n}>\left(  1-\delta\right)  ^{-\left(  n-1\right)  }\text{,}%
\]
we have (\ref{95}).
\end{proof}

\subsection{Ergodic Schr\"{o}dinger operators}

This section provides several basic facts on 1D Schr\"{o}dinger operators with
ergodic potentials, which are necessary in this paper.

Let $\left(  \Omega,\mathcal{F},P\right)  $ be a probability space and
$\left\{  \theta_{x}\right\}  _{x\in\mathbb{R}}$ be a one parameter group of
$\mathcal{F}$-measurable transformations on $\Omega$ which satisfies%
\begin{equation}
P\left(  \theta_{x}^{-1}A\right)  =P\left(  A\right)  \text{ for any }%
x\in\mathbb{R}\text{ and }A\in\mathcal{F}\text{. (stationarity)} \label{37}%
\end{equation}
$\left(  \Omega,\mathcal{F},P,\left\{  \theta_{x}\right\}  _{x\in\mathbb{R}%
}\right)  $ is called ergodic if it satisfies%
\begin{equation}
P\left(  \theta_{x}^{-1}A\circleddash A\right)  =0\text{ for any }%
x\in\mathbb{R}\Longrightarrow P\left(  A\right)  =0\text{ \ or \ }1\text{.}
\label{59}%
\end{equation}
For an $\mathcal{F}$-measurable real valued function $Q$ on $\Omega$ set%
\[
q_{\omega}\left(  x\right)  =Q\left(  \theta_{x}\omega\right)  \text{,
\ \ }\omega\in\Omega\text{.}%
\]
Then we obtain an ergodic potential $\left\{  q_{\omega}\right\}  _{\omega
\in\Omega}$. A simple but important example is quasi-periodic potentials. Set
$\Omega=\mathbb{R}^{n}/\mathbb{Z}^{n}$ and for $\alpha\in\mathbb{R}^{n}$%
\[
\theta_{x}\omega=x\alpha+\omega\text{, \ }P=\text{the Lebesgue measure on
}\mathbb{R}^{n}/\mathbb{Z}^{n}\text{.}%
\]
This $\left(  \Omega,\mathcal{F},P,\left\{  \theta_{x}\right\}  _{x\in
\mathbb{R}}\right)  $ is ergodic if $\alpha$ is rationally independent and the
resulting $q_{\omega}\left(  x\right)  $ is a quasi-periodic function. If
$n=1$, we have a periodic function and for $n=\infty$ in a certain sense we
have an almost periodic function. One has more random ergodic potentials. For
a technical reason we assume%
\begin{equation}
\mathbb{E}\left(  \left\vert Q\right\vert \right)  =\int_{\Omega}\left\vert
Q\left(  \omega\right)  \right\vert P\left(  d\omega\right)  <\infty\text{
\ and }Q\left(  \omega\right)  \geq\lambda_{0}\text{ \ for any }\omega
\in\Omega\text{.} \label{87}%
\end{equation}
$\mathbb{E}$ denotes the expectation by $P$. Then one can consider the
associated Schr\"{o}dinger operator%
\[
L_{\omega}=-\partial_{x}^{2}+q_{\omega}\text{.}%
\]
Under the condition (\ref{87}) it is known that $\inf$ \textrm{sp }$L_{\omega
}\geq\lambda_{0}$ and the boundaries $\pm\infty$ are of limit point type for
$L_{\omega}$ for a.e. $\omega\in\Omega$. One can apply the Weyl spectral
theory to each $L_{\omega}$.

The \textbf{Floquet exponent} is defined by%
\begin{equation}
w(z)=\mathbb{E}\left(  m_{\pm}\left(  z,\omega\right)  \right)  \text{ (the
two expectations coincide),} \label{115}%
\end{equation}
by which the \textbf{Lyapounov exponent} and \textbf{integrated density of
states} are obtained by%
\[
\gamma\left(  z\right)  =-\operatorname{Re}w(z)\text{ }\left(  \geq0\right)
\text{, \ }N\left(  \lambda\right)  =\dfrac{1}{\pi}\operatorname{Im}%
w(\lambda)\text{ }\left(  \lambda\in\mathbb{R}\right)  \text{.}%
\]
$N\left(  \lambda\right)  $ is non-negative, continuous and non-decreasing on
$\mathbb{R}$. \cite{k1} found an identity%
\[
\dfrac{\gamma\left(  z\right)  }{\operatorname{Im}z}-\operatorname{Im}%
w^{\prime}\left(  z\right)  =\dfrac{1}{4}\mathbb{E}\left(  \left(  \dfrac
{1}{\operatorname{Im}m_{+}\left(  z,\omega\right)  }+\dfrac{1}%
{\operatorname{Im}m_{-}\left(  z,\omega\right)  }\right)  \left\vert R\left(
z,\omega\right)  \right\vert ^{2}\right)
\]
for $z\in\mathbb{C}_{+}$. Set%
\[
\chi\left(  z\right)  =\dfrac{\gamma\left(  z\right)  }{\operatorname{Im}%
z}-\operatorname{Im}w^{\prime}\left(  z\right)  \geq0\text{.}%
\]
Then applying the Schwarz inequality we have%
\begin{align}
\mathbb{E}\left(  \left\vert R\left(  z,\omega\right)  \right\vert \right)
&  \leq\sqrt{4\chi\left(  z\right)  }\sqrt{\mathbb{E}\left(  \dfrac
{1}{\operatorname{Im}m_{+}\left(  z,\omega\right)  }+\dfrac{1}%
{\operatorname{Im}m_{-}\left(  z,\omega\right)  }\right)  ^{-1}}\nonumber\\
&  \leq\sqrt{4\chi\left(  z\right)  }\sqrt{\mathbb{E}\left(  \dfrac
{\operatorname{Im}m_{+}\left(  z,\omega\right)  +\operatorname{Im}m_{-}\left(
z,\omega\right)  }{4}\right)  }\nonumber\\
&  =\sqrt{2\chi\left(  z\right)  \operatorname{Im}w\left(  z\right)  }\text{
\ \ (due to (\ref{115})).} \label{90}%
\end{align}
It is also known that%
\begin{equation}
\Sigma_{ac}^{\omega}=\left\{  \lambda\in\mathbb{R}\text{; \ }\gamma\left(
\lambda\right)  =0\right\}  =\left\{  \lambda\in\mathbb{R}\text{; \ }R\left(
\lambda,\omega\right)  =0\right\}  \text{ for a.e. }\omega\in\Omega\text{.}
\label{114}%
\end{equation}

\bigskip

\bigskip

Acknowledgement: \emph{This research was partly supported by JSPS (grant no.
26400128).}

\end{document}